\documentclass[a4paper,12pt]{amsart}

\usepackage{a4wide}
\usepackage[driverfallback=dvipdfm]{hyperref}
\usepackage[all]{xy}
\usepackage{tikz}
\usepackage{amsaddr}
\usepackage{amsthm, amssymb}
\usepackage{color}
\usepackage{enumerate}

\newtheorem{thm}{Theorem}[section]
\newtheorem*{thm*}{Theorem}
\newtheorem{cor}[thm]{Corollary}
\newtheorem{lem}[thm]{Lemma}
\newtheorem{prop}[thm]{Proposition}

\theoremstyle{definition}
\newtheorem{defn}[thm]{Definition}

\newtheorem{nt}[thm]{Notation}
\newtheorem{rem}[thm]{Remark}
\newtheorem{ex}[thm]{Example}
\newtheorem*{rem*}{Remark}

\newtheorem*{theoremaux}{Theorem \theoremauxnum}
\gdef\theoremauxnum{1}

\def \p{{\mathbb P}}          

\def\N{{\mathbb N}}           
\def\R{{\mathbb R}}           
\def\O{{\mathcal O}}          
\def\Og{{\mathcal O}_{\text{gr}}}          
\def\G{{\mathcal G}}          
\def\e{\varepsilon}           
\def\f{\varphi}           

\def\Aut{\operatorname{Aut}}
\def\Out{\operatorname{Out}}
\def\Lip{\operatorname{Lip}}
\def\Hor{\operatorname{Hor}}
\def\PL{\operatorname{Str}}
\def\opt{\operatorname{opt}}
\def\wopt{\operatorname{weakopt}}
\def\vol{\operatorname{vol}}
\def\Min{\operatorname{Min}}
\def\LocMin{\operatorname{LocMin}}
\def\TT{\operatorname{TT}}
\def\TTo{\operatorname{TT_0}}
\def\rank{\operatorname{rank}}
\def\core{\operatorname{core}}

\def\simfk{\langle \sim_{f^k}\rangle}
\def\wt{\widetilde}
\def\ul{\underline}

\title[Displacement of automorphisms I]{Displacements of automorphisms of free groups
I: Displacement functions, minpoints and train tracks}
\author{Stefano Francaviglia}
\address{Dipartimento di Matematica of the University of
Bologna}
\email{stefano.francaviglia@unibo.it}
\author{Armando Martino}
\address{Mathematical Sciences, University of Southampton }
\email{A.Martino@soton.ac.uk}
\begin{document}


\begin{abstract}
This is the first of two papers in which we investigate the properties of the displacement
functions of automorphisms of free groups (more generally, free products) on Culler-Vogtmann
Outer space  and its simplicial bordification - the free splitting complex - with respect to the Lipschitz metric. The theory for irreducible automorphisms being well-developed, we
concentrate on the reducible case. Since we deal with the bordification, we develop
all the needed tools in the more general setting of deformation spaces, and their associated free splitting complexes.

In the present paper we study the local properties of the displacement function. In particular,
we study its convexity properties and the behaviour at bordification points, by 
geometrically characterising its continuity-points. We prove that the global-simplex-displacement
spectrum of $Aut(F_n)$ is a well-ordered subset of $\mathbb R$, this being helpful for
algorithmic purposes. We introduce a weaker notion of train tracks, which we call {\em partial train tracks} (which
coincides with the usual one for irreducible automorphisms) and we prove that, for any automorphism, points of minimal displacement - minpoints - coincide with the marked metric graphs that support partial train tracks. We show that any
automorphism, reducible or not, has a partial train track (hence a minpoint) either in the outer space or its bordification. We show that, given an automorphism, any of its invariant free factors is seen in a partial train track map. 

In a subsequent paper we will prove that level sets of the displacement functions are
connected, and we will apply that result to solve certain decision problems.
\\\ \\
Note: the two papers where originally packed together in the preprint 
\verb|arxiv:1703.09945.| 
We decided to split that paper following the recommendations of a referee.  
\\\ \\
Current date: \today
\end{abstract}
\maketitle
\tableofcontents

\section{Introduction}

\subsection{Overview}

Let $F_n$ denote the free group of rank $n$ and $\Out(F_n)=\Aut(F_n)/\operatorname{Inn}(F_n)$ be
the group of outer automorphisms. The natural space upon which this acts is $CV_n$, Culler-Vogtmann Space, which in turn admits a (non-symmetric) metric, the Lipschitz metric (see \cite{FM11}, \cite{{MR3342683}}).

The motivation for this paper is to extend the main results of \cite{FM13}, and in particular the result that for irreducible elements, $\phi$, of  
$\Out(F_n)$, the points that are minimally displaced by $\phi$ in $CV_n$ (with respect to the Lipschitz metric) coincide exactly with the points that support train track representatives for $\phi$. 

In order to do this, we need to extend our space to the free splitting complex, $\mathcal{FS}_n$ -  see \cite{handelmosher} and \cite{hatcher} for more details, and a complete exposition. However, our methods rely on inductive arguments which require us to deal with a more general setting, of a deformation space and its free splitting complex, as follows.

Let $\Gamma$ be a free product of groups with a specified free splitting. That is, abusing notation, $\Gamma$ is a group $G_1* \ldots G_k*F_r$, where the $F_r$ is a free group and the $G_i$ are specified. This is not, necessarily, the Gruschko decomposition, and we allow the $G_i$ to be freely decomposable, or even free. Equivalently, we can think of $\Gamma$ as given by a specific graph of groups with trivial edge groups. (We shall suppress this here, but we also need to allow $\Gamma$ to be disconnected in general). 

The space $\O(\Gamma)$ is then the deformation space of this splitting. That is, the space of all the simplicial, edge-free trees on which $\Gamma$ acts with the same elliptic elements as the defining splitting. We endow these trees with $\Gamma$-metrics, and identify two trees when there is an equivariant isometry between them (or an equivariant homothety if we don't insist on having volume 1). The space $\overline{\O(\Gamma)}$ is the free splitting complex associated to this; the space of all edge-free, simplicial metric trees whose elliptic elements contain all the elliptic elements from the original splitting (but may contain more). Again, trees are identified if they are equivariantly isometric. 

More concretely, we work almost entirely with the natural simplicial structure on $\O(\Gamma)$, since each point is given by a graph of groups and induces a cone on an open simplex by varying the lengths of edges. Technically, varying the lengths of edges in the graph of groups produces a simplex if we impose the condition that the volume is one, or if we work projectively; these correspond to the spaces $\O_1(\Gamma)$ and $\p\O(\Gamma)$. If the number of edges in the graph of groups is $k$, this produces an open $k-1$ simplex. In $\O(\Gamma)$ we get a cone on this open $k-1$ simplex. (See \ref{outerspace} and \ref{simplices}.)  Since the cone on an open simplex is homeomorphic to a open simplex of one dimension higher, we shall abuse notation and simply refer to simplicies in $\O(\Gamma)$.

One can then view $\overline{\O(\Gamma)}$  as the simplicial closure of $\O(\Gamma)$; a simplex
in $\O(\Gamma)$ has faces which correspond to collapsing various subgraphs. When such a
subgraph carries a hyperbolic element, the resulting quotient object defines a tree (graph of
groups) with more elliptic elements, and hence a point of $\overline{\O(\Gamma)}$. All points
of $\overline{\O(\Gamma)}$ arise in this way. We use the notation $\partial \O(\Gamma)$ to
denote the points in $\overline{\O(\Gamma)}$ which are not in $\O(\Gamma)$. Since points of
$\partial\O(\Gamma)$ can be considered as {\em points at infinity} of $\O(\Gamma)$, we often
use the notations $$\partial_\infty\O(\Gamma)\text{ and }\overline{\O(\Gamma)}^\infty.$$
(See Definitions~\ref{defn2020_2.29} and~\ref{boundary} for more details).
Any such $\Gamma$ has an associated group, $\Out(\Gamma)$ of (outer) automorphisms of the group which preserve the elliptic elements, and this groups acts on $\O(\Gamma)$ by isometries with respect to the Lipschitz metric.

In the case that $\Gamma=F_n$; that is, the trivial splitting of the free group where every non-identity element is hyperbolic, we obtain ${\O(\Gamma)}=CV_n$ and $\overline{\O(\Gamma)} = \mathcal{FS}_n$, and the associated automorphism group is $\Out (F_n)$.

\bigskip

Any $\phi\in\Out(F_n)$ acts on
$\O(F_n)$ and induces a displacement function $\lambda_{\phi}:\O(F_n)\to [1,\infty)$
$$\lambda_\phi(X)=\Lambda(X,\phi X)$$
where $\Lambda$ denotes the (multiplicative, non symmetric) Lipschitz distance. 

\medskip

We extend this function to $\overline{\O(F_n)}^\infty$. 

\medskip

If 
$X\in\O(F_n)$ exhibits a $\phi$-invariant sub-graph $A$, the collapse of $A$ defines a point
$X/A\in\partial_\infty \O(F_n)$, whose displacement is finite, since $\phi(X/A)$ is a well defined
point, again in $\partial_\infty\O(F_n)$ (and carrying the same set of hyperbolic/elliptic elements). In other words, if we let $\Gamma$ denote the induced free splitting of $F_n$ arising from the collapse of $A$, then both $X/A$ and $\phi(X/A)$ belong to $\O(\Gamma)$, and hence are at finite Lipschitz distance. 

By setting
$\lambda_\phi(X/A)=\infty$ for those points $X/A\in\partial_\infty \O(F_n)$ whose collapsed part is not $\phi$ invariant, we have $\lambda_\phi$
defined on the whole $\overline{\O(F_n)}^{\infty}$ (although some points have infinite displacement).

\medskip

The same process works for any $\Gamma$ in place of $F_n$ and we study these all at the same time. The advantage of this is that we can apply inductive arguments, which turn out to be key in understanding the properties of $\lambda_\phi$. In particular, we prove that minimally displaced points are characterised in terms of (partial) train-track maps, and any automorphism
has a minpoint in $\overline{\O(\Gamma)}^\infty$, though not necessarily in $\O(\Gamma)$ - Theorem~\ref{Theoremtt} and Theorem~\ref{corttE}.

We also study the (failure of the) continuity of the function $\lambda_{\phi}$ on $\overline{\O(\Gamma)}^\infty$ and characterize the points at which it is not continuous - the `jumping' points, Theorem~\ref{newjump} and Theorem~\ref{corlalx}. We describe some of these results in more detail below. 

\bigskip

%


\subsection{Anticipating the results}

The main tool for studying $\lambda_\phi$ is to use good representatives for $\phi$. Namely, given
$X\in\O(F_n)$ (or in any $\O(\Gamma)$), we need  to find the best Lipschitz maps representing $\phi$ (that is to say $f:X\to
X$ so that $f_*=\phi$ on $\pi_1(X)$). All maps we use will be {\em straight}, meaning that have
constant speed on edges (hence they are determined by the image of vertices).
It is classical that one may always find an optimal map
$f:X\to \phi X$, whose Lipschitz constant satisfies $$\Lip(f)=\Lambda(X,\phi X).$$
However the usual proof, by means of Ascoli-Arzel\`a, is not constructive, nor quantitative. Our
first result is Theorem~\ref{Lemma_opt} which can be stated as follows, and gives a constructive proceedure - via a flow - for making a straight map optimal, and crucially adds a quantative bound to the process. 

\bigskip
\parbox{0.9\textwidth}{{\bf Theorem} (Optimization). Given $X,Y\in\O(F_n)$ and $f:X\to Y$ a Lipschitz map,
  there is a  Lipschitz map $g:X\to Y$ so that $\Lip(g)=\Lambda(X,Y)$ and so that
  $d_\infty(g,f)\leq \vol(X)(\Lip(f)-\Lambda(X,Y))$.}


\bigskip
The estimate arising from this theorem will be crucial in many proofs. For any {\em straight} map $f:X\to
Y$, the {\em tension graph} of $f$, denoted by $X_{\max}$, is the sub-graph of $X$ whose edges
are maximally stretched. We introduce the notion of {\em partial train track} map as a straight
map $f:X\to \phi X$ such that there is an invariant sub-graph $A\subseteq X_{\max}$ (not
necessarily proper) so that the restriction of $f$ to $A$ is a train track map in the usual sense.
Our study of displacement functions is based on the use of partial train tracks. The first result
on partial train tracks is that they characterise minimally displaced points (See
Theorem~\ref{Theoremtt} for a precise statement):

\bigskip
\parbox{0.9\textwidth}{{\bf Theorem}. For any automorphism $\phi$, local minima for
  $\lambda_\phi$ are global minima and consist exactly of those points supporting a partial
  train track.}

\bigskip

One of the main problems is that $\lambda_\phi$ is not continuous at the boundary points of
$\O(F_n)$. We say that $X\in\partial_\infty\O(F_n)$ {\em has not jumped} if there is a sequence
$X_i\to X$ of points $X_i\in\O(F_n)$, all contained in a single simplex, such that $\lambda_\phi(X_i)\to \lambda_\phi(X)$ (precise
definitions are given in Definition~\ref{defnojumpr}). In
Sections~\ref{Sectionbordification} and~\ref{SectionTTatinfinity} we give a complete description of jumping and non-jumping
points.  For instance, if we set
$\lambda(\phi)=\inf_{X\in\O(F_n)}\lambda_\phi(X)$ we get (Theorem~\ref{corlalx}, see also
Theorem~\ref{thmjump} for related statements)

\bigskip
\parbox{0.9\textwidth}{{\bf Theorem}. $X\in\partial_\infty\O(F_n)$ has not jumped if and only
  if $\lambda_\phi(X)\geq \lambda(\phi)$.
              }

\bigskip

In particular,

\bigskip
\parbox{0.9\textwidth}{{\bf Theorem} (Corollary~\ref{corminpoindonjump}). For any automorphism
  $\phi$, if $X\in\partial_\infty \O(F_n)$ is a min-point for $\lambda_\phi$ (i.e.
  satisfies $\lambda_\phi(X)=\lambda(\phi)$) then $X$ has not jumped.
}

\bigskip

For any $\phi$ we give the notion of {\em partial train track at infinity} as points
$X\in\partial_\infty\O(F_n)$ which have not jumped {\bf and} are partial train tracks for the
induced automorphism in the deformation space of $X$.  In Section~\ref{SectionTTatinfinity} we
prove that partial train tracks at infinity exist and are min-points. We prove in particular
the existence of (non-jumping) min-points in the bordification of outer space. 

\bigskip
\parbox{0.9\textwidth}{{\bf Theorem} (Theorem~\ref{corttE}). Any automorphism has a partial
  train track in $\overline{\O(F_n)}^\infty$. Partial train tracks (at infinity or otherwise) are min-points
  for the displacement function.  In particular, (non-jumping) min-points always exist.
}

\bigskip

The existence of partial train tracks also give information on invariant free factors:

\bigskip
\parbox{0.9\textwidth}{{\bf Theorem} (Theorems~\ref{strongcorred}). 
For any automorphism $\phi$, any $\phi$-invariant free factor of $F_n$ is visible in some partial train track.
}

\bigskip

And as in the irreducible case, existence of partial train tracks allows one to easily deduce that for
any automorphism we have $\lambda(\phi^n)=\lambda(\phi)^n$. (Corollary~\ref{corollaryphin}.)

\bigskip

\begin{rem*}[Connection with Relative Train Track Maps]
	There is a connection between relative train track maps and partial train track maps as follows: given the automorphism, $\phi$, one constructs a relative train track map as in \cite{BestvinaHandel}. Suppose that $\lambda$ is the maximum Perron-Frobenius eigenvalue for any stratum, and that the highest stratum in which it occurs is the $r^{th}$ one. (That is, that the Perron-Frobenius eigenvalue is strictly greater than that of any higher stratum, and at least as great as that of any lower stratum). Now collapse the invariant subgraph $G_{r-1}$ - the union of all the strata below the $r^{th}$ stratum. This defines a point of the free splitting complex, where $\phi$ admits a representative supporting an invariant subgraph on which it is train track with expansion factor, $\lambda$. By making the volume of this subgraph sufficiently small, we can ensure that the Lipschitz constant of every other edge is strictly less than $\lambda$, and this is our partial train track at infinity. It then follows that $\lambda = \lambda(\phi)$.

	However, the important difference between the two objects is that partial train tracks characterise exactly the minimally displaced set, Theorem~\ref{Theoremtt}, whereas relative train tracks do not. 
\end{rem*}

\bigskip

The objects with which we work are usually not locally compact. This makes all
convergence arguments technically difficult. For controlling the  convergence and minimisation processes, in particular those
of Section~\ref{SectionTTatinfinity}, we make crucial use of the following result on
displacements. For any simplex of $\overline{\O(F_n)}^\infty$ define
$\lambda_\phi(\Delta)=\inf_{X\in\Delta}\lambda_\phi(X)$. Then we prove;

\bigskip
\parbox{0.9\textwidth}{{\bf Theorem}(Theorem~\ref{conj}).
For any $F_n$ (and in fact for any deformation space) the global simplex-displacement spectrum
$$\operatorname{spec}(F_n)= \Big\{\lambda_\phi(\Delta):[\phi]\in\Out(F_n), \Delta\text{ a
  simplex of }\overline{\O(F_n)}^\infty \text{s.t. } \lambda_\phi(\Delta)<+\infty \}$$ is well-ordered as a subset of $\mathbb R$.
In particular, for any $[\phi]\in\Out(F_n)$ the spectrum of possible minimal displacements $$\operatorname{spec}(\phi)= \Big\{\lambda_\phi(\Delta):\Delta\text{ a simplex of }
\overline{\O(F_n)}^\infty\text{such that } \lambda_\phi(\Delta)<+\infty\}$$ is well-ordered as a subset of $\mathbb R$.

}

\bigskip

Finally, we want to mention also Section~\ref{section_conv}, in which we give a detailed
description of useful convexity properties of displacement functions, for instance proving that the displacement function is quasi-convex along Euclidean segments - see Lemma~\ref{lconvexity} and Lemma~\ref{lconv2}.

\bigskip

\noindent
\textsc{Acknowledgements:} We would like to thank both the Universit\`{a} di Bologna and the
Universitat Polit\'{e}cnica de Catalunya, for their hospitality during several visits.

We would also like to thank the referee of the earlier version of this paper (when it was a
single paper together with~\cite{partII}) for many very helpful
comments, as well of the referee of the present paper for the very useful
comments and suggestions.

\section{Setting, notation, and general definitions}
\subsection{Motivation for new definitions}
First, we wish to motivate our definitions and the general setting. Our aim is to study automorphisms of free groups which are possibly reducible. (Although our results will apply to free products more generally). If $\Gamma$ is a marked graph with $\pi_1(\Gamma)=F$ a free
group, and $\phi\in\Out(F)$, then $\phi$ can be represented by a simplicial map (that is, a continuous map on the graph, sending vertices to vertices and edges to edge paths) $f:\Gamma\to\Gamma$. That is, $f$ represents $\phi$ if there is an isomorphism $\tau: F_n \to \pi_1(\Gamma)$ such that $\phi = \tau^{-1} f_* \tau$. (The reason we are working with outer automorphisms is that we do not keep track of basepoints). 

If $\phi$ is reducible, then it is possible that we may find a collection of
disjoint connected sub-graphs $\Gamma_1,\dots,\Gamma_k$ such that $f$ permutes the $\Gamma_i$'s. (We are guaranteed to find such a collection in {\em some} $\Gamma$). In order to study the properties of $\phi$ it may help to collapse such an invariant
collection. (In other words, in the study of reducible automorphisms, we are naturally led 
to study the simplicial bordification of the Culler-Vogtmann Outer space $CV_n$.)

If we want to keep track of all the relevant information, we will be faced with the study of some
particular kind of moduli spaces. Namely, moduli spaces of actions on trees with possibly
non-trivial vertex stabilizers (when we collapse the $\Gamma_i$'s) and the product of such spaces
(when we consider the restriction to $\phi$ to the $\Gamma_i$'s.)

The typical topological object we are concerned with is a disjoint union of
metric trees, where $G$ acts with finitely many orbits, but with possibly non-trivial
vertex-stabilizers. Therefore we will develop the paper in this general - free splitting - setting, but the reader is invited to restrict attention to the case of $CV_n$ and its bordification.

\subsection{Notation for free splittings}

Let $G=G_1*\dots* G_p* F_n$ be any free product of groups, where $F_n$ denotes the free group of
rank $n$ (we allow $n$ to be zero, in that case we omit $F_n$). We do not assume that the $G_i$'s are
indecomposable. The theory we are going to develop is general, but we are mainly interested in
the case where $G$ is itself a free group. (Thus, in general, this free product decomposition is not unique, since $G$ has many different
splittings as a free product.) 

\begin{defn}[Free Splittings]
  Given a group $G$, a {\em free splitting} $\G$ of $G$ is a pair $(\{G_i\},n)$ where $\{G_i\}$ is a
collection of subgroups  of
$G$ and $n$ is natural a number such that $G=G_1*\dots*G_p*F_n$.
Two splittings $(\{G_1,\dots,G_p\},n)$ and $(\{H_1,\dots,H_p\},m)$ of $G$ are considered to be {\em of the same type} if $m=n$ and, up to reordering factors, each $H_i$ is conjugate to $G_i$.
\end{defn}

\begin{rem}
  We admit the trivial splitting $G=F_n$, $(\emptyset, n)$. That is the splitting with no free factors groups. In
  this case our discussion will amount to considering the free group $F_n$ and the classical Culler-Vogmtann Outer space $CV_n$.
\end{rem}

\begin{rem}
	Free splittings are also referred to as {\em free factor systems} in the literature - originally introduced in \cite{BFH2}, and also used in \cite{handelmosher}, \cite{HM2} and \cite{horbez}.

	Our point of view here is to take a fixed free factor system - a free splitting - and form the deformation space of that. This consists of trees equipped with edge-free actions whose vertex stabilisers are the (conjugacy classes of) the elements of the free factor system. That is, one can form the space of all trees which give rise to the same free factor system. 
	
	One can also form the what is known as the free splitting graph or complex, which consists of all possible free splittings (and one can also make this relative to a base free splitting). This relative version is what we have in mind when we come to define our {\em simplicial bordification} (see~\ref{freesplitcomplex}).

	  Note that a ``splitting" in general refers to any action on a tree and the induced graph of groups decomposition, but no confusion should arise since all of the splittings we consider are ``free", in the sense that the edge stabilisers in the tree are trivial (equivalently, the splitting which arises is a free factor system). 
\end{rem}

\begin{defn}[Sub-splittings]\label{roas}
Let $\G=(\{G_1,\dots,G_p\},n)$  and $\mathcal S=(\{H_1,\dots,H_q\},r)$ be two
free splittings of $G$. We say that $\mathcal S$ is a {\em sub-splitting} of
$\G$ if each $H_i$ decomposes as $$H_i=G_{i_1}*\dots G_{i_{l}}*F_{s_i},$$ and $r+ \sum_i s_i=n$. 
\end{defn}

\begin{defn}[Kurosh rank of a free splitting] \label{kurosh}
  The Kurosh rank of the splitting $G=G_1*\dots*G_p*F_n$ is $n+p$. 
\end{defn}

\subsection{$G$-graphs and $G$-trees}

Given a group $G$, a simplicial $G$-tree is a 
  simplicial tree $T$ endowed with a faithful simplicial action of $G$. $T$ is minimal if it
  has no proper $G$-invariant sub-tree. In 
  particular, if $T$ is minimal then $G$ acts without global fixed points and $T$ has no leaves
  (valence one vertices).

We next define $\G$-trees and $\G$-graphs. For those familiar with Bass-Serre theory, these are
the trees dual to a given splitting and the corresponding graphs of groups.
  \begin{defn}[$\G$-trees and $\G$-graphs]\label{def20202.6}
    Let $G$ be a group $G$ and  $\G=(\{G_1,\dots,G_p\},n)$ be a splitting of $G$. A
      $\G$-tree is a metric simplicial $G$-tree $T$ such that
   \begin{itemize}
   \item $G$ acts isometrically on $T$,
  \item For every $G_i$ there is exactly one orbit of vertices whose stabilizer is conjugate to
    $G_i$. Such vertices are called {\em non-free}. Remaining vertices have trivial
    stabilizer and are called {\em free} vertices.
  \item $T$ has trivial edge stabilizers.
  \end{itemize}
A $\G$-graph is a finite connected metric graph of groups $X$ such that
\begin{itemize}
\item $X$ is {\em marked}; that is, there is a fixed isomorphism between the fundamental group of $X$ (as a graph of groups) and $G$. 
\item $X$ has trivial edge-groups;
\item the fundamental group of $X$ as a topological space is $F_n$;
\item the splitting given by the vertex groups is equivalent to $\G$.
\end{itemize}

We note that Bass-Serre theory gives a correspondence between $\G$-trees and $\G$-graphs. 

Two $\G$-trees are equivalent if there is an equivariant isometry between them; there is an analogous equivalence at the level of $\G$-graphs, which is harder to state but comes down to what is called a graph of groups morphism, and arises as the quotient map one gets from an equivariant isometry. 

  \end{defn}

\begin{rem}
	Recall that for an action on a (simplicial) tree, every group element either fixes a point or has an axis of minimal displacement. In the former case the element is called elliptic, and in the latter case hyperbolic. 
\end{rem}

\begin{nt}
  Throughout the paper, if $G$ has a splitting $\G$ which is clear from the context,
  then any $G$-tree is required to be a $\G$-tree. (And the same for graphs.)
\end{nt}

\begin{ex}
If $X$ is a finite connected graph of groups with trivial edge-groups, then denote the splitting induced by the vertex groups of $X$ by $\G$. It is clear that $\G$ is a splitting for $\pi_1(X)$, $X$ is a  $\G$-graph, and the Bass-Serre tree associated to $X$ is a $\G$-tree. 
\end{ex}

\begin{defn}[Core graph]
  A {\em core-graph} is a graph of groups whose leaves have
  non-trivial vertex-group. Given a graph $X$ we define $\core(X)$ to be the maximal core
  sub-graph of $X$. (If the vertex groups are all trivial, so that $X$ is simply a graph, then a core graph has no valence one vertices).
Note that $\core(X)$ is obtained by recursively cutting edges ending at leaves.
\end{defn}

Given a splitting $\G=(\{G_i\},n)$ of a group $G$ and $T$ a $\G$-tree, the
quotient $X=G\backslash T$ is a connected $\G$-graph. $T$ is minimal if and only if $X$ is a
finite core graph. Since in the paper we are dealing with both $\G$-graphs and $\G$-trees, we
introduce what we call the tilde-underbar notation.

\begin{nt}[Tilde-underbar notation]\label{tildeunderbar}
 Let $\G$ be a splitting of a group $G$.
If $X$ is a $\G$-graph, then $\wt X$  denotes its universal covering, which is a $\G$-tree. As usual, if $x\in X$ then $\wt x$ will denote  a lift of $x$ in $\wt X$. 

We will also often want to lift ``loops"; that is, given an element of $\pi_1(X)$ - the fundamental group of $X$ {\em as a graph of groups} - we lift the loop to a line in the universal cover. Concretely, this requires a description of the loop as a sequence of edges and vertex group elements which can be mirrored in the Bass-Serre tree; note that edges ``downstairs" are orbits of edges, so the vertex group element serves to determine which edge ``upstairs" in the orbit to follow. Similarly, we may lift paths, using the same process.

Equivalently, a lift of a group element is the axis, and the collection of lifts is the orbit of this axis (if one forgets basepoints, then there is no favourite lift, and one just has the orbit, of which one chooses a component.) 

For subsets, $A \subseteq X$, we wil usually mean a collection of edges and vertex groups - that is, the entire vertex group is a part of $A$ if at all. We can lift $A$ by taking all the lifts of the loops realised in $A$, and if $A$ is connected we usually refer to a component of this as a ``lift" of $A$.


Conversely, if $T$ is a minimal $\G$-tree we denote by $\ul T$ the
quotient $\G$-graph. We mirror this notation for points and subsets. 

\noindent
Hence, $\widetilde{\ul X}=X$ for both graphs and trees.
\end{nt}

\begin{defn}[$X$-graphs, trees and forests]\label{def2020_2.12}
	
Let $\G$ be a splitting of a group $G$.  
\begin{itemize}
\item 
 If $X$ is a $\G$-graph
 (resp. $\G$-tree), 
  then a $X$-graph (resp. $X$-tree) is just a $\G$-graph (resp. tree).
  Unless otherwise specified, given a finite connected graph of groups $X$ with trivial edge-groups, an $X$-graph is a
  $\G$-graph (and an $X$-tree is a $\G$-tree). 
\item 
 If $\Gamma=\sqcup \Gamma_i$ is a
  disjoint finite union of finite graphs of groups 
  with trivial edge-groups, a $\Gamma$-graph is a disjoint finite union $X=\sqcup X_i$ of
  $\Gamma_i$-graphs (and a $\Gamma$-forest is a union of $\Gamma_i$-trees). 
\item A $\Gamma$-sub-forest $\wt A$ of a $\Gamma$-forest $\wt X$ is the lift of a sub-graph of a
  $\Gamma$-graph $\ul X$. Here sub-graphs and lifts are ``full" in the sense that if a sub-graph
  contains a vertex, it contains the entire vertex group, as a sub-graph of groups --- and ``lift''
  means the full pre-image under the map from $\wt X\to \ul X$.
\item A sub-graph $A$ of a $\Gamma$-graph is {\bf non-trivial} if the fundamental group (as a
  graph-of groups) of any of its components contains a hyperbolic element. That is, if it contains a non-trivial element which is not just a vertex group element. A non-trivial $\Gamma$-sub-forest
  of a $\Gamma$-forest is a sub-forest obtained as a lift of a non-trivial sub-graph.
\end{itemize}
\end{defn}

\begin{defn}[Immersed loops]\label{dimmlor}
   A path $\gamma$ in a $\G$-graph $X$ is called {\em immersed} if it is has a lift $\wt \gamma$
   in $\wt X$ which is embedded. (Note that
$\gamma$ might not be topologically immersed in $X$ near non-free vertices.)
\end{defn}

\subsection{Outer spaces}
We briefly recall the definition of the outer space 
of a group $G$ corresponding to a splitting $\G$, referring to
\cite{FM13,GuirardelLevitt} for a detailed discussion of definitions and general properties.

\begin{defn}[Outer space]
	\label{outerspace}
  Let $G$ be a group and $\G$ be a splitting of $G$. The (projectivized) outer space of
  $G$, relative to the splitting $\G=(\{G_1,\dots,G_p\},n)$, consists of (projective)
  classes of minimal, simplicial, metric $\G$-trees, $X$ with no redundant vertex
  (that is, no valence two vertex is allowed to be free) and such that the $G$-action is by isometries. 

We use the notation $\O(G;\G)$ or simply $\O(\G)$ to denote the outer space
of $G$ relative to $\G$.
We use $\p\O(G;\G)$ (or simply  $\p\O(\G)$) to denote the 
projectivized outer space.  

For $X\in\O(\G)$ we define its (co-)volume $\vol(X)$ as the sum of lengths of edges in
$G\backslash X$. 
On occasion, 
we will need to work with the co-volume one slice of $\O(\G)$, which we denote by
$\O_1(\G)$.

\end{defn}

\begin{rem}
  If $\G$ is the trivial splitting of $G=F_n$, then $\O(\G)=CV_n$. 
\end{rem}

We stress here that the distinction between $\O(\G)$ and $\p\O(\G)$ is
not crucial in our setting as we will mainly work with scale-invariant functions.

\begin{rem}
The equivalence relation that defines $\p\O(\G)$ is the following: $X$ and $Y$ are equivalent if
there is an homothety (isometry plus a rescaling by a positive number) $X\to Y$ conjugating the
actions of $G$ on $X$ and $Y$.  In particu;ar, since $G$ acts isometrically on a metric $\G$-tree, the inner automorphisms of $G$ act trivially on $\O(\G)$ and $\p\O(\G)$.
\end{rem}

\begin{rem}
  If $G$ has a the simple splitting $G=G_1$, then $\O(\G)$ consists of a single element: a
  point stabilized by $G_1$, and in this case the equivalence relation is trivial. 
\end{rem}

\begin{rem}
If $X\in\O(\G)$, the quotient $\ul{X}$ is a metric core $\mathcal
G$-graph. Conversely, if $X$ is a core metric $\G$-graph with no redundant vertex, then $\wt X\in\O(\G)$.
\end{rem}

In the paper we will work with both graphs and trees. Strictly speaking we have defined $\O(\G)$ as
a space of trees, but we it will be often convenient to use graphs $X$ so that $\wt X\in
\O(\G)$. Clearly the two viewpoints are equivalent and we shall have occasion to abuse notation and switch between graphs and trees. However, when we wish to make the distinction clear, we will 
add a ``$gr$'' subscript to indicate that we are working with graphs. To illustrate:
$$\Og(\G)=\{\G\text{-graph } X: \wt X\in\O(\G)\}$$
The spaces $\Og(\G)$ and $\O(\G)$ are naturally identified via $\ul{X}\leftrightarrow \wt X$.

\begin{nt}
If $\ul X$ is a finite connected graph of groups with trivial edge-groups, and $\mathcal S$ is the
splitting of $\pi_1(\ul X)$ given by vertex-groups, then we set 
$$\O(\ul X)=\O(\pi_1(\ul X),\mathcal S).$$
  \end{nt}
Clearly, if $\ul X$ is a metric core graph of groups with no redundant vertices, then $\widetilde X\in\O(\ul X)$.


\begin{defn}\label{not:gamma}
We define $\Gamma$ and $\O(\Gamma)$ as follows:
	\begin{itemize}
		\item  $\Gamma=\sqcup \Gamma_i$ will always mean that $\Gamma$ is a finite disjoint
		union of {\em connected} finite graphs of groups $\Gamma_i$, each with trivial edge-groups and non-trivial fundamental group $H_i=\pi_i(\Gamma_i)$, each
		$H_i$ being equipped with the splitting given by the vertex-groups.
		\item Then $\O(\Gamma)$ is defined to be the product of the
                  $\O(\Gamma_i)$. That is, a point in $\O(\Gamma)$ is a tuple of minimal,
                  simplicial, edge-free isometric actions of the corresponding $H_i$ on metric
                  trees, up to equivariant isometry. The elliptic elements are precisely the
                  vertex groups in $\Gamma$.  
		\item There is a natural action of $\R^+$ on $\O(\Gamma)$ given by scaling each component by the same
		amount. The quotient of $\O(\Gamma)$ by such action is the projective outer space of $\Gamma$
		and it is denoted by $\p\O(\Gamma)$. (Thus $\p\O(\Gamma)$ is not the product of the $\p\O(\Gamma_i)$'s.)
		\item 
		The notion of co-volume extends to $\Gamma$-trees: If
		$X=(X_1,\dots,X_k)\in\O(\Gamma)$ we set $\vol(X)=\sum_i\vol(X_i)$, and $\O_1(\Gamma)$ denotes
		the co-volume $1$ slice of $\O(\Gamma)$.
              \item We tacitly identify $X=(X_1,\dots,X_k)\in \O(\Gamma)$ with the labelled
                disjoint union $X=\sqcup_i X_i$. So an element of $\O(\Gamma)$ can be
                interpreted as a metric $\Gamma$-forest (See Definition~\ref{def2020_2.12}).
	\end{itemize}

%
%


\end{defn}

\begin{rem}
  If $X$ is a $\G$-tree, then $\O(\ul X)=\O(\G)$. In other words, $\O(\G)$ is a particular case of
  $\O(\Gamma)$ with connected $\Gamma$. In the following we will therefore develop the
  theory for general $\O(\Gamma)$, as this includes the ``connected'' case $\O(\G)$ (and in
  particular, the $CV_n$ case).
\end{rem}


\begin{defn}[Rank]\label{pr4_rank}
  The Kurosh rank of a finite graph of groups with trivial edge-groups is
  the Kurosh rank of the splitting\footnote{See Definition~\ref{kurosh}.} induced on its fundamental group by the vertex-groups. If $\Gamma=\sqcup \Gamma_i$ we set $$\rank(\Gamma)=\sum_i\rank(\Gamma_i).$$
\end{defn}
By definition, the rank is a natural number greater or equal to one.
Note the the rank of a graph of groups $X$ is simply the rank of its fundamental group as
a topological space plus the number of non-free vertices.

We will also consider moduli spaces with marked points.
\begin{nt}
  Let $\G$ be a splitting of $G$. The moduli space of $\G$-trees with $k$ labelled points $p_1,\dots,p_k$ (not necessarily
  distinct) is denoted by 
  $\O(G;\G,k)$ or simply $\O(\G,k)$. If $X$ is a finite graph of groups with
  trivial edge-groups we set $\O(X,k)=\O(\pi_1(X),k)$.
If $\Gamma=\sqcup_{i=1}^s \Gamma_i$, given $k_1,\dots,k_s\in\N$ we
  set $$\O(\Gamma,k_1,\dots,k_s)=\Pi_i\O(\Gamma_i,k_i).$$
\end{nt}

\begin{nt}\label{not:Ogamma}
Let $X$ be a finite, connected graph of groups with trivial edge-groups and let $\G$ be the
splitting of $\pi_1(X)$ given by vertex groups. Let $A\subseteq X$ be a non-trivial sub-graph
(Definition~\ref{def2020_2.12}). Then $A$ induces a sub-splitting $\mathcal S_A$ of $\G$ where
the factor groups are either the fundamental groups of the components of $A$, or vertex groups
in $X\setminus A$.
\begin{itemize}
\item We denote by $X/A$ the graph of groups where components of $A$ are collapsed to points
  (a point for each component), and a vertex group is inserted - the fundamental group of the collapsed component.
\item We denote by $\O(X/A)$ the outer space $\O(\pi_1(X),\mathcal S_A)$.
\item We extend this definition to the case of $\Gamma=\sqcup_i\Gamma_i$ as in
  Definition~\ref{not:gamma} and $A\subseteq \Gamma$, by considering the splittings induced by
  $A$ on any component $\Gamma_i$, and we use notations $\Gamma/A$ and $\O(\Gamma/A)$ for the
  resulting (disjoint union of) graphs and its outer space.  
\item We borrow the same notation when speaking of $\Gamma$-forests (or $\G$-forests) and non-trivial sub-forests.
\end{itemize}
\end{nt}

\subsection{Simplicial structure}
The simplicial structure we are going to use is the one familiar to experts - see \cite{cv} and \cite{GuirardelLevitt}. Since we want to study
the simplicial bordificiation of our outer spaces, we need to introduce faces ``at infinity''
and a suitable notation for distinguish them from usual finitary faces. Faces ``at infinity''
of $\O(\G)$, will be in fact simplices in the outer space of some sub-splitting of $\G$.

\begin{defn}[Open simplices]
	\label{simplices}
Given a $\G$-tree $X$,  the open simplex $\Delta_X$ is the set of $\G$-trees
equivariantly homeomorphic to $X$. If $X$ is a $\G$-graph, then we agree that $\Delta_X$ is the
set of graphs obtained by quotients of elements of $\Delta_{\wt X}$.
The Euclidean topology on $\Delta_X$ is given by assigning a
$G$-invariant positive length $L_X(e)$ to each edge $e$ of $X$. Therefore, if $X$ has $k$ orbits of
edges, then $\Delta_X$ is isomorphic to the
standard open $(k-1)$-simplex if we work in $\p\O(\G)$ or $\O_1(\G)$, and to the positive cone over it if we
work on $\O(\G)$. 

\medskip

Since the positive cone of a $k-1$-simplex is homeomorphic to an open $k$-simplex, we shall abuse notation and simply refer to simplices in $\O(\Gamma)$. 

\medskip

Given two elements $X,Y$ in the same simplex
$\Delta\subset \O(\G)$ we define the {\bf
  Euclidean} sup-distance $d_\Delta^{Euclid}(X,Y)$ ($d_\Delta(X,Y)$ for short)
$$d_\Delta^{Euclid}(X,Y)=d_\Delta(X,Y)=\max_{e\text{ edge}}|L_X(e)-L_Y(e)|.$$
\end{defn}

Such definitions extend to the case of $\Gamma=\sqcup_i\Gamma_i$. We refer the  Definition~\ref{not:gamma}.

\begin{defn}[Euclidean topology]
    If $X=(X_1,\dots,X_k)\in\O(\Gamma)$, the simplex $\Delta_X$ is the set of $\Gamma$-forests 
    equivariantly homeomorphic to $X$ (component by component).
The Euclidean topology and  distance on $\Delta_X$ are defined by $$d_\Delta(X,Y)=\sup_i d_{\Delta_{X_i}}(X_i,Y_i).$$
\end{defn}

We note that the simplicial structure of $\p\O(\Gamma)$ is not the product of the structures of
$\p\O(\pi_1(\Gamma_i))$.

\begin{defn}[Faces and closed simplices]\label{def:face}
Let $X$ be a $\Gamma$-graph and let $\Delta=\Delta_X$ be the corresponding
open simplex. Let $F\subset X$ be a forest whose trees each contains at most one non-free
vertex. The collapse of $F$ in $X$ produces a new $\Gamma$-graph, whence a
simplex $\Delta_F$. Such a simplex is called a {\em face} of $\Delta$.

  The {\em closed} simplex $\overline{\Delta}$ is defined by
$$\overline{\Delta}=\Delta\cup\{\text{all the faces of $\Delta$}\}.$$
\end{defn}

\subsection{Simplicial bordification}

  There are two natural topologies on $\O(\Gamma)$, the simplicial one and the equivariant Gromov
  topology, which are in general different. Here we will
  mainly use the simplicial topology. We notice that if $\Delta$ is an open simplex, then
  the simplex $\overline\Delta$ is not the standard simplicial closure of $\Delta$, because  not all
  its simplicial faces are {\em faces} according to Definition~\ref{def:face}. This is because
  some simplicial faces of $\Delta$ are not in $\O(\Gamma)$ as defined. Such faces are somehow ``at
  infinity'' and describe limit points of sequences in $\O(\Gamma)$.
We now give precise definitions
  to deal with these limit points.

  It will be convenient to start with describing the free splitting complex, relative to $\Gamma$. 
  
  \begin{defn}[Free Splitting Complex]
  	\label{freesplitcomplex}
  	Let $\Gamma$ be as in Definition~\ref{not:gamma}. That is, $\Gamma$ is a finite disjoint union of graphs of groups. Then the free splitting complex relative to $\Gamma$ is the set of tuples of minimal, simplicial, isometric, edge-free, actions on $H_i$-trees, where the set of elliptic elements {\em includes} all the vertex groups of $\Gamma$, and up to equivariant isometry.  
  	
  	\medskip
  	
  	\noindent
  	We denote this set, $\overline{\O(\Gamma)}$. 
  	
  	\medskip
  	
  	\noindent
  	As in, Definition~\ref{simplices}, by varying the edge lengths on a given (tuple of) trees we can produce an simplex in $\overline{\O(\Gamma)}$. (Recall that we are abusing notation and calling the positive cone on an open simplex, an open simplex).  
  \end{defn}

  \begin{defn}[Simplices and Faces]\label{defn2020_2.29}
  	Let $\Delta$ be an open simplex in $\O(\Gamma)$. 
  	\begin{itemize}
  		\item $\overline{\Delta}$ is the closure of $\Delta$ in $\O(\Gamma)$, as in Definition~\ref{def:face},
  		\item $\partial_\O\Delta=\partial_\O\overline\Delta$ is the set-difference, $\overline\Delta\setminus\Delta$. We call this the {\em finitary boundary} of $\Delta$. The {\em finitary faces} of $\Delta$ are the faces which appear in $\partial_\O\Delta$.
  		\item $\overline\Delta^\infty$ is the closure of $\Delta$ in $\overline{\O(\Gamma)}$, 
  		\item $\partial_\infty\Delta=\partial_\infty\overline\Delta$ is the set-difference, $\overline\Delta^\infty \setminus \overline{\Delta}$. 
  	\end{itemize}
  	
  \end{defn}

%
%

Let $X$ be a $\Gamma$-graph and $\Delta=\Delta_X$. Let $A$ be a proper subgraph of $X$ having
at least a component  
which is not a tree with at most one non-free vertex. Equivalently, $\wt A$ contains the axis of a hyperbolic element.

Let $Y$ be the graph
of groups obtained by collapsing each component of $A$ to a point (different components to
different points). Then, $Y\in \O(X/A)$. The corresponding simplex $\Delta_Y$ is a simplicial
face of $\Delta_X$ obtained by setting  the edge-lengths of $A$ to zero. Note that
$\Delta_Y$ belongs to  $\O(X/A)$ and not to $\O(X)$. However, the simplicial topology naturally
defines a topology on $\Delta_X\cup \Delta_Y$, which we still name simplicial topology.

Also note that any point in $\partial_\infty\Delta$ is obtained by collapsing such a sub-graph, $A$. We may, in general, assume that all the components of $A$ are core-graphs (each component contains the axis of a hyperbolic element), if we are willing to replace $\Delta$ with a finitary face of $\Delta$.

%
%
%
%

\begin{defn}[Boundary at infinity]\label{boundary}
  We define the boundary at infinity and the simplicial bordification of $\O(\Gamma)$ as
$$\partial_\infty\O(\Gamma)=\bigcup_{\Delta\text{ simplex}}\partial_\infty\Delta.$$

Note that, 
$$ \overline{\O(\Gamma)}^\infty:= \overline{\O(\Gamma)}=\O(\Gamma)\cup\partial_\infty\O(\Gamma).$$
\end{defn}

\begin{rem}
	Note that all these operations can be carried out in the projective spaces, with the definitions essentially unchanged.
\end{rem}

\begin{rem}
	We note that when $\Gamma=F_n$, that is the splitting of the free group where every non-trivial element is hyperbolic, then we get that $\O(\Gamma)$ is simply Culler-Vogtmann space, $CV_n$ and the bordification, $\overline{\O(\Gamma)}$ is the free splitting complex, $\mathcal{FS}_n$.
\end{rem}

\subsection{Horoballs and regeneration}\label{sechor}
We recall Definition~\ref{not:gamma}.
\begin{defn}[Horoballs]

 Given $X\in\partial_\infty\O(\Gamma)$, 
$\Hor(X)$ is the set of marked metric trees, $Y \in \O(\Gamma)$,  such that $\ul X$ is obtained
from $\ul Y$ by collapsing a proper family of core sub-graphs. We set $\Hor(X)=X$ for $X \in
\O(\Gamma)$, by convention (and use $\Hor(\ul X)$ for graphs).

%
%
%
%
%
%
%
\end{defn}
In other words,  a metric graph $\ul Y$ is in the horoball of $\ul X$ if $\ul X$ is obtained
from $\ul Y$ by setting to zero the length of edges of a proper family
of core sub-graphs.  On the other hand, $\Hor(X)$ can be regenerated from $X$ as follows.

Suppose $X\in\partial_\infty\O(\Gamma)$. Thus there is a $\Gamma$-graph $\ul Y$ and a sub-graph
$\ul A=\sqcup_i\ul A_i\subset Y$ whose components $\ul A_i$ are core-graphs, and such that $\ul
X=\ul Y/\ul A$. Let
$v_i$ be the non-free vertex of $\ul X$ corresponding to $\ul A_i$. In order to recover a generic point
$Z\in\Hor(X)$, we need to replace each $v_i$ with an element $\ul V_i\in\O(\ul A_i)$. Moreover,
in order to completely define the marking on $\ul Z$, we need to know where to attach to $\ul
V_i$ the edges of $\ul X$
incident to $v_i$, and this choice has to be done in the universal covers $\wt{V_i}$. No more
is needed. Therefore, if $k_i$ denotes the valence of the vertex $v_i$ in $X$, we have
$$\Hor(X)=\Pi_i\O(A_i,k_i).$$
(Note that some $k_i$ could be zero, e.g.  if $A_i$ is a connected component of $Y$.)
There is a natural projection $\Hor(X)\to\O(A)$ which forgets the marked points. We will be mainly
interested in cases when we collapse $A$ uniformly, for that reason we will use the
projection to $\p\O(A)$: $$\pi:\Hor(X)\to\mathbb P\O(A)$$
where $\Hor(X)$ is intended to be not projectivized.

Note that if $[P]\in\mathbb P\O(A)$, then $\pi^{-1}(P)$ is connected because it is just
$\Pi_{i}(A_i^{k_i})$. Since $\O(A)$ is connected (as a product of connected spaces), then $\Hor(X)$ is connected.

\begin{rem}
Note that {\em the same} graph of groups $X$ can be considered as a point at
infinity of {\em different} spaces. If we need
to specify the space in which we work, we shall write $\Hor_\Gamma(X)$ (or $\Hor_\G(X)$.)
\end{rem}

\subsection{The groups $\Aut(\Gamma)$ and $\Out(\Gamma)$}
We are going to introduce the groups of automorphisms that preserve splittings, and their
generalizations to the case of non-connected graphs. 

\begin{defn}[Automorphism-groups of splittings]
  Let $G$ be endowed with the splitting $\G:G=G_1*\dots*G_p*F_n$.
The group of automorphisms of $G$
  that preserve the set of conjugacy classes of the $G_i$'s is
  denoted by $\operatorname{Aut}(G;\G)$. We set
  $\operatorname{Out}(G;\G)=\operatorname{Aut}(G;\G)/\operatorname{Inn}(G)$\footnote{Clearly $\operatorname{Inn}(G)\subset\Aut(G;\G)$.}.
\end{defn}

The group $\Aut(G,\G)$ acts on $\O(G)$ by
changing the marking (i.e. the action), and $\operatorname{Inn}(G)$ acts trivially. Hence
$\operatorname{Out}(G;\G)$ acts on
$\O(G;\G)$. If $X\in\O(G;\G)$ and
$\phi\in\operatorname{Out}(G;\G)$ then $\phi X$ is the same metric tree as $X$, but the action
is $(g,x)\to \phi(g)x$. The action is simplicial and continuous w.r.t. both simplicial and
equivariant Gromov topologies.

\medskip
We now extend the definition of $\Aut(G,\G)$ to the case of $\Gamma=\sqcup_i\Gamma_i$. We
denote by
$\mathfrak S_k$ the group of permutations of $k$ elements.
\begin{defn}[Splitting isomorphism-groups]
  Let $G$ and $H$ be two isomorphic groups endowed with splitting $\G:G=G_1*\dots G_p*F_n$ and
$\mathcal H:H=H_1*\dots H_p*F_n$. The set of isomorphisms from $G$ to $H$ that map each $G_i$
to a conjugate of one of the $H_i$'s is denoted by $\operatorname{Isom}(G,H;\G,\mathcal H)$. If
splittings  are clear from the context we write simply
$\operatorname{Isom}(G,H)$. 

\end{defn}

\begin{defn}[$\Aut(\Gamma)$]
	\label{aut}
 For $\Gamma=\sqcup_{i=1}^k\Gamma_i$ as in Definition~\ref{not:gamma},
  we set $$\Aut(\Gamma)=\{\phi=(\sigma,\phi_1,\dots,\phi_k):\ \sigma\in\mathfrak S_k\text{
and } \phi_i\in\operatorname{Isom}(H_i,H_{\sigma_i})\}.$$
\end{defn}

The composition of $\Aut(\Gamma)$ is component-wise, defined as follows. Given
$\phi=(\sigma,\phi_1,\dots,\phi_k)$ and $\psi=(\tau,\psi_1,\dots,\psi_k)$ we have
$$\psi\phi=(\tau\sigma,\psi_{\sigma(1)}\phi_1,\dots,\psi_{\sigma(k)}\phi_k)$$
\begin{rem}
  Not all permutations appear. For instance, if the groups $H_i$ are pairwise non-isomorphic,
  then the only possible $\sigma$ is the identity.
\end{rem}

\begin{defn}[$\operatorname{Inn}(\Gamma)$ and $\Out(\Gamma)$]
  We set:
$$\operatorname{Inn}(\Gamma)=\{(\sigma,\phi_1,\dots,\phi_k)\in\Aut(\Gamma): \sigma=id,
\phi_i\in\operatorname{Inn}(H_i)\}$$
$$\Out(\Gamma)=\Aut(\Gamma)/\operatorname{Inn}(\Gamma).$$
\end{defn}

\begin{ex}
If $X$ is a $\G$-graph and $f:X\to X$ is a homotopy equivalence which leaves invariant a core subgraph $A$,
then $f|_A$ induces and element of $\Aut(A)$, and its free homotopy class an element of $\Out(A)$.
\end{ex}

The group $\Out(\Gamma)$ acts on $\O(\Gamma)$ as follows. If $X=(X_1,\dots,X_k)\in\O(\Gamma)$,
then each $X_i$ is an $H_i$-tree. If $(\sigma,\phi_1,\dots,\phi_k)\in\Aut(\Gamma)$ then
$X_{\sigma(i)}$ becomes an $H_i$-tree via the pre-composition of $\phi_i:H_i\to H_{\sigma(i)}$
with the $H_{\sigma(i)}$-action. We denote such an $H_{i}$-tree by
$\phi_iX_{\sigma(i)}$. With that notation we have
$\phi(X_1,\dots,X_n)=(\phi_{1}X_{\sigma(1)},\dots,\phi_{k}X_{\sigma(k)}).$
(We remark that despite the left-positional notation, this is a right-action.)
Since $\operatorname{Inn}(\Gamma)$ acts trivially on $\O(\Gamma)$, then the
$\Aut(\Gamma)$-action descends to an $\Out(\Gamma)$-action.

\section{Straight maps, gate structures, and optimal maps.}
In this section we describe the theory of maps between trees (or graphs) representing points in outer
spaces. We will simultaneously deal with the ``connected'' case $\O(\G)$ (for instance the classical $CV_n$) and the general case
$\O(\Gamma)$. 

In this section $G,\G$ and $\Gamma$ will be as in Definitions\ref{def20202.6} and~\ref{not:gamma}.

\subsection{Straight maps} Now we will mainly work with trees.

\begin{defn}[$\O$-maps in $\O(\G)$]
Let $X,Y\in\O(\G)$. A map $f:X\to  Y$ is called an $\O$-map if it is Lipschitz-continuous and
$G$-equivariant. The Lipschitz constant of $f$ is denoted by $\Lip(f)$.  
\end{defn}

We recall that we tacitly identify $X=(X_1,\dots,X_k)\in\O(\Gamma)$ with the labelled disjoint
union $\sqcup_i X_i$. Hence, if $X,Y\in\O(\Gamma)$, a continuous map $f:X\to Y$ is a collection
of continuous maps $f_i:X_i\to Y_{\sigma(i)}$, where $\sigma\in\mathfrak S_k$.

\begin{defn}[$\O$-maps in $\O(\Gamma)$]\label{defOmap}
Let $X=(X_1,\dots,X_k)$ and $Y=(Y_1,\dots,Y_k)$ be two elements of $\O(\Gamma)$.
A map $f=(f_1,\dots,f_k):X\to Y$ is called an $\O$-map if for each $i$ the map
 $f_i$ is an $\O$-map from $X_i$ to $Y_i$. (No index permutation here. Compare with Definition~\ref{maprepf}.)
\end{defn}

\begin{defn}[Straight maps\footnote{In previous papers of the authors, a straight map is called
  $PL$-map. As a referee pointed out, piece-wise linearity is a well-established notion in
  literature, which is slightly different from our notion (we don't allow subdivisions). For
  that reason we  decided to change our previous terminology.}]
Let $X,Y$ be two metric trees.
A Lipschitz-continuous  map $f:X\to Y$ is 
{\em straight} if it has constant speed on edges, that is to say, for any edge $e$
of $X$ there is a non-negative number $\lambda_e(f)$ such that for any $a,b\in e$ we have
$d_Y(f(a),f(b))=\lambda_e(f)d_X(a,b)$. If $X,Y\in\O(\G)$ then we require any straight map to be an
$\O$-map. A straight map between elements of $\O(\Gamma)$ is an $\O$-map whose components are
straight. If $X,Y$ are metric graphs, we understand that $f:X\to Y$ is a straight map if its
lift to the universal covers is straight.

\end{defn}

\begin{rem}\label{rem:24} $\O$-map always exist and the images of non-free vertices is
  determined a priori by  equivariance
  (see~\cite{FM13}).
  For any $\O$-map $f$ there is a unique straight map denoted by $\PL(f)$, which is homotopic,
  relative to vertices, to $f$. We have $\Lip(\PL(f))\leq \Lip(f)$.
\end{rem}

\begin{defn}[$\lambda_{\max}$ and tension graph]
  Let $f:X\to Y$ be a $\PL$-map. We set $$\lambda(f)=\lambda_{\max}(f)=\max_{e}\lambda_e(f)=\Lip(f).$$  We
  define the {\em tension  graph} of $f$ by
$$X_{\max}(f)=\{e \text{ edge of } X : \lambda_e(f)=\lambda_{\max}\}.$$
If there are no ambiguities on the map, we write $\lambda_{\max}$ instead of
$\lambda_{\max}(f)$ and $X_{\max}$ for $X_{\max}(f)$.
\end{defn}

\begin{defn}[Stretching factors]
For $X,Y\in\O(\Gamma)$ we define
$$\Lambda(X,Y)=\min_{f:X\to Y\ \O\text{-map}}\Lip(f)$$
\end{defn}

The theory of stretching factors is well-developed in the connected case (i.e. for $CV_n$ or
general free products),
but one can readily see that connectedness of trees plays no role, and the theory extends without
modifications to the non-connected case. In fact,
 $\Lambda$ is well-defined, (see~\cite{FM11,FM13}for details) and it satisfies the
 multiplicative triangular inequality: $$\Lambda(X,Z)\leq\Lambda(X,Y)\Lambda(Y,Z)$$
It can be used to define a non-symmetric metric $d_R(X,Y)=\log(\Lambda(X,Y))$ and its
symmetrized version $d_R(X,Y)+d_R(Y,X)$ (see~\cite{FM11,FM12,FM13} for details) which
induces the Gromov topology. The group $\operatorname{Out}(\Gamma)$ acts by isometries on
$\O(\Gamma)$.

Moreover, there is an effective way to compute $\Lambda$, via the so-called
``sausage-lemma'' (see~\cite[Lemma 3.14]{FM11},\cite[Lemma 2.16]{FM12} for the classical case,
and~\cite[Theorem 9.10]{FM13} for the case of trees with non-trivial vertex-groups). We briefly
recall here how it works.

Let $X,Y$ be metric $\Gamma$-graphs. 
Any non-elliptic element $\gamma\in\pi_1(\Gamma)$ (i.e. an element not in a
vertex-group) is represented by an immersed loop $\gamma_X$ in $X$ and one $\gamma_Y$ in
$Y$ (see Definition~\ref{dimmlor}). The loop $\gamma_X$ (or, rather, its lift to $\wt X$) is usually called {\bf axis} of $\gamma$ in $X$
(or in $\wt X$) and corresponds to the points of minimal translation of $\gamma$ in $\wt X$.
The  lengths $L_X(\gamma_X)$ and $L_Y(\gamma_Y)$ are then the minimal 
translation lengths of the element $\gamma$ acting on $\wt X$ and $\wt Y$, respectively. (So
$L_X(\gamma_X)=L_X(\gamma)$ and $L_Y(\gamma_Y))=L_Y(\gamma)$.)  We can define the stretching
factor of $\gamma$ as $L_Y(\gamma)/L_X(\gamma)$. Then $\Lambda(\wt X,\wt Y)$ is the minimum of
the stretching factors of all non-elliptic elements. (Recall we are using the tilde-underbar notation~\ref{tildeunderbar}.)

\begin{thm}[Sausage Lemma~{\cite[Theorem 9.10]{FM13}}]\label{sausagelemma}
  Let $X,Y,\in\Og(\Gamma)$. The stretching factor $\Lambda(X,Y)$ is realized by a loop
  $\gamma\subset X$ having  one of the following forms:
  \begin{itemize}
  \item Embedded simple loop $O$;
  \item embedded ``infinity''-loop $\infty$;
  \item embedded barbel $O$--- $O$;
  \item singly degenerate barbel $\bullet$---$O$;
  \item doubly degenerate barbel $\bullet$---$\bullet$.
  \end{itemize}
(the $\bullet$ stands for a non-free vertex.) Such loops are usually named ``candidates''.
\end{thm}

\begin{rem}
  The stretching factor $\Lambda(X,Y)$ is defined on $\O(\Gamma)$ (or in the co-volume slice $\O_1(\Gamma)$) and not in
  $\p\O(\Gamma)$. However, we will mainly interested in computing factors of type
  $\Lambda(X,\phi X)$ (for $[\phi]\in\Out(\Gamma)$) and that factor is scale invariant.
\end{rem}

\begin{defn}[Gate structures] Let $X$ be any graph.
	A {\em turn} is a pair of germs of edges incident to the same vertex.
	
  A {\em gate} structure on $X$ is an equivalence relation on germs of edges, generated by (and in fact equal to) some collection of turns. Equivalence classes of germs are called {\em gates}.
 A turn is {\em
  illegal} if the two germs are in the same gate, it is {\em legal} otherwise. An immersed path
in $X$ is legal if it has only legal turns.

We can also consider gate structures on trees, where we require that the gate structures are equivariant, and hence descend to a gate structure on the quotient graph.

If $X=(X_1,\dots,X_k)\in\O(\Gamma)$ we require the equivalence relation to be $H_i$-invariant
on each $X_i$.
\end{defn}

Any straight map induces a gate structure as follows.
\begin{defn}[Gate structure induced by $f$]\label{pr4_gate} Given $X,Y\in\O(\Gamma)$ and a
  straight map $f:X\to Y$,
  the gate structure induced by $f$, denoted by $$\sim_f$$
  is defined by declaring equivalent two germs that have the
  same non-degenerate $f$-image.
\end{defn}

\begin{rem}[See~\cite{FM13}]\label{rem:2gated}
Let $X,Y\in\O(\Gamma)$ and be $f:X\to Y$ a straight map. If $v$ is a non-free vertex of $X$ and $e$
is an edge incident to $v$, then $e$ and $ge$ are in different gates for any $id\neq
g\in\operatorname{Stab}(v)$. (If $e$ is collapsed by $f$, then it is not equivalent to any other
edge by definition.)
\end{rem}

\begin{defn}[Optimal maps] Given $X,Y\in\O(\Gamma)$,
  a map $f:X\to Y$ is {\em weakly  optimal} if it is straight and $\lambda(f)=\Lambda(X,Y)$.

 A map $f:X\to Y$ is {\em optimal} if each vertex of $X_{\max}$ has at least two gates in
 $X_{\max}$ with respect to the gate
  structure induced by $f$.
\end{defn}

\begin{prop}\label{propELr}
	 A straight map between two $\Gamma$-forests is weakly optimal if and only if there is a
	periodic embedded legal line in the tension graph (i.e. a legal immersed loop in the quotient
	graph). In particular, optimal maps are  weakly optimal.

\end{prop}
	\proof First note that the Lipschitz constant of any straight map $f$ from $X$ to $Y$ provides an upper bound for the stretching factor of a loop. Hence, for any loop, $\gamma$,
	$$
	\frac{L_Y(\gamma)}{L_X(\gamma)} \leq \Lambda(X,Y)\leq\Lip(f). 
	$$
	
	Let $f:X \to Y$ be our straight map. Suppose first that we have an embedded legal line, $L \subseteq \widetilde{X}_{max}$. To say that $L$ is periodic means that $L$ is the axis of a hyperbolic element, $g$. Moreover, the axis of $g$ in $Y$ is contained in $f(L)$, and since $L$ is legal, the axis is exactly equal to $f(L)$ (as $f|_L$ is an embedding). Hence the stretching factor for $g$ is exactly the Lipschitz constant for $f$. Thus, 
	$$
	\Lambda(X,Y) \leq \Lip(f) = \frac{L_Y(g)}{L_X(g)} \leq \Lambda(X,Y). 
	$$
	Thus $f$ is weakly optimal.

	Conversely, suppose that  $f:X \to Y$ is weakly optimal. By the Sausage Lemma~\ref{sausagelemma}, we may find a loop, $\gamma$, whose stretching factor equals $\Lambda(X,Y)$.  Then, 
	$$
	\Lambda(X,Y) = \frac{L_Y(\gamma)}{L_X(\gamma)}  = \Lip(f).
	$$
	Let $L$ be the axis of $\gamma$ in $\widetilde{X}$. If either (i) $L$ is not legal or, (ii), $L$ is not a subset of $\widetilde{X}_{max}$, then $\frac{L_Y(\gamma)}{L_X(\gamma)}  < \Lip(f)$. Thus $L$ is our required line.

	\qed

In general optimal maps are neither unique nor do they form a discrete set,
 even if $X_{\max}=X$, as the following example shows. (If $X_{\max}\neq X$ then one can use
 the freedom given by the lengths of edges not in $X_{\max}$ to produce examples.)

\begin{ex}[A continuous family of optimal maps with $X_{\max}=X$]

Consider $G=F_2$. Let $X$ be a graph with three edges $e_1,e_2,e_3$ and two free vertices
$P,Q$, as in Figure~\ref{fig:ex3.14}. Set the length of $e_2$ to be $2$, name $x$ the length of
$e_1$, and $1+\delta$ that of $e_3$. The parameters $x,\delta$ will be determined below. For
any $t\in [0,1]$ consider the point $P_t$ at distance $1+t$ from $P$ along $e_2$, and the
point $Q_t$  at distance $1-t$ from $P$ along $e_3$. $P_t$ divides $e_2$ in oriented segments
$a_t,c_t$. $Q_t$ divides $e_3$ into $b_t,d_t$.
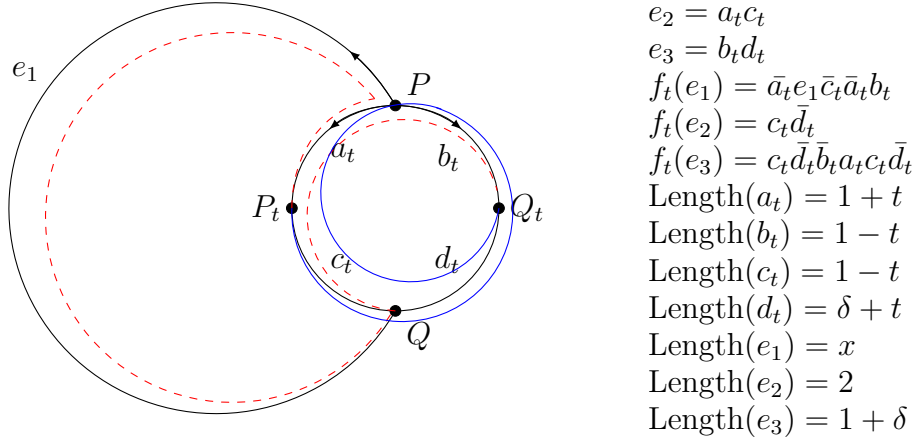
\begin{figure}[htbp]
  \centering
  \begin{tikzpicture}[x=5ex,y=5ex]
    \draw (0,0) circle [radius =1.5];
    \draw (0,1.5) arc (30:330:3) ;
    \draw [arrows=-{latex}] (0,1.5) arc (30:50:3) ; 
    \draw [arrows=-{latex}] (0,1.5) arc (90:130:1.5) ; 
    \draw [arrows=-{latex}] (0,1.5) arc (90:50:1.5) ; 
    \filldraw (0,1.5) circle(2pt) node[above right] {$P$};
    \filldraw (0,-1.5) circle(2pt) node[below right] {$Q$};
    \filldraw (-1.5,0) circle(2pt) node[left] {$P_t$};
    \filldraw (1.5,0) circle(2pt) node[right] {$Q_t$};
    \draw (-1.1,1.1) node[below right] {$a_t$};
    \draw (1.1,1.1) node[below left] {$b_t$};
    \draw (-1.1,-1.1) node[above right] {$c_t$};
    \draw (1.1,-1.1) node[above left] {$d_t$};
    \draw (-5,2) node[left] {$e_1$};
    \draw (3.5,-3.5) node[above right] {\parbox{20ex}{$e_2=a_tc_t$ \\ $e_3=b_t d_t$ \\ $f_t(e_1)=\bar{a_t}e_1\bar
        c_t\bar a_t b_t$\\ $f_t(e_2)=c_t \bar d_t$\\ $f_t(e_3)=c_t \bar d_t \bar b_t a_t c_t \bar d_t$\\
        Length$(a_t)=1+t$\\ Length$(b_t)=1-t$\\ Length$(c_t)=1-t$\\Length$(d_t)=\delta+t$\\
        Length$(e_1)=x$\\Length$(e_2)=2$\\Length$(e_3)=1+\delta$}};
    \draw[red, dashed] (-1.5,0) to [out=90,in=190] (-.3,1.6);
    \draw[red, dashed] (-.3,1.6) arc (40:330:2.7);
    \draw[red, dashed] (0,-1.5) arc (-95:-345:1.4);
    \draw[blue] (-1.5,0) arc(-182:72:1.6);
    \draw[blue] (1.5,0) arc(-10:-290:1.3);
  \end{tikzpicture}
  \caption{A continuous family of optimal maps with $X_{\max} = X$. The red dashed line is
    $f(e_1)$ and the blue line is $f(e_3)$ ($f(e_2)$ is not depicted).}
  \label{fig:ex3.14}
\end{figure} 
Consider the straight map $f:X\to X$ defined as in the figure, sending $P$ to $P_t$ and $Q$ to $Q_t$. If we collapse $e_3$, and we
homotope $P_t$ to $P$ along $a$, this corresponds to the automorphism $e_1\mapsto
e_1\overline{e_2}, e_2\mapsto \overline{e_2}$. 

The following direct calculation shows that if we set $\delta=1+2\sqrt2$ and $x=2\sqrt2$, the
map $f_t$ is optimal for any $t$ and all the three edges are stretched by the same amount.   

\medskip

The edges
 $e_1$ and $e_2$ are in different gates at $P$ and $e_1$ and $e_3$ are in different gates at
 $Q$. In order to check that $f_t$ is optimal it suffices to check that every edge is stretched
by the same amount. 
$$\lambda_{e_1}(f_t)=\frac{x+4}{x}\qquad \lambda_{e_2}(f_t)=\frac{1+\delta}{2}\qquad\lambda_{e_3}(f_t)=\frac{4+2\delta}{1+\delta}.$$
In particular they do not depend on $t$. If we set $x=2\sqrt2$ and $\delta=1+2\sqrt2$ we get
$$\lambda_{e_1}(f_t)=\frac{2\sqrt2+4}{2\sqrt2}\qquad \lambda_{e_2}(f_t)=\frac{2+2\sqrt2}{2}\qquad\lambda_{e_3}(f_t)=\frac{6+4\sqrt2}{2+2\sqrt2}$$
which are all equal to $1+\sqrt2$. \qed
\end{ex}

However, given a straight map, we can choose an optimal map which is in some sense the closest
possible. Given two $\O$-maps $f,g:X\to Y$ we define $$d_\infty(f,g)=\max_{x\in X} d_Y(f(x),g(x)).$$

\begin{thm}[Optimization]\label{Lemma_opt}
Let $X,Y\in\O(\Gamma)$ and let $f:X\to Y$ be a straight map. There is a map\footnote{We describe
  an algorithm to find the map $\wopt(f)$, but the algorithm will depend on certain choices, hence
  the map $\wopt(f)$ may be not unique in general.} $\wopt(f):X\to Y$
which is weakly optimal and such that
$$d_\infty(f,\wopt(f))\leq \vol(X)(\lambda(f)-\Lambda(X,Y))$$
Moreover, for any weakly optimal map $\f:X\to Y$ and for any $\varepsilon >0$ there is an optimal map
$g:X\to Y$ such that  $d_\infty(g,\f)<\varepsilon$.
\end{thm}


\proof[Proof of Theorem~\ref{Lemma_opt}.] By arguing component by component, we may assume without loss of generality that
$\Gamma$ is connected, hence that we can work in $\O(\G)$.
For this proof it will be convenient to work with both graphs and trees.
(Recall the tilde-underbar Notation~\ref{tildeunderbar}:  $\ul X=G\backslash X$, and
similarly for vertices and edges). By Remark~\ref{rem:2gated} a
non-free vertex will never be considered one-gated.

Let us concentrate on the first claim.

Let $\lambda=\Lambda(X,Y)$. Since straight maps are uniquely
determined by their value on vertices, we need only to define $\wopt(f)$ (and $g$) on vertices of
$X$. By Remark~\ref{rem:24} the image of non-free vertices is fixed.
We define straight maps $f_t$ for $t\in[0,\lambda_f-\lambda]$ by moving the images of all
one-gated vertices of $X_{\max}(f_t)$, in the direction given by the gate, so that $$\frac{d}{dt}\lambda(f_t)=-1.$$

Let us be more precise on this point. We define a flow which is piecewise linear, depending on
the geometry of the tension graph at time $t$.
The key remark to have in mind is that if an edge is not in
$X_{\max}(f)$, then it remains in the complement of the tension graph for small perturbations
of $f$. Therefore, we can restrict our attention to the tension graph.

Suppose we are at time $t$. We inductively define sets of vertices and edges as follows:

\begin{itemize}
\item $V_0$ is the set vertices of $X_{\max}(f_t)$ which are one-gated in
  $X_{\max}(f_t)$;
\item $E_0$ is the set of edges of $X_{\max}(f_t)$ incident to vertices in $V_0$. We agree that
  such edges contain the vertices in $V_0$ but not others. (If an edge has both vertices in
  $V_0$ then it contains both, otherwise it contains only one of its vertices.)
\end{itemize}
Having defined $V_0,\dots,V_i$ and $E_0,\dots,E_i$ , we define $V_{i+1}$ and $E_{i+1}$ as follows:
\begin{itemize}
\item $V_{i+1}$ is the set of
  one-gated vertices of $X_{\max}(f_t)\setminus\cup_{j=0}^iE_j$;
\item $E_{i+1}$ is the set of edges of $X_{\max}(f_t)\setminus\cup_{i=0}^iE_i$ incident to vertices in
  $V_{i+1}$. (As above such edges contain vertices in $V_{i+1}$ but not others.)
\end{itemize}
We notice  that since $G\backslash X$ is a finite $\G$-graph, we have only finitely many sets
$V_i$, say $V_0,\dots,V_k$ (each one formed by finitely many $G$-orbits).
\begin{lem}
  If $f_t$ is not weakly optimal, then $X_{\max}(f_t)\setminus\cup_{i=0}^kE_i$ is a (possibly empty)
collection of vertices, that we name terminal vertices.
\end{lem}
\proof
Note that no vertex in
$X_{\max}(f_t)\setminus\cup_{i=0}^kE_i$ can be one-gated, hence any vertex in
$X_{\max}(f_t)\setminus\cup_{i=0}^kE_i$ is either isolated or has at least two gates in
$X_{\max}(f_t)\setminus\cup_{i=0}^kE_i$. Thus if there is an edge $e$ in $X_{\max}(f_t)\setminus\cup_{i=0}^kE_i$,
the component of $\ul{X_{\max}(f_t)\setminus\cup_{i=0}^kE_i}$ containing $\ul e$
must also contain an immersed legal loop and so $f_t$ is weakly optimal.\qed

\medskip

By convention we denote the set of terminal
vertices  by $V_\infty$.

\begin{rem}\label{rem:vivj}
Any $e\in E_i$ has by definition at least one endpoint in
$V_i$, and  the
other endpoint is in some $V_j$ with $j\geq i$.
\end{rem}
Our flow is defined by equivariantly moving the images $f_t(v)$ of vertices in $X_{\max}(f_t)$. We need to define a
direction and a speed $s(v)\geq 0$ for any $f_t(v)$.

For $i<\infty$ each vertex in $V_{i}$ has a preferred gate: the one that survives in
$X_{\max}(f_t)\setminus\cup_{j=0}^{i-1}E_j$ ({\em i.e.} the unique gate of the map
$f_t|_{X_{\max}(f_t)}$ such that some edge of that gate is in $X_{\max}(f_t)\setminus \cup_{j=0}^{i-1}E_j$). That gate gives us the direction in which we move $f_t(v)$.

The idea is the following. Since a vertex in $V_0$ is one-gated, we can define the flow so as
to reduce the Lipschitz constant for every edge in $E_0$ (shrinking the image of each $E_0$
edge). Similarly, every vertex in $V_1$ is one gated in $X_{\max}(f_t)\setminus E_0$, so we
define the flow to reduce the Lipschitz constants of edges in $E_1$ and so on. We have only to
set speeds properly.

\begin{lem}
  There exists $G$-equivariant speeds $s(v)\geq 0$ such that if we move the
  images of any $v$ at speed $s(v)$ in the direction of its preferred gate, then
  for any edge $e\in X_{\max}(t)$ $$\frac{d}{dt}\lambda_e(f_t)\leq -1.$$ Moreover,
  for any $i$, and for any $v\in V_i$, either $s(v)=0$ or there is an edge $e\in E_i$ incident to
  $v$ such that $$\frac{d}{dt}\lambda_e(f_t)=-1.$$
\end{lem}
\proof
We start by choosing a total order on the set orbits of vertices of $X_{\max}(f_t)$
(i.e. on the set of vertices of $\ul{X_{\max}(f_t)}$) with
the only requirement that orbits of vertices  in $V_i$ are bigger than those in $V_j$ whenever $i>j$. This define a
partial order on vertices by declaring $w>v$ when $\ul w>\ul v$. Now, we define
speeds  recursively starting from the the biggest vertex and going down through the order.

The speed of terminal vertices is set to zero. Let $v$ be a
vertex of $X_{\max}(f_t)$ and suppose that we already defined the speed $s(w)$ for all
$w>v$.

The vertex $v$ belongs to some set $V_i$. For any edge $e\in E_i$ emanating from $v$ let
$u_e$ be the other endpoint of of $e$, and define a sign $\sigma_e(u_e)=\pm1$ as follows:
$\sigma_e(u_e)=-1$ if the germ of $e$ at $u_e$ is in the preferred gate of $u_e$, and
$\sigma_e(u_e)=1$ otherwise. (So, for example, $\sigma_e(u_e)=1$ if $u_e$ is terminal, and
$\sigma_e(u_e)=-1$ if $\ul v=\ul{u_e}$, or if $u_e\in V_i$.)

With this notation, if we move $f(v)$ and $f(u_e)$ in the direction given by their gates, and  at
speeds $s(v)$ and $s(u_e)$ respectively, then the derivative of $\lambda_e(f_t)$ is given by
$$-\left( \frac{s(v)-\sigma_e(u_e)s(u_e)}{L_X(e)} \right)$$

If $u_e>v$ we already defined its speed. We set
$$s(v)=\max\{0,\max_{u_e>v}\{L_X(e)+\sigma_e(u_e)s(u_e)\},\max_{\ul{u_e}=\ul v}\frac{L_X(e)}{2}\}$$

where the maxima are taken over all edges $e\in E_i$ emanating from $v$. Note that there may
exist some such edge with $\ul{u_e}<\ul{v}$. (By Remark~\ref{rem:vivj} in this case $u_e\in V_i$ (same $i$ as $v$), $\sigma_e(u_e)=-1$ and the
derivative of $\lambda_e$ will be settled later, when defining the speed of $u_e$.)

With the speeds defined in this way, we are sure that for any edge $e$ we have $d/dt
\lambda_e(f_t)\leq -1$ and, if $s(v)\neq 0$, then the edges that realize the above maximum
satisfy $d/dt \lambda_e(f_t)=-1$.
\qed

\medskip
The first consequence of this lemma is that if we start moving then $\lambda(f_t)$ decreases.
Locally in $t$, when we start moving,
the tension graph may lose some edges. However, the above lemma 
ensures  that any vertex $v$ with $s(v)\neq 0$ is incident to an edge $e$ which is maximally
stretched and $d/dt\lambda_e=-1$. Hence such an edge remains in the tension graph when we start
moving. Since $d/dt\lambda_e\leq -1$ for any edge in the tension graph, it follows that when we
start moving, the tension graph stabilizes. So our flow is well defined in $[t,t+\epsilon]$ for some
$\epsilon>0$. If at a time $t_1>t$ some edge that was not previously in $X_{\max}(f_t)$ becomes
maximally stretched, then we recompute speeds and we start again. A priori we may have to recompute
speeds infinitely many times $t<t_1<t_2<\dots$ but the control on $d/dt\lambda(f_t)$
ensures that $\sup t_i=T\leq \lambda_f-\lambda$. Since the speeds, $s(v)$, are uniformly bounded (one can take the number of edges in $G\backslash X$ 
multiplied by the maximum length of an edge, as an upper bound) the flow has a limit
for $t\to T$. More precisely, for any monotone sequence $t_n \to T$ as above, and any vertex $v$, the sequence $f_{t_n}(v)$ must be a Cauchy sequence and hence convergent, 
since all our trees are complete. Thus we can define $f_T(v) = \lim_{n \to \infty} f_{t_n}(v)$ for each vertex. This is enough 
to define a straight map, and then we can restart our flow from $T$. Therefore the set of times
$s\in[0,\lambda_f-\lambda]$ for which the flow  is well-defined for $t\in[0,s]$ is closed and
open and thus is the whole $[0,\lambda_f-\lambda]$.

With these speeds, we have $d/dt(\lambda(f_t))=-1$. Therefore
for $t=\lambda(f)-\lambda$, and not before, we have $\lambda(f_t)=\lambda$ hence
$f_t$ is weakly optimal. We define
$$\wopt(f)=f_{\lambda(f)-\lambda}.$$

We prove now the claimed estimate on $d_\infty(f,f_t)$. The $d_\infty$-distance between straight maps is bounded by the $d_\infty$-distance of their restriction to vertices.

We first estimate the speed at which the images of
vertices move. Let $S$ be the maximum speed of vertices,
i.e. $S=\max_v|s(v)|$. Let $v$ be a fastest vertex. Since it moves, it belongs
to $V_s$ for some $s<\infty$.
Let $v=v_1,v_2\dots,v_m$ be a maximal sequence of vertices such that:
\begin{enumerate}
\item $s(v_i)>0$ for $i<m$;
\item there is an edge $e_i$ between $v_i$ and $v_{i+1}$ such that $e_i\in E_a$ if $v_i\in V_a$;
\item $\sigma_{e_i}(v_{i+1})=1$ for $i+1<m$;
\item  $d/dt(\lambda_{e_i}(f_t))=-1$.
\end{enumerate}
By the above lemma, we have that either $s(v_m)=0$ or $\sigma_{e_{m-1}}(v_m)=-1$. Moreover, by
$(2)-(3)$ and Remark~\ref{rem:vivj} we have that $v_i<v_{i+1}$ and therefore the edges $\ul{e_i}$
are all distinct.

Let $\gamma$ be the path obtained by concatenating the $e_i$'s. By $(2)-(3)$,
$\gamma$ is a legal path in the tension graph. So let

$$L=\sum_i L_X(e_i)=L_X(\gamma)\qquad L_t=\sum_i L_Y(f_t(e_i))=L_Y(f_t(\gamma)).$$

Since the $e_i$'s are in the tension graph and by condition $(4)$ we have $$L_t=\lambda(f_t) L\qquad \frac{d}{dt}L_t=-L$$

On the other hand $-\frac{d}{dt}L_t\geq S$ because by $(3)$ the contributions of the speeds of $v_i$
do not count for $i=2,\dots,m-1$ and $f(v_m)$ either stay  or moves towards $f(v_1)$. It follows
that $$S\leq L\leq \vol(X).$$

It follows that for any vertex $v$ we have
$$d_Y(f(v),f_t(v))\leq\int_0^t\left|\frac{d}{ds}f_s(w)\right|ds\leq\int_0^tS=tS\leq t\vol(X)$$
hence
$$d_{\infty}(\wopt(f),f))=d_\infty(f_{\lambda(f)-\lambda},f)\leq (\lambda(f)-\lambda)\vol(X).$$

We prove the last claim of Theorem~\ref{Lemma_opt}. If $\f$ is optimal then we are
done. Otherwise, there is 
some one-gated vertex in $X_{\max}$. We start moving the one-gated vertices as described above,
by an arbitrarily small amount. Let $g$ be the map obtained, clearly we can make
$d_\infty(g,\f)$ arbitrarily small.
Since $\f$ is optimal, we must have $\lambda(g)=\lambda(\f)$. It follows that there is a
core sub graph of $\ul{X_{\max}}$ which survives the moving. In particular, every vertex of
$X_{\max}(g)$ is at least two-gated, hence $g$ is optimal.
\qed

\medskip

\begin{defn}
We denote by $\opt(f)$ any optimal map obtained from $\wopt(f)$ as described in the proof of Theorem~\ref{Lemma_opt}.
\end{defn}

\medskip

We want to stress the fact that Theorem~\ref{Lemma_opt} holds in a general context for $X,Y$
metric one-dimensional complexes where the notions of straight and optimal maps are generalized in
the obvious way.

\begin{prop}
  Let $A,B$ be metric one-dimensional simplicial complexes and let $f:A\to B$ a straight map. Then
  there is a weakly optimal map $\wopt(f)$ which is homotopic to $f$ relatively to $\partial
  A$, such that $$d_\infty(f,\wopt(f))\leq\vol(A)(\Lip(f)-\Lip(\wopt(f))).$$ Moreover, for any
  $\varepsilon >0$, there is an
  optimal map $g:A\to B$ homotopic to $f$ relatively to $\partial A$ such that $d_\infty(g,\wopt(f))<\varepsilon$.
\end{prop}
The proof is basically the same as that of Theorem~\ref{Lemma_opt} and it is left to the reader.
(We do not use this generalization in what follows, and simply register the result as it may be interesting to the reader.)
\medskip

Let $X,Y\in\O(\Gamma)$ and let $f:X\to Y$ be an optimal map. Let $v$ be a vertex of $X$ having an
$f$-illegal turn $\tau=(e_1,e_2)$. Since $f(e_1)$ and $f(e_2)$ share an initial segment, we can
identify an initial segment of $e_1$ and $e_2$. We obtain a new element $X'\in \O(\Gamma)$, with an
induced map, still denoted by $f$, from $X'$ to $Y$. This is a particular case of Stallings
fold (\cite{Sta}). We refer to~\cite{FM13} for further details.

\begin{defn} 
	\label{foldpath}
  We call the above operation a {\bf simple fold directed by $f$}. By a folding path directed by $f$
  we mean a sequence $(X_0,f_0),\dots,(X_n,f_n)$, where $f_0=f$, and
  $(X_{i+1},f_{i+i})$ is obtained by a simple fold directed by $f_i$.  
\end{defn}

 We finish this section by proving the existence of optimal maps with an additional
property, that will be used in the sequel.

\begin{defn}[Minimal optimal maps]
  Let $X,Y\in\O(\Gamma)$. An optimal map $f:X\to Y$ is {\em minimal} if its tension graph
  consists of the union of axes of maximally stretched elements it contains.
  In other words, if any edge $e\in X_{\max}$ is contained in the axis of some element in
  $\pi_1(X_{\max})$ which is maximally stretched by $f$.
\end{defn}
Note that not all optimal maps are minimal, as the following illustrates.

\begin{ex}
Let $X$ be the graph consisting of two barbels joined by an edge, as in
Figure~\ref{fig:ex3.22}. All edges have length one except the two lower loops that have length two.
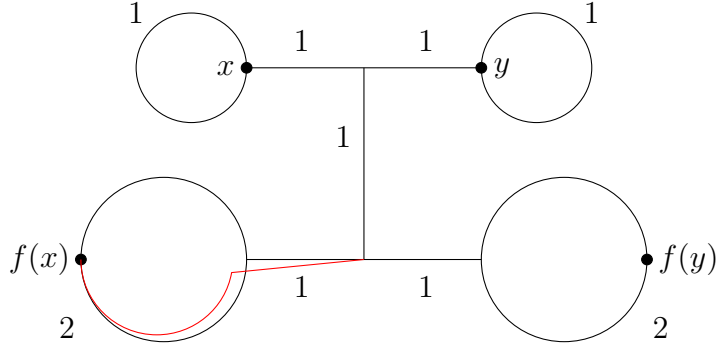
\begin{figure}[htbp]
  \centering
  \begin{tikzpicture}[x=1ex,y=1ex]
   \foreach \x in {0,25} \draw (\x,0) circle (4);
   \foreach \x in {0,29} \draw (\x-2,-14) circle (6);
   \draw (4,0) -- (21,0);
   \draw (4,-14) -- (21,-14);
   \draw (12.5,0) -- (12.5,-14);
   \foreach \x in {(-4,4), (29,4), (8,2), (17,2), (11,-5), (8,-16), (17,-16)} \draw \x node
   {$1$};
   \foreach \x in {(-9,-19), (34,-19)} \draw \x node {$2$};
   \foreach \x in {(4,0),(21,0),(-8,-14),(33,-14)} \filldraw \x circle(2pt);
   \draw (2.5,0) node {$x$} (22.5,0) node {$y$} (-11,-14) node {$f(x)$} (36,-14) node {$f(y)$};
   \draw[red] (-8,-14) arc (-180:-10:5.5) -- (12.5,-14);   
  \end{tikzpicture}
  \caption{A non-minimal optimal map. The dots $f(x)$ and $f(y)$ are not vertices, all other
    crossings are. The red line is the image of the left ``bar-edge'' of the top barbell. }
  \label{fig:ex3.22}
\end{figure}

Let $f:X\to X$ be the straight map that exchanges the the top and bottom barbells (preserving left and right) and maps $x$ to the middle point
of the lower left loop, and $y$ to the middle point of the lower right loop (see the figure).

The restriction of $f$ to the lower barbell is $1$-Lipschitz (each loop is shrunk and the bar is the same length as its image), while the stretching factor of
all top edges is two. Hence the tension graph $X_{\max}$ is the top barbel. The map is optimal
because all vertices of $X_{\max}$ are two gated, but the ``bar-edges'' of the top barbel are
not in the axis of any maximally stretched loop. This is because the only legal loops in
$X_{\max}$ are the two lateral loops of the barbell. Clearly this map can be homotoped to a map
with smaller tension graph. As the next theorem shows this is always the case for non-minimal
optimal maps.
\qed
\end{ex}

\begin{thm}\label{thmminopt}
  Let $X,Y\in\O(\Gamma)$ and let $f:X\to Y$ be an optimal map. If $f$ locally minimizes the
  tension graph amongst all optimal maps $X\to Y$, then $f$ is minimal. Moreover, given $g:X\to
  Y$ optimal, for any $\e>0$ there is a minimal optimal map $f:X\to Y$ with
  $d_\infty(g,f)<\e$.
\end{thm}
\proof The first claim clearly implies the second, because the tension graph is
combinatorially finite, hence the set of possible tension graphs is finite and we can always
locally minimize it.

We will prove the contrapositive, that if $f$ is not minimal then we can decrease the tension
graph by perturbations as small as we want. The spirit is similar to that of the proof of Theorem~\ref{Lemma_opt}.

As above, connectedness plays no role an we can work in $\O(\G)$ without loss of generality.
We will work with graphs rather than trees. For the ease of the reader we omit the underlines,
and we declare that $X,Y$ are $G$-graphs. Also we choose an orientation on edges, using the
classical bar-notation to indicate the inverse.

 At the level of 
graphs, the non-minimality of $f$ translates to the fact that there is an edge
$\alpha$ in the tension graph which is not part of any legal loop in $X_{\max}$.

Let $x$ be the terminal vertex of the oriented edge $\alpha$. 
We say that a path starting at $x$ is $\alpha$-legal, if it is a legal path in the tension graph, whose initial edge, $e$, is not in the same gate as  $\overline{\alpha}$. We say a loop at $x$ is $\alpha$-legal if, considered as paths, both the loop and its inverse are $\alpha$-legal.

If the terminal vertex of $\alpha$ admits an $\alpha$-legal loop and the initial point of
$\alpha$ also admits an  $\overline{\alpha}$-legal loop, then we can form the concatenation of
these loops with $\alpha$ to get a legal loop in the tension graph crossing $\alpha$ and
contradicting our hypothesis. 
(Note that an $\alpha$-legal loop need not be legal as a loop; that is, the lift of the loop to the tree need not be a legal line. We simply require that the loops can be concatenated in this way with $\alpha$ to form a legal loop.)

Hence, by reversing the orientation of $\alpha$ if necessary, we may assume that the endpoint $x$ (rather than the initial point) admits no $\alpha$-legal loops.

We will show
that it is possible to move the $f$-image of $x$ a small amount (and possibly some other vertices) so that we obtain an optimal map with smaller
tension graph.  Let $\e$ be small enough so that if an edge is not in $X_{\max}$, than it remains
outside the tension graph for any perturbation of $f$ by less than $\e$.

From now on, we restrict ourselves to the tension graph. We say that a vertex $v$ is legally seen from
$x$ if there is an $\alpha$-legal path $\gamma$ from $x$ to $v$. Note that in this case $v$ is free. Indeed, otherwise the path $\gamma$ followed by its
inverse can in fact be turned into  an $\alpha$-legal loop thanks to the action of the vertex
group 
($\gamma\overline\gamma$ has a legal ``lift'' to $\wt X$ defined by using the action of the
stabilizer of $\tilde v$). Since $v$ is free, we can move $f(v)$.

We want to chose a direction to move the images of vertices $\alpha$-legally seen from
$x$. First, the direction we choose for $f(x)$ is given by the gate of $\alpha$. That is, we
move $f(x)$ so as to reduce the length of $f(\alpha)$. For any vertex, $v$, $\alpha$-legally
seen from $x$, via a path $\gamma$, we move $f(v)$ backwards via the last gate of
$\gamma$. That is, we move $f(v)$ so as to retrace $\gamma$. Note that this direction depends
only on $v$ and not on the choice of $\gamma$. This is because, were there to be another
$\alpha$-legal path from $x$ to $v$, $\gamma'$, then the concatenation $\gamma
\overline{\gamma'}$ would define an $\alpha$-legal loop at $x$ unless the terminal edges of
$\gamma$, $\gamma'$ lie in the same gate. Hence directions are well defined.

Observe that if
the initial point $x_0$ of $\alpha$ is legally seen from $x$, and $\gamma$ is an $\alpha$-legal path
from $x$ to $x_0$, then the last edge of $\gamma$ must be in the same gate as $\alpha$, because
otherwise the concatenation of $\gamma$ and $\alpha$ would form a legal loop containing
$\alpha$. It follows that, whether $x_0$ can be legally seen from $x$ or not, in either case the
length of $f(\alpha)$ decreases when we move $x$ and possibly $x_0$. 

Next we move by $\e$ all the images of vertices legally seen from $x$, in the directions given
above. Consider an edge, $\beta$ (not equal to $\alpha$ or its inverse) in the tension
graph. If neither vertex of $\beta$ is $\alpha$-legally seen from $x$, then the image of
$\beta$ is unchanged and it remains in the tension graph. Otherwise, suppose that the initial
vertex of $\beta$ is $\alpha$-legally seen from $x$, via a path $\gamma$, whose terminal edge
is $\overline{\eta}$. If $\eta $ and $\beta$ are in different gates, then the terminal vertex
of $\beta$ is also $\alpha$-legally seen from $x$ and both vertices are moved the same amount,
such that the length of the image of $\beta$ remains unchanged. If, conversely, $\eta $ and
$\beta$ are in the same gate then either the length of the image of $\beta$ is reduced (if the
terminal vertex is not $\alpha$-legally seen) or it remains unchanged (if it is. For instance
if $\eta=\beta$.) Moreover, by our above observation, the length of the image of $\alpha$ must strictly decrease. In particular, $\alpha$ itself is no longer in the tension graph. 

On the other hand, since the tension graph has no one-gated vertices, there is at least one $\alpha$-legal
path emanating from $x$, an so some part of the tension graph survives.
Since $f$ is optimal, our assumption on $\e$ implies that the
new map is optimal  and  it  has a tension graph strictly smaller than $f$.\qed

\section{Displacement function and train track maps for automorphisms}
 For the rest of the section we fix $G,\G$ and
 $\Gamma=\sqcup_i\Gamma_i$ as in Definitions\ref{def20202.6} and~\ref{not:gamma}.
  (Recall that $CV_n$
is a particular case of $\O(\Gamma)$.) 
If not specified otherwise, $\phi=(\sigma,\phi_1,\dots,\phi_k)$ will be an element of $\Aut(\Gamma)$ - recall Definition~\ref{aut}.

This section is devoted to the study of train track maps, and related objects, from a metric point of view. In
particular, we prove that the points which are minimally displaced by $\phi$ are
exactly those admitting a {\em partial train track map} for $\phi$, see Definition~\ref{defttm}. (In the irreducible case, 
this amounts to showing that points of minimal displacement are precisely train track maps, in the usual sense. We broaden the class of maps to allow for the reducible case as well.)

The spirit of our analysis is that of~\cite{BestvinaBers,FM13}. We will recall the main facts proved
in~\cite{FM13} for irreducible elements of $\Out(G)$, and generalize such facts to the case
of $\Out(\Gamma)$. Connectedness does not really play a
crucial role, and most of the arguments of~\cite{FM13} transfer without requiring embellishment. The main
contribution of this section is to generalize from irreducible to reducible automorphisms.


\begin{defn}[Maps representing $\phi$]\label{maprepf}
  Let $X\in\O(\Gamma)$ and $\phi=(\sigma,\phi_1,\dots,\phi_k)$ be an automorphism. We say that a (straight) map $f:X\to X$ {\em represents $\phi$} if $f$
  maps $X_i$ to $X_{\sigma(i)}$, by a map we denote by $f_i$ which is equivariant in the following sense:

 For each $i=1 \ldots,k$ the map 
  
  $$f_i : X_i \to X_{\sigma(i)},$$ 
  
  is equivariant with respect to the isomorphism $\phi_i : H_i \to H_{\sigma(i)}$. This 
 means that for each $x \in X_i$ and each $h \in H_i$ we have,
  	
  	$$f_i(h \cdot x)  = \phi_i(h) \cdot f_i(x)$$.

  We also require that each $f_i$ is a is a straight map. We say that $f$ is optimal if each $f_i$ is optimal. 
  
%

If $X$ is a $\Gamma$-graph, then a map $f:X\to X$ represents $\phi$ if it has a lift 
$\wt f:\wt X\to \wt X$ representing $\phi$.
\end{defn}

Note that a map $f:X\to X$ representing $\phi$ can be viewed as an $\O$-map $f:X\to \phi X$.


\begin{defn}[Displacements]\label{defdispl}
  For any $[\phi]\in\Out(\Gamma)$ we define the function
$$\lambda_\phi:\O(\Gamma)\to\R\qquad \text{by} \qquad \lambda_\phi(X)=\Lambda(X,\phi X).$$
If $\Delta$ is a simplex of $\O(\Gamma)$ we define
$$\lambda_\phi(\Delta)=\inf_{X\in\Delta}\lambda_\phi(X)$$
If there is no ambiguity we write simply $\lambda$ instead of $\lambda_\phi$.
Finally, we set
$$\lambda(\phi)=\inf_{X\in\O(\Gamma)}\lambda_\phi(X)$$
\end{defn}

\begin{defn}[Minimally displaced points]
For any automorphism $\phi$ we define sets:
$$\Min(\phi)=\{X\in\O(\Gamma):\lambda(X)=\lambda(\phi)\}$$
$$\LocMin(\phi)=\{X\in\O(\Gamma):\exists U\ni X\text{ open s.t. } \forall Y\in U\  \lambda(X)\leq\lambda(Y)\}$$
\end{defn}

\begin{rem}[Fold-invariance of $\Min(\phi)$]
  A fold directed by a weakly optimal map does not increase $\lambda$. This is because the fold naturally induces a map with the same Lipschitz constant, and it is easy to see that a legal loop in the tension graph for the original map becomes a legal loop for the folded map in the tension graph. See~\cite{FM13} for more details. 
  
  In particular,
  $\Min(\phi)$ is   invariant by folds directed by weakly optimal maps.
\end{rem}

\begin{defn}[Reducibility]\label{defn_reducibility}

  An  automorphism $\phi$ is called {\em reducible} if there is an $X \in \O(\Gamma)$ 
  and $f:X\to X$
  representing $\phi$ having a proper, non-trivial, $f$-invariant $\Gamma$-sub-forest
  (See Definition~\ref{def2020_2.12}).

  We say $\phi$ is {\em irreducible} if it is not reducible.
  
%
%
%
%
\end{defn}
\begin{rem}
In the  connected case, if $G=F_n$ then this definition coincides with the usual definition of irreducibility.  For irreducible
automorphisms we have $\Min(\phi)\neq\emptyset$, but the converse is not true in
general. (See~\cite{FM13} for more details.)

\end{rem}

%

\begin{rem}
If $\phi$ is irreducible, then any closed simplex has a min-point for $\lambda$. (See for
instance~\cite[Section~$8$]{FM13}. See also Proposition~\ref{propfinvariance} below.)
In~\cite{BestvinaBers,FM13}  automorphisms so that $\Min(\phi)\neq\emptyset$ and $\lambda>1$ 
are called hyperbolic.
\end{rem}

\begin{defn}[Train track between trees]
  Let $\sim$ be a gate structure on a (not necessarily connected) tree $X$. A map
  $f:X\to X$ is a {\bf train track map w.r.t. $\sim$} if
  \begin{enumerate}
  \item any vertex has at least two gates w.r.t. $\sim$;
  \item $f$ maps edges to legal paths (in particular, $f$ does not collapse edges);
  \item for any vertex $v$, if $f(v)$ is a vertex, then $f$ maps inequivalent germs at $v$ to
    inequivalent germs at $f(v)$.
  \end{enumerate}
\end{defn}

We already defined the gate structure $\sim_f$ induced by a straight map (Definition~\ref{pr4_gate}).
\begin{defn}[Gate structure $\langle\sim_{f^k}\rangle$]
  Let $X$ be a (not necessarily connected) tree, and let $f:X\to X$ be a map whose components
  are straight. We define the gate structure $\langle \sim_{f^k}\rangle$ as the equivalence relation on germs generated
  by all $\sim_{f^k}, \ k\in \N$.
\end{defn}

\begin{lem}\label{Lemmatt1}
  Let $\phi\in\Aut(\Gamma)$, $X\in\O(\Gamma)$ and $\sim$ be a gate structure on $X$. Let
  $f:X\to X$ be a straight map representing $\phi$. If
  $f:X\to X$ is a train track map w.r.t. $\sim$, then relation $\sim$ is stronger than (\em{i.e.} it
  contains) $\simfk$. In particular
  if $f$ is a train track map w.r.t. some $\sim$ then it is a train track map w.r.t $\simfk$.
\end{lem}
See~\cite[Section~8]{FM13} for a proof (where it is proved  in the connected case, but
connectedness plays no role).

Now we give a definition of partial train track map
representing an automorphism. Our definition is given at once for both reducible and
irreducible automorphisms. In the irreducible case coincides with the standard one. For reducible
automorphisms there already exist notions of relative and absolute train tracks (see~\cite{BestvinaHandel}). Our notion
is different from that of relative train tracks; absolute train tracks are train tracks in
our setting but not vice versa.\footnote{Our present
definition of partial train track map coincides with the notion of {\em optimal} train track
map given in~\cite{FM13} for irreducible automorphisms in the connected case.}

The main motivation for this new definition is that it 
well-behaves with respect to the displacement function, as we will see that it characterise minimally displaced points.

\begin{defn}[Partial train track maps for automorphisms]\label{defttm}
  Let $[\phi] \in\Out(\Gamma)$. Let $X\in\O(\Gamma)$ and
  let $f:X\to X$   be a straight map representing
  $\phi$. Then we say that $f$ is a
   \begin{itemize}
  \item {\bf partial train track map with one-step gates} if there is a (not necessarily
    proper) $f$-invariant $\Gamma$-sub-forest 
  $A\subseteq X_{\max}(f)$ such that
  \begin{enumerate}
  \item  $f|_A$ is a train track map w.r.t. $\sim_f$, and 
  \item  $A$ is homotopically non-trivial; that is, $A$ contains the axis of a hyperbolic element.
  \end{enumerate}  
  \item {\bf partial train track map} if there is a (not necessarily proper) $f$-invariant, homotopically non-trivial, $G$-sub-forest
  $A\subseteq X_{\max}(f)$ such that $f|_A$ is a train track map w.r.t. $\simfk$. 
  \end{itemize}
At level of graphs, a map $X\to X$ is said a partial train track if it is the projection of a
partial train track map $\wt X\to \wt X$ (that is, if there is a non-trivial invariant sub-graph
of the tension graph so that $f$ restricted to that graph is a train-track map). We stress the
fact that no requirements are made outside the tension graph.
\end{defn}

Here some more remarks are needed, since the {\em metric} theory of train tracks maps, first introduced
in~\cite{BestvinaHandel}, does not have a completely standard treatment. 
That is, train tracks can be defined topologically, from a simplicial viewpoint, and a metric is subsequently introduced. 
Our point of view is to always have a metric, and deduce the topological properties from certain minimizing conditions. Additionally, it should be noted that the standard definition requires train track maps (or representatives in general) to send vertices to vertices, whereas we do not. While this condition is extremely useful, and can often be recovered, our arguments are based on continuous deformations where it is more natural to relax this condition. These are sometimes called {\em simplicial} train-tracks and are useful 
for computation purposes. This is not a big issue as the closure of any simplex containing a
partial train track also contains a simplicial one. (See~\cite{FM13}.)

In the case that $\phi$ is irreducible there is not much difference between topological and metric train track maps. Indeed if $f:X\to X$ is a topological
train track map representing $\phi$, then one can rescale the edge-lengths of $X$ so that $f$
is a train track map for Definition~\ref{defttm}. And the same holds true if
$f$ has no proper invariant sub-graphs. This is because train track maps do not collapse edges,
hence edge-lengths can be adjusted so that every edge is stretched by the same amount. In particular, the
following two results are proved in~\cite{FM13} for irreducible automorphisms and $\Gamma$
connected. The proofs for general automorphisms are essentially the same (details are left to the reader).

\begin{lem}\label{Lemmatt2}
  Let $[\phi]\in\Out(\Gamma)$, $X\in\O(\Gamma)$, and $f:X\to X$ be a straight map representing
  $\phi$. Then $f$ is partial train track if and only if there is an embedded periodic line $L$
  in $\widetilde{X}_{\max}$  such that $f^k(L)\subseteq \widetilde{X}_{\max}$  and $f^k|_L$ is
  injective for all $k\in\N$ (here $A=\cup_kf^k(L)$ is the invariant sub-forest). In particular if $f$ is partial train track then
  \begin{enumerate}
  \item $f^k$ is a partial train track;
  \item $\Lip(f)=\Lambda(X,\phi X)$ (hence $f$ is weakly optimal);
  \item $\Lip(f)^k=\Lip(f^k)=\Lambda(X,\phi^kX)$.
  \end{enumerate}
\end{lem}

\begin{cor}
    Let $\phi\in\Aut(\Gamma)$, $X\in\O(\Gamma)$, and $f:X\to X$ be a map
    representing  $\phi$. Suppose that there is an embedded periodic line $L$ in $\widetilde{X}$ such that
  $f^k|_L$ is injective for all $k\in\N$.  Suppose moreover that $\cup_kf^{k}(L)=\widetilde{X}$.
  Then there is $X'$ obtained by rescaling edge-lengths of $X$ so that $\PL(f):X'\to X'$ is a
  train track map.
\end{cor}

In general, if $\cup_kf^k(L)$ is just an $f$-invariant subtree $Y$ of $X$, we can adjust edge
lengths so that every edge of $Y$ is stretched the same, but we cannot guarantee a priori that
$Y\subset X_{\max}$. 

\begin{defn}[Train track sets]
For any $[\phi]\in\Out(\Gamma)$ we define:
$$\TT(\phi)=\{X\in\O(\Gamma):\exists f:X\to X \text{ partial train track}\}$$
$$\TTo(\phi)=\{X\in\O(\Gamma):\exists f:X\to X \text{ partial train track with one-step gates}\}$$
If we need to specify the map we write $(X,f)\in \TT(\phi)$ or $(X,f)\in\TTo(\phi)$.\footnote{We remark that, since in the irreducible case our present definition of train track
  map corresponds to that of {\em optimal} train track map of~\cite{FM13}, the two definitions of
  $\TT$ and $\TTo$ coincide with those given in~\cite{FM13}.}

\end{defn}

\begin{thm}\label{Theoremtt}
Let $[\phi]\in\Out(\Gamma)$. Then $$\overline{\TTo(\phi)}=\TT(\phi)=\Min(\phi)=\LocMin(\phi)$$
where the closure is made with respect to the simplicial topology.
\end{thm}
\proof If $\phi$ is irreducible and $\Gamma$ connected, the proof is given in~\cite{FM13} and goes through the
following steps:
\begin{enumerate}
\item $\TTo(\phi)\subseteq \TT(\phi)\subseteq\Min(\phi)\subseteq \LocMin(\phi)$.
\item If $X$ locally minimizes $\lambda_\phi$ in $\Delta_X$, and $f:X \to X$ is an optimal map
  representing $\phi$ then $X_{\max}$ contains a homotopically non-trivial $f$-invariant sub-forest, $A$.
\item $\TTo(\phi)$ is dense in $\LocMin(\phi)$.
\item $\TT(\phi)$ is closed.
\end{enumerate}
We now adapt the proof so that it works also for  $\phi$ reducible and general $\Gamma$. Clearly $\Min(\phi)\subseteq \LocMin(\phi)$. By Lemma~\ref{Lemmatt1} $\TTo(\phi)\subseteq \TT(\phi)$. We
now show that $\TT(\phi)\subseteq\Min(\phi)$, arguing by contradiction. If $X\in\TT(\phi)$ and $\lambda_\phi(X)>\lambda(\phi)$ then
there is $Y\in\O(\Gamma)$ such that $\lambda_\phi(Y)<\lambda_\phi(X)$. By Lemma~\ref{Lemmatt2}
$\Lambda(X,\phi^kX)=\lambda_\phi(X)^k$ but then
$$\lambda_\phi(X)^k=\Lambda(X,\phi^kX)\leq\Lambda(X,Y)\Lambda(Y,\phi^kY)\Lambda(\phi^kY,\phi^kX)$$
$$=\Lambda(X,Y)\Lambda(Y,\phi^kY)\Lambda(\phi^kY,\phi^kX)
\leq\Lambda(X,Y)\Lambda(Y,X)\lambda_\phi(Y)^k$$
thus $(\frac{\lambda_\phi(X)}{\lambda_\phi(Y)})^k$ is bounded for any $k$, which is impossible if
$\frac{\lambda_\phi(X)}{\lambda_\phi(Y)}>1$.

Thus we have
$$\TTo(\phi)\subseteq\TT(\phi)\subseteq\Min(\phi)\subseteq\LocMin(\phi).$$

\begin{lem}\label{Lemmatt3}
  Suppose $(X,f)$ locally minimizes $\lambda_\phi$ in $\Delta_X$. Then there is a hmotopically non-trivial $ A
  \subseteq   X_{\max}$ which is $f$-invariant.
\end{lem}
\proof For every $\epsilon >0 $, consider the $\epsilon$-neighbourhood of $X$ in
$\Delta_{X}$. For each point in this neighbourhood, $f$ induces a map on it via rescaling. We
optimize that map by using Theorem~\ref{Lemma_opt}, and consider the tension graph,
$A_{\epsilon}$ with respect to that optimal map. By abuse of notation, we think of
$A_{\epsilon}$ as a subgraph of $X$ (since all we have done is rescale edges). $A_\epsilon$ is
homotopically non-trivial because of Proposition~\ref{propELr}.

Now for each (sufficiently small) $\epsilon$, choose a particular  $X_{\epsilon} $   in the $\epsilon$-neighbourhood of $X$,  such that 
\begin{itemize}
	\item $X_{\epsilon}$ minimizes $\lambda_\phi$ in the  $\epsilon$-neighbourhood of $X$ in $\Delta_{X}$ (we allow that $X_{\epsilon}$ could be $X$ and in particular, we have that $\lambda_\phi(X) = \lambda_\phi(X_{\epsilon})$), and 
	\item the tension graph, $A_{\epsilon}$ is smallest, with respect to inclusion,
          amongest all possible choices, subject to the previous condition.
\end{itemize}

In particular, since there are only finitely many subgraphs, by taking sufficiently small $\epsilon$ we may assume that  $A:=A_{\epsilon}$ does not depend on $\epsilon$.

Let $g_{\epsilon}$ denote the optimal map on $X_{\epsilon}$ (obtained as above) and
$f_{\epsilon}$ denote the map on $X$ obtained by rescaling $g_{\epsilon}$. (That is,
$g_{\epsilon}$ is simply the optimization --- via Theorem~\ref{Lemma_opt} --- of $f$, when thought of as a map on $X_{\epsilon}$, and $f_{\epsilon}$ is $g_{\epsilon}$, thought of as a map on $X$.)

Then it is clear, by Theorem~\ref{Lemma_opt},  that $\lim_{\epsilon \to 0} d_{\infty}( f_{\epsilon},f)=0$ and hence $f_{\epsilon} \to f$, uniformly. 

If $A_{\epsilon}$ (thought of as a subgraph of $X_{\epsilon}$) contains an edge $e$
whose image (under $g_{\epsilon}$) is not in $A_{\epsilon}$, then by shrinking (the orbit of)
such an edge, either 
we reduce  $\lambda_\phi(X) = \lambda_\phi(X_{\epsilon})$ --- which is impossible --- or we
reduce the tension graph --- which is also 
impossible. Thus $A=A_{\epsilon}$ is $g_{\epsilon}$ and hence $f_{\epsilon}$ invariant. 

But now the fact
that $f_{\epsilon} \to f$, implies that $A$ is $f$-invariant.  Moreoever, by
Theorem~\ref{Lemma_opt}, for sufficiently small $\epsilon$, $A_{\epsilon}$ --- hence $A$ ---
will be a subgraph of $X_{\max}$  since $d_{\infty}(g_\e,f)\to 0$.


\qed

\begin{lem}\label{Lemmatt4}
  $\LocMin(\phi)\subseteq \overline{\TTo(\phi)}$. More precisely, let $X\in\O(\Gamma)$ and fix
  $f:X\to X$ an optimal map representing $\phi$. Suppose $X$ has an open
  neighbourhood $U$ such that for any $Y\in U$ obtained from $X$ by a sequence of simple folds
  directed by $f$, we have $\lambda_\phi(X)\leq\lambda_\phi(Y)$. Then there
  is a sequence $Y_n\in U$, all contained in the same simplex, with $Y_n\to X$ and $\wt Y_n\in
  \TTo$, each equipped with a partial train track map, $f_n$ such that $f_n \to f$ uniformly.
\end{lem}
\proof The proof is basically the same as in~\cite{FM13}.
When $Y$ obtained from $X$ by folds directed by $f$, then we
let $f_Y$ denote the induced optimal map. First we remark that if $Y$ is obtained from $X$ by
folds directed by $f$ then $\lambda_\phi(Y)\leq\lambda_\phi(X)$ and by minimality of $X$ we have $\lambda_\phi(Y)=\lambda_\phi(X)$.
We consider the gate structure induced by $f_Y$.
We call a vertex of $Y_{\max}$ {\em foldable} if it
has at least two edges of $Y_{\max}$ in the same gate.

Locally, by using arbitrarily small folds
in $X_{\max}$, directed by $f$, we find $Y\in U$ such that
\begin{enumerate}
\item $\lambda_\phi(Y)=\lambda_\phi(X)$;
\item the simplex $\Delta=\Delta_Y$ maximizes the dimension among simplices reachable from $X$ via
  folds directed by $f$;
\item $Y$ minimizes $Y_{\max}$ among points of $\Delta$ satisfying $(1)$;
\item $Y$ maximizes the number of orbit of foldable vertices of $Y_{\max}$ among points of
  $\Delta$ satisfying $(1),(3)$.
\end{enumerate}
Let $A\subseteq {Y_{\max}}$ be a homotopically non-trivial $f_Y$-invariant sub-forest given by Lemma~\ref{Lemmatt3}. We
claim that ${f_Y}|_A$ is a train track map with one-step gates.
Indeed, otherwise there is either an edge $e$ or a legal turn $\tau$ in $A$ having
       illegal image. Let $v$ be the vertex of $\tau$.
       \begin{itemize}
       \item If $f_Y(e)$ contains an illegal turn $\eta$ then by folding (the orbit of) it a
         little, we would 
         reduce the tension graph, contradicting $(3)$. (Note that $\eta\subset Y_{\max}$
         because $A\subseteq Y_{\max}$ is $f_Y$-invariant, thus by folding $\eta$ we do not
         change simplex of $\O(\Gamma)$ thanks to $(2)$.)
       \item If $f_Y(\tau)$ is an illegal turn $\eta$ then we fold it a little.  Either
         $Y_{\max}$ becomes one-gated at $v$, and in this case the optimization process reduces
         the tension graph, contradicting $(3)$, or $v$ was not foldable at $Y$ and becomes
         foldable, thus contradicting $(4)$.
       \end{itemize}
       Finally, note that given such an $Y$, the sequence $Y_n$ can be chosen in $\Delta_Y$.
       \qed

In particular, since $\Min(\phi)$ is clearly closed, we now have:
$$\LocMin(\phi)\subseteq
\overline{\TTo(\phi)}\subseteq\overline{\TT(\phi)}\subseteq\overline{\Min(\phi)}=\Min(\phi)\subseteq\LocMin(\phi)$$
hence all inclusions are equalities.

\begin{lem}
$\overline{\TT(\phi)}=\TT(\phi)$.
\end{lem}
\proof Let $X\in\overline{\TT(\phi)}=\Min(\phi)$. Let $f:X\to X$ be an optimal map representing
$\phi$. By Lemma~\ref{Lemmatt4} there is $Y_n\to X$ and $f_n\to f$ so that
$(Y_n, f_n)\in\TTo(\phi)$. By Lemma~\ref{Lemmatt2} there is an embedded periodic line
$L_n$ in $( Y_n)_{\max}$ such that $f_n^k(L_n)\subset ( Y_n)_{\max} $ is embedded for all $k\in\
N$. 

We will argue that, up to taking subsequences, there is a single embedded periodic line $L$ which is $f_n$-legal for all $n$. (More precisely, there is a single hyperbolic element whose axis in each $Y_n$ is both legal and contained within the tension graph of $Y_n$. We can think of this as a single topological line, since the $Y_n$ all belong to the same simplex.)

First of all, we may assume that each $f_n$ has the same (topological) tension graph. Next we claim that, for some $n$, the line $L_n$ is $f_m$-legal for infinitely many $m$.

Suppose this is not the case. Then for each $n$, there is a turn crossed by $L_n$ which is $f_m$-illegal for infinitely many $m$; this is because each $L_n$ only crosses finitely many orbits of turns, so one of these orbits must be illegal for infinitely many $m$ (recall that gate structures are equivariant) since otherwise $L_n$ would be $f_m$-legal for infinitely many $m$.

Starting with $Y_1$, we can now define a subsequence $Y_n$ as follows: Let $\tau_1$ be a turn
in $L_1$ which is illegal for infinitely many $m$. Now choose the first $Y_2$ where $\tau_1$ is
$f_2$-illegal and let $\tau_2$ be a turn crossed by $L_2$ such that there are infinitely many
$m$ where both $\tau_1$ and $\tau_2$ are $f_m$-illegal. (If there were no such turn, then
amongst the infinitely many $Y_m$ where $\tau_1$ is illegal, we would get infinitely many $m$
in which $L_2$ is $f_m$-legal and we would be done.) We may continue in this fashion to obtain (after renumbering) a sequence $Y_n$ such that $\tau_n$ is a turn crossed by $L_n$ (hence is $f_n$-legal) and  all the turns $\tau_1, \ldots, \tau_{n-1}$ are $f_n$-illegal.

By the equivariance of the $f_n$, the $\tau_n$ must all be in distinct orbits. If there were infinitely many of the $\tau_n$ we should be able to find a non-free vertex, $v$, and two edges $e_1, e_2$ so that the turn defined by $e_1, e_2$ is in the same orbit as some $\tau_r$ (for some $r$) and the turn defined by $g e_1, he_2$ is in the orbit of some $\tau_s$ (for some $s$), where $g, h$ are group elements stabilising $v$ and $g \neq h$. (That is, if we project to the quotient graph, then $\tau_r$ and $\tau_s$ look the same. The fact that they are in different orbits now implies that they are based at a non-free vertex as described.)

But now for some (almost all) $Y_m$, both $\tau_r$ and $\tau_s$ are $f_m$-illegal, which contradicts the fact that $f_m$ is equivariant and the action on the edge set is free. Therefore the process of choosing our subsequence must terminate at some point, at which time we will have produced a single embedded periodic line $L$ which is $f_m$-legal for infinitely many $m$. Without loss of generality, we may assume that $L$ is $f_n$-legal for all $n$.

Since $f_n\to f$ and the maps are all straight, $L\subset  X_{\max}$
and  $f^k(L)\subset  X_{\max}$.
Moreover, if $f^k$ were not injective on $L$ for some $k$, then we
could find $\varepsilon >0$ and points $p,q$ with $d_X(p,q)=\varepsilon$ and
$f^k(p)=f^k(q)$. Now the fact that $f_n\to f$ would contradict the fact that $f_n^k|_L$ is a
homothety of ratio $\lambda(\phi)$. Thus $f^k|_L$ is embedded for any $k$, $f$ is a partial train
track map and so $X\in\TT(\phi)$.\qed

This completes the proof of Theorem~\ref{Theoremtt}.\qed


We end this section by proving a lemma which is basically a rephrasing of
Lemma~\ref{Lemmatt4} with a language which will be more usable. (For example we will use it in forthcoming part II of the present paper.)

\begin{defn}[Exit points]\label{exitp}
    Let $[\phi] \in\Out(\Gamma)$. A point $X\in\O(\Gamma)$ is called an {\em exit point} of
    $\Delta_X$ if for any neighbourhood $U$ of $X$ in $\O(\Gamma)$, there exists an optimal map $f: X \to X$, representing $\phi$, a point $X_E\in U$,  and a folding path (Definition~\ref{foldpath}) directed by $f$, $X=X_0,X_1,\dots,X_m=X_E$ in $U$, such that $\Delta_{X_i}$ is finitary face
     of $\Delta_{X_{i+1}}$, $\Delta_{X}$ is a proper face of
    $\Delta_{X_E}$, and     such that $$\lambda_\phi(X_E)<\lambda_\phi(X)$$ (strict inequality).
\end{defn}

\begin{rem*}
	The idea is that an exit point is one which allows one to decrease the displacement by arbitrarily small folds. The complication is that a simple fold may not be sufficient, so we allow folding paths. 
	
	The principal application of this is as below; when $X$ does not admit a partial train track representing $\phi$, but is minimally displaced within its simplex, then it will be an exit point. This idea of folding to decrease the displacement goes back to \cite{BestvinaHandel}. 
\end{rem*}

\begin{lem}\label{LemmaX}
  Let $[\phi] \in\Out(\Gamma)$ and $X\in\O(\Gamma)$ such that
  $\lambda_\phi(X)$ is a local minimum for $\lambda_\phi$ in $\Delta_X$.
  Suppose $X\notin\TT(\phi)$.

Then, for any open neighbourhood $U$ of $X$ in $\O(\Gamma)$, there exists an optimal map $f:X \to X$, representing $\phi$, points $Z, X' \in U$, and a folding path, $X=X_0, \ldots, X_m=Z, X_{m+1}, \ldots, X_n=X'$, directed by $f$ and such that: 

\begin{itemize}
	\item $X_0, \ldots, X_m \in U \cap \Delta_X$,
	\item $\lambda_{\phi}(Z) = \lambda_{\phi}(X)$,
	\item $\Delta_X$ is a proper face of $\Delta_{X'}$,
	\item $\lambda_\phi(X')<\lambda_\phi(X)$.
\end{itemize}

%
  
  (See Figure~\ref{fig:lemmaX}.)


In particular $X$ is an exit point of $\Delta_X$.
\end{lem}
\setlength{\unitlength}{1ex}
\begin{figure}[htbp]
  \centering
  \begin{picture}(52,17)
    \put(10,2){\line(1,0){30}}
    \put(10,2){\line(1,1){15}}
    \put(40,2){\line(-1,1){15}}
    \put(30,2){\line(0,1){3}}
    \put(25,2){\makebox(0,0){$\bullet$}}
    \put(30,2){\makebox(0,0){$\bullet$}}
    \put(30,5){\makebox(0,0){$\bullet$}}
    \put(18,2){\makebox(0,0){$($}}
    \put(32,2){\makebox(0,0){$)$}}
    \put(21,0){\makebox(0,0){$U \cap \Delta_X$}}
    \put(25,0){\makebox(0,0){$X$}}
    \put(30,0){\makebox(0,0){$Z$}}
    \put(30,7){\makebox(0,0){$X'$}}
    \put(22,10){\makebox(0,0){$\Delta_{X'}$}}
    \put(8,2){\makebox(0,0){$\Delta_{X}$}}
  \end{picture}
  \caption{Graphical statement of Lemma~\ref{LemmaX}}
  \label{fig:lemmaX}
\end{figure}
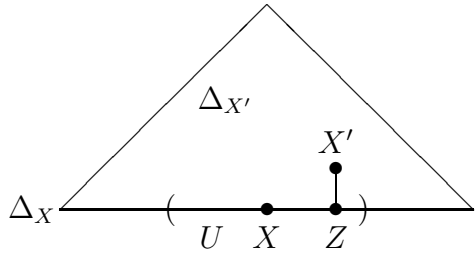
\proof Let's prove the first claim. Since $X\notin \TT(\phi)$, by Theorem~\ref{Theoremtt} there is a neighbourhood of $X$
in $\Delta_X$ which is contained in the complement of $\TTo(\phi)$. Without loss generality we may assume
that $U \cap \Delta_X$ is contained in such neighbourhood.

Let $f:X\to X$ be an optimal map representing $\phi$. If there is a non-trivalent foldable
vertex in $X_{\max}$ then we set $Z=X$ and we are done. Otherwise, consider $Z\in U \cap \Delta_X$ obtained
from $X$ by a fold directed by $f$ (we still denote by $f:Z\to Z$ the map induced by $f$). We have $\lambda(Z)\leq \lambda(X)$. Since
$\lambda(X)$ is a local minimum in $\Delta_{X}$, we must have $\lambda(Z)=\lambda(X)$. Let $\ul A\subset \ul{Z_{\max}}$ be
an $f$-invariant sub-graph given by Lemma~\ref{Lemmatt3}. Since $Z\notin \TTo(\phi)$, the
restriction $f|_A$ is not a train-track with one-step gates. That is, $f:Z \to Z$ is not a partial train track map with one-step gates.

It follows that by using
folds directed by optimal maps we can either
\begin{itemize}
\item[$a)$] reduce the tension graph; or
\item[$b)$] increase the number of foldable vertices; or
\item[$c)$] create a non-trivalent foldable vertex.
\end{itemize}
So far $Z$ is generic. We choose $Z\in U \cap \Delta_X$ so that, in order:
\begin{enumerate}
\item it locally minimizes the tension graph;
\item it locally maximizes the number of foldable vertices among points satisfying $(1)$.
\end{enumerate}

For such a $Z$ the only possibility that remains in the above list of alternatives is $c)$, which therefore admits a fold directed by an optimal map into a simplex of strictly larger dimension. 

We continue this process inductively, and we get the result due to the fact that the simplicial dimension of $\O(\Gamma)$ is bounded.

\qed

\section{Behaviour of $\lambda$ at bordification points}\label{Sectionbordification}
For the rest of the section we fix $G,\G$ and
$\Gamma=\sqcup_i\Gamma_i$ as in Definitions\ref{def20202.6} and~\ref{not:gamma}. We also fix $\phi\in\Aut(\Gamma)$ and
if there is no ambiguity we understand that $$\lambda=\lambda_\phi.$$ In this section we discuss the behaviour of $\lambda$
at boundary points of outer space, that is to say, when we reach points in
$\partial_\infty\O(\Gamma)$. As above, we remind the reader that the results of this section hold true
in particular for $CV_n$ and its simplicial bordification.

We will see that the function $\lambda$
is not continuous and we will provide conditions that assure continuity along particular
sequences. We will also focus on the behaviour of $\lambda$ on
horoballs. In this section we will often work with $\Gamma$-graphs. We recall
that we are denoting by $\ul X$ the $\Gamma$-graph corresponding to $X\in\O(\Gamma)$.

Points near the boundary at infinity have some sub-graph that is almost collapsed. This is
usually referred to as the ``thin'' part of outer space. We introduce now more quantified notions
of ``thinness''.

\begin{defn}[$\e$-thinness]
  Let $\e>0$. A point $X\in\O(\Gamma)$ is {\em $\e$-thin} if there is a nontrivial loop
  $\gamma$ in $\ul X$ such that $L_X(\gamma)<\e\vol(X)$.
\end{defn}

\begin{defn}[$(M,\e)$-collapsed points]
  Let $M,\e>0$. A point $X\in\O(\Gamma)$ is {\em $(M,\e)$-collapsed} if there is a non-trivial loop
  $\gamma$ in $\ul X$ such that $L_X(\gamma)<\e\vol(X)$ and for any other loop $\eta$
  such that 
  $L_X(\eta)\geq\e\vol(X)$ we have $L_X(\eta)>M\vol(X)$.
\end{defn}

\begin{defn}[$\e$-thin part]
  Let $\e>0$. For any $X\in\O(\Gamma)$ we define $X_\e$ the $\e$-thin part of $X$ as the 
  sub-forest formed by the union of the axes of elements $\gamma$ with
  $L_X(\gamma)<\e\vol(X)$. (Note that $\ul{X_\e}$ is a core graph.)
\end{defn}

\begin{defn}[$\phi$-invariance]
  Let $X\in \O(\Gamma)$. A $\Gamma$-sub-forest $A\subset X$ is called {\em $\phi$-invariant} if there is a
  straight map $f:X\to X$ representing $\phi$ such that $f(A)\subseteq A$. 
  The quotient, $\ul A$ is called a $\phi$-invariant subgraph of $\ul X$.
  
\end{defn}

We now state some easy facts, the first of which can be found in \cite{BestvinaBers}.

\begin{prop}
  For any $C>\lambda(\phi)$ there is $\e>0$ such that for any $X\in\O(\Gamma)$, if
  $\lambda_\phi(X)<C$ and  $X_\e\neq\emptyset$ then $\ul X$  contains a non-trivial\footnote{In the
    sense of Definition~\ref{defn_reducibility}.} $\phi$-invariant
  subgraph. 
\end{prop}
For a proof in the case $\Gamma$ is connected see~\cite[Section 8]{FM13} (connectedness plays
in fact no role).

However, we will need a slightly more precise statement, in order to be able to determine a particular invariant subgraph.

\begin{prop}\label{propfinvariance}
  Let $C \geq 1$ and $M>0$. Let $D$ be the maximal number of orbits of edges for any tree in  $\O(\Gamma)$.

\noindent  
  Let $\e = 1/2 \min \{ M/CD, 1/D \}$.   Then, for  $X\in\O(\Gamma)$, if $\lambda_\phi(X)<C$
  and $X$ 
  is $(M,\e)$-collapsed, then $X_\e$ is not the whole $X$ and it is $\phi$-invariant.
\end{prop}
\proof By definition any edge in $X_\e$ is shorter than $\e\vol(X)$. Thus we have $\vol(X_\e)<\e\vol(X)D$. In particular, since $\e D<1$ then
$X_\e\neq X$ (and thus there exists a loop $\eta$ with $L_X(\eta)>\e\vol(X)$, whence $L_X(\eta)>M\vol(X)$), since $X$ is $(M,\e)$-collapsed).

Let $f:X\to X$ be an optimal straight map representing $\phi$.  By picking a maximal tree in
$\ul X$, we may find a generating set of the fundamental group of (each component of) $\ul X_\e$ whose elements have length at most $2 \vol(X_\e)$. For any such generator, $\gamma$, we have that $L_X(f(\gamma))/L_X(\gamma) \leq C$ and hence, $L_X(f(\gamma)) \leq C L_X(\gamma) \leq 2C \vol(X_\e) < 2C D \e \vol(X) \leq M \vol(X)$. But since $X$ is $(M,\e)$-collapsed, we get that $L_X(f(\gamma)) < \e \vol(X)$. Hence $f(\gamma)$ is homotopic to a loop in $X_\e$. 

Varying $\gamma$ we deduce that $X_\e$ is $\phi$-invariant. 
\qed

\begin{prop}\label{oss}
  Let $X\in\Og(\Gamma)$ and
  $\phi\in\Aut(\Gamma)$. Suppose that $A\subset X$ is a 
  $\phi$-invariant core graph. Then $\lambda_{\phi|_A}(A)\leq \lambda_\phi(X)$.
\end{prop}
\proof Let $f: X\to X$ be a straight map representing $\phi$. Since $A$ is $\phi$-invariant,
$f(A)\subset A$ up to homotopy. By passing to the universal covering we see that $f|_A:A\to
X$ retracts to a map $f_A:A\to A$ representing $\phi$ with $\Lip(f_A)\leq\Lip(f)$, hence
$\lambda_{\phi|_A}(A)\leq\Lip(f_A)\leq \Lip(f)=\lambda_\phi(X)$.\qed

\begin{thm}[Lower semicontinuity of $\lambda$]\label{fatto1}
  Fix $\phi\in\Aut(\Gamma)$ and $X\in\Og(\Gamma)$. Let $(X_i)_{i\in\N}\subset \Delta_X$ be a
  sequence such that there is $C$ such that for any $i$,  $\lambda_\phi(X_i)<C$. Suppose that $X_i\to
  X_\infty\in\partial_\infty\Delta_X$ which is  obtained from $X$  by collapsing a sub-graph
  $A\subset X$. Then $\phi$ induces an element of $\Aut(X/A)$, still denoted by
  $\phi$.

  Moreover $\lambda_\phi(X_\infty)\leq \liminf_{i\to\infty} \lambda_\phi(X_i)$, and
  if strict inequality holds, then there is a sequence of
  minimal optimal maps $f_i:X_i\to X_i$ representing $\phi$ such that eventually on $i$ we have
  $(X_i)_{\max}\subseteq \core(A)$.
\end{thm}
\proof Let $M$ be the ``systole'' of $X_\infty$, that is to say the shortest length
of simple non-trivial loops in $X_\infty$. For any $M/\vol(X)>\e>0$, eventually on $i$, $X_i$ is
$(M/2\vol(X),\e)$-collapsed and $(X_i)_\e=\core(A)$. By Proposition~\ref{propfinvariance} $A$ is
$\phi$-invariant, thus $\phi\in\Aut(X/A)$.

For any loop $\gamma$ the lengths $L_{X_i}(\gamma)$ and $L_{X_i}(\phi(\gamma))$ converge to
$L_{X_\infty}(\gamma)$ and $L_{X_\infty}(\phi(\gamma))$ respectively.  Therefore, if $\gamma$
is a candidate in $X_\infty$ that realizes $\lambda_\phi(X_\infty)$, we have that
$\lambda_\phi(X_i)\geq L_{X_i}(\phi(\gamma))/L_{X_i}(\gamma)\to \lambda_\phi(X_\infty)$ whence
the lower semicontinuity of $\lambda$.

On the other hand, by Theorems~\ref{thmminopt} and~\ref{Lemma_opt}, for any $i$ there is a minimal optimal map
$f_i:X_i\to X_i$ representing $\phi$. Let $\gamma_i$ be a candidate that realizes
$\lambda_\phi(X_i)$, i.e. a $f_i$-legal candidate in $(X_i)_{\max}$. Since $X$ is
combinatorically finite, we may assume w.l.o.g. that $\gamma_i=\gamma$ is the same loop for any
$i$. We have
$$\lambda_\phi(X_i)=\frac{L_{X_i}(\phi(\gamma))}{L_{X_i}(\gamma)}\to \frac{L_{X_\infty}(\phi(\gamma))}{L_{X_\infty}(\gamma)}$$

Thus if $L_{X_\infty}(\gamma)\neq 0$ we have
$\lambda_\phi(X_\infty)=\liminf\lambda_\phi(X_i)$. It follows that if there is a jump in $\lambda$
at $X_\infty$, then any legal candidate is contained in $A$. Since $f_i$ is minimal this
implies that $\core(A)$ contains the whole tension graph.\qed

\begin{rem}
 A comment on Theorem~\ref{fatto1} is required. To avoid cumbersome notation, we have decided to denote by
 $\phi$ both the element of $\Aut(X)$ and the one induced in $\Aut(X/A)$. So when we write
 $\lambda_\phi(X_\infty)$ we mean $\Lambda(X_\infty,\phi X_\infty)$ as elements in $\O(X/A)$. In
 particular, $\lambda_\phi=\inf_X\lambda_\phi(X)$ can be different if computed in $\O(X)$ or in
 $\O(X/A)$. When this will be crucial we will specify in which space we take the infimum.

Moreover, if $\phi|_A$ is the restriction of $\phi$ to $A$, then $\lambda_{\phi|_A}$ is
 calculated in the space $\O(A)$. While the simplex $\Delta_{X_\infty}$ is a simplicial face
 of $\Delta_X$, $\Delta_A\in\O(A)$ does not have the same meaning. One could argue that $\Delta_A$ is the
 simplex ``opposite'' to $\Delta_{X_\infty}$ in $\Delta$, but $\phi$ does not necessarily
 produce an element of $\Aut(X/(X\setminus A))$ as the complement of $A$ may be not invariant.
\end{rem}

Clearly, if $A\subset X$ is $\phi$-invariant then $\lambda_\phi(X/A)<\infty$. On the other
hand, if $A$ is not $\phi$-invariant, its collapse makes $\lambda$ explode. Thus we can extend
the function $\lambda$ as follows.

\begin{defn}\label{dli}
  Let $X_\infty\in\partial_\infty\Og(\Gamma)$. We say that $\lambda_\phi(X_\infty)=\infty$ if
  $X_\infty$ is obtained from a $\Gamma$-graph $X$ by collapsing a sub-graph $A\subset X$
  which is not $\phi$-invariant. (Note that $X_\infty$ since is not in $\O(\Gamma)$, then $A$
  must have some non-trivial component).
\end{defn}

In general, the function $\lambda$ is not uniformly continuous with respect to the Euclidean metric,
even in region where it is bounded, and so we cannot extend it to the simplicial closure of
simplices. However we see now that the behaviour of $\lambda$ is controlled on segments.

We recall the description of horoballs given in~\ref{sechor}. Suppose that $X_\infty$ is
obtained from a $\Gamma$-graph $X$ by
collapsing a $\phi$-invariant core sub-graph $A=\cup_i A_i$. Let $k_i$ be the number of
germs of edges incidents to $A_i$ in $X\setminus A$. Then $\Hor(X_\infty)$
is a product of outer spaces with marked points $\O(A_i,k_i)$. 
\begin{nt}
We denote
$\pi:\Hor(X_\infty)\to\mathbb P \O(A)$ the projection that forgets marked points.  
\end{nt}
 Note that we
chosen $X_\infty$ to not be projectivized and $\mathbb P\O(A)$ to be projectivized.
For any
$Y\in\mathbb P\O(A)$ if
$Z\in\pi^{-1}(Y)$, then there is a scaled copy of $Y$ in $Z$. We denote by $\vol_Z(Y)$ the
volume of $Y$ in $Z$. With this notation in place, we can now prove a key
regeneration lemma.

\begin{lem}[Regeneration of optimal maps]\label{lemma9}
  Fix $\phi\in\Aut(\Gamma)$ and $X\in\Og(\Gamma)$. Let $X_\infty\in\partial_\infty\Delta_X$ be
  obtained from $X$ by collapsing a $\phi$-invariant core sub-graph $A$.
  Then, for any straight map $f_A:A\to A$ representing $\phi|_A$, and for any $\e>0$
  there is $X_\e\in\Delta_X$ such that $$\lambda_\phi(X_\e)\leq
  \max\{\lambda_\phi(X_\infty)+\e, \Lip(f_A)\}.$$

More precisely, for any $Y\in \mathbb P\Og(A)$ and map $f_Y:Y\to Y$ representing $\phi|_A$, for
any map $f:X_\infty\to X_\infty$ representing $\phi$,
for any $\widehat X\in\pi^{-1}(Y)$, and for any $\e>0$; there is
$0<\delta=\delta(f,f_Y,X_\infty,\Delta_{\widehat X})$,
such that for any  $Z\in\Delta_{\widehat X}\cap\pi^{-1}(Y)$, if $\vol_Z(Y)<\delta$  there is a
straight map $f_Z:Z\to Z$
representing $\phi$ such that $f_Z=f_Y$ on $Y$ and
$$\Lip(f_Z)\leq\max\{\lambda_\phi(X_\infty)+\e, \Lip(f_Y)\}$$ (hence the optimal map
  $\opt(f_Z)$ satisfies the same inequality\footnote{We notice that while $f_Z=f_Y$ on $Y$,
  this may no longer be true for $\opt (f_Z)$}).
\end{lem}
\proof We denote by $\sigma:X\to
X_\infty$ the map 
that collapses $A$. If $A_i$ is a component of $A$, we denote by $v_i$ the non-free vertex
$\sigma(A_i)$.  Let $k_i$ be the valence of $v_i$ in $X_\infty$.
For any $v_i$ let $E_i^1,\dots,E_i^{k_i}$ be the half-edges incident to $v_i$ in $X_\infty$.

Let $Y_i$ be the components of $Y\in\p\Og(A)$. Points in $\Delta_{\widehat X}\cap\pi^{-1}(Y)$ are built by inserting a
scaled copy of each $Y_i$  at the $v_i$ as follows. (Now we need to pass to the universal coverings.)

For every half-edge $E_i^j$ of $X_\infty$ we choose a lift in $\wt X_\infty$.
The tree $\wt{\widehat X}$ is given by attaching $\wt E_i^j$ to a point $\wt y_i^j$ of $\wt
Y_i$, and  then equivariantly attaching any other lift of the $E_i^j$. At the level of graphs
this is equivalent to choosing $y_i^j\in Y_i$. Two different choices at the level of universal
coverings differ, at the level of graphs, by  closed paths in $Y_i$ and based
at $y_i^j$. The choice of the simplex $\Delta_{\widehat X}$ fixes such ambiguity. Moreover for
any two graphs in $\pi^{-1}(Y)\cap \Delta_{\widehat X}$ the points $y_i^j$ are attached to the
same edge of $Y_i$. Let $Z\in\pi^{-1}(Y)\cap\Delta_{\widehat X}$.

Given $f_Y:Y\to Y$, consider its lift to $\wt Y$ and set $\wt z_i^j=\wt f_Y(\wt y^j_i)$. There is a
unique embedded arc $\tilde\gamma_i^j$ from $\wt z_i^j$ to $\wt y_i^j$. Let $L_i$ be the
number of edges crossed $\gamma_i^j$; so $L_i$ is an integer by including partial edges crossed in the count. $L_i$ depends only on $f_Y$ and the choices of $\wt y_i^j$,
hence it depends only on $f_Y$ and $\Delta_{\widehat X}$.

Now, given $f:X_\infty\to X_\infty$, there exists a
continuous map  $g:Z\to Z$ representing $\phi$, which agrees with $f_Y$ on
$Y$ and which is obtained by a perturbation of $f$ on edges of $X_\infty$. Namely on $E_i^j$ when we consider the path $f(E_i^j)$ we need to lift this from a path in $X_{\infty}$ to a path in $Z$ by ``filling in the missing parts". That is, we can consider $f(E_i^j)$ to be a sequence of edges in $Z$ whose completion to a path requires insertion of suitable paths which lie entirely in (the lifts of) $Y$, along with the path $\gamma_i^j$ which tells us where to ``attach" $f(E_i^j)$. 

An accurate and detailed discussion on the properties of such a map will be carried on in~\cite{partII}.

For the present purpose it is sufficient to note that there is a constant $C$ such that $g$ can
be obtained so that $\Lip(\PL(g))\leq\max\{\Lip(f)+C\vol_Z(Y),\Lip(f_Y)\}$. Moreover the constant $C$ depends only on
the $L_i$'s, the paths added to ``fill in the missing parts", and the edge-lengths of $X_\infty$. Hence it depends
only on $f_Y,\Delta_{\widehat X},X_\infty$. (We are using  $\vol_Z(Y)$ as a crude estimate for the maximum length of an edge in $Y$.)

The result follows by setting $\delta <\e/C$ and $f_Z=\PL(g)$.\qed

\medskip

\begin{defn}\label{defnojumpr}
  Fix $\phi\in\Aut(\Gamma)$. Let $X_\infty\in\partial_\infty\Delta\subset\partial_\infty\O(\Gamma)$. We say
  that $X_\infty$ has {\em not jumped in $\Delta$} if there is a sequence of points $X_i\in\Delta$
  such that $X_i \to X_{\infty}$ and $\lambda_{\phi}(X_\infty)=\lim_i\lambda_\phi(X_i)$.  We say that
  $X_\infty\in\partial_\infty\O(\Gamma)$ {\em has not jumped} if there is a simplex
  $\Delta$ intersecting $\Hor(X_\infty)$ such that $X_\infty$ has not jumped in $\Delta$.
\end{defn}
The above definition is for points in $\partial_\infty\O(\Gamma)$. By convention, we say that $X$ has not jumped for any $X\in\O(\Gamma)$.

Notice that even if $X_\infty$ has not jumped, there may exist a simplex
$\Delta$ intersecting $\Hor(X_\infty)$ such that $X_\infty$ has jumped in $\Delta$. This is because
if $A$ is the collapsed part and $\phi|_A$ does not have polynomial growth, then we can choose a point in
$\O(A)$ with arbitrarily high $\lambda_{\phi|_A}$. Moreover, even if $X_\infty$ has not jumped in
$\Delta$ it may happen that $X_\infty$ is not a continuity point of $\lambda$. For example if
the collapsed part $A$ has a sub-graph $B$ which is not invariant, then the collapse of $B$
forces $\lambda$ to increase due to Proposition~\ref{propfinvariance}, and thus we can approach $X_\infty$ with arbitrarily high $\lambda$.

Also, note that if $\lambda>1$ at some point $X$, then $\lambda$
  is in fact unbounded on $\Delta_X$. This is because if $X$ contains a loop which is not $\phi$-invariant, then by collapsing that loop we force
$\lambda$ to explode. On the other hand, if any loop is $\phi$-invariant then by
Theorem~\ref{sausagelemma} we get  $\lambda=1$.

\begin{thm}\label{newjump}
  Let $\phi\in\Aut(\Gamma)$. Let $X\in\Og(\Gamma)$ containing an invariant
  sub-graph $A$. Let $X_\infty=X/A$ and $C=\core(A)$. Then
  $$\lambda_{\phi|_C}(\Delta_C)\leq\lambda_\phi(\Delta_X).$$ Moreover the following are
  equivalent:
  \begin{enumerate}
  \item $X_\infty$ has not jumped in $\Delta_X$;
  \item $\lambda_\phi(X_\infty)\geq\lambda_{\phi}(\Delta_X)$;
  \item $\lambda_\phi(X_\infty)\geq\lambda_{\phi|_C}(\Delta_C)$.
  \end{enumerate}
 In particular, $\lambda_\phi(X_\infty)$ cannot belong to the (potentially empty)  interval
$(\lambda_{\phi|_C}(\Delta_C),\lambda_\phi(\Delta_X))$. Moreover, points realising
$\lambda_\phi(\Delta_X)$  do not jump in $\Delta_X$.
\end{thm}
\proof The first claim is a direct consequence of Proposition~\ref{oss}.
 Let $X_i\in\Delta_X$ with $X_i\to X_\infty$ without jump. Then
 $$\lambda_{\phi}(\Delta_X)\leq\lambda_\phi(X_i)\to \lambda_\phi(X_\infty)\qquad \text{ hence
 }\qquad (1)\Rightarrow (2).$$
If $\lambda_\phi(X_\infty)\geq\lambda_{\phi}(\Delta_X)$, then first claim implies $(3)$, so $(2)\Rightarrow (3)$. 

Finally, suppose $\lambda_\phi(X_\infty)\geq\lambda_{\phi|C}(\Delta_C)$. For any $\e>0$
there is $C_\e\in\Delta_C$ and a straight map $f_{C_\e}:C_\e\to C_\e$ representing $\phi|_C$ such
that $\Lip(f_{C_\e})<\lambda_{\phi|_C}(\Delta_C)+\e$. By Lemma~\ref{lemma9} there is a point
$X_\e\in X$ and a map $f_\e:X_\e\to X_\e$ representing $\phi$ such that $X_\e\to X_\infty$ as
$\e\to 0$ and $\Lip(f_\e)\leq \lambda_\phi(X_\infty)+\e$. This, plus lower semicontinuity (Theorem~\ref{fatto1}), implies 
$\lambda_\phi(X_\e)\to\lambda_\phi(X_\infty)$. So $(3)\Rightarrow (1)$.

The final statement is now a consequence of the fact that  $(2)\Rightarrow (1)$. 
\qed

\begin{lem}[Constant before jumping]\label{lconst}
  Let $\phi\in\Aut(\Gamma)$. Let $X\in\Og(\Gamma)$ be a point with a $\phi$-invariant
  sub-graph $A$. Let $X_\infty=X/A$ and let $C=\core(A)$. Let $$X_t=(1-t)X_\infty+tX$$ and let
  $C_t$ be the metric  version of $C$ in $X_t$. If $\lambda_\phi(X_\infty)<\liminf
  \lambda_\phi(X_t)$ then for small enough $t>0$, the function $\lambda_\phi(X_t)$ is locally
  constant on $t$; more precisely we have $$\lambda_\phi(X_t)=\lambda_{\phi|_C}(C_1).$$

  In particular, this is the case if displacement has jumped in $\Delta=\Delta_X$ along the segment $XX_\infty$.
\end{lem}
\proof
By Theorem~\ref{fatto1} for $t$ small enough there is an optimal map $f_t:X_t\to X_t$ whose
tension graph is contained in $C_t$.
Since $C_t$ is $\phi$-invariant, $f_t(C_t)\subset C_t$ up to homotopy. Since the vertices of
$(X_{t})_{\max}$ are at least two gated, $f((X_t)_{\max})\subset C_t$.
Therefore $\lambda_{\phi|_C}(C_t)=\Lip(f_t)$ and
$\lambda_\phi(X_t)=\Lip(f_t)=\lambda_{\phi|_C}(C_t)=\lambda_{\phi|_C}(C_1)$ (where the last
equality follows from the fact that $[C_t]=[C_1]\in\mathbb P\Og(C)$).

The last claim follows because by Theorem~\ref{newjump}, and since $X_\infty$ has jumped in
$\Delta$, we have $$\lambda_\phi(X_\infty)<\lambda_{\phi|_C}(\Delta_C)\leq\lambda_\phi(\Delta)\leq\lambda_\phi(X_t)$$
hence $\lambda_\phi(X_\infty)<\liminf_t\lambda_\phi(X_t)$.
\qed

\begin{lem}\label{lemmanewlemma}
    Let $\phi\in\Aut(\Gamma)$. Let $X\in\Delta\subseteq \Og(\Gamma)$ be a point with a
    $\phi$-invariant 
  sub-graph $A$. Let $X_\infty=X/A$ and let $C=\core(A)$. Suppose that the displacement jumps
  at $X_\infty$ along all segments of $\Delta$. Then
  \begin{enumerate}
  \item[$(I)$] $\forall Y\in\Delta_C\subseteq\O(C)$ there is $X^Y\in\Delta$ such that
    $\lambda_{\phi|_C}(Y)=\lambda_\phi(X^Y)> \lambda_\phi(X_\infty)$;
  \item[$(II)$] $\lambda_{\phi|_C}(\Delta_C)=\lambda_\phi(\Delta)\geq\lambda(X_\infty)$;
  \item[$(III)$] if $X_\infty$ does not jump\footnote{See Example~\ref{exjumpseg} for an
      explicit case of a non-jumping min-point that jumps along segments.} in $\Delta$ then it is a minpoint of $\Delta$, {\em i.e.} $\lambda_\phi(X_\infty)=\lambda_\phi(\Delta)$.
  \end{enumerate}
\end{lem}
\proof 
$(I)$. For any $Y\in\Delta_C$, let $A^Y$ be a metric version of $A$ so that $\core(A^Y)=Y$ and let 
$X_\infty^Y$ be a graph obtained by inserting a copy of $A^Y$
in the collapsed part of $X_\infty$. Since $X_\infty$ jumps along segments, by Lemma~\ref{lconst}
there is a point $X^Y$ in the segment $X_\infty^Y X_\infty$ so that
$\lambda_\phi(X^Y)=\lambda_{\phi|_C}(Y)$. Inequality $\lambda_\phi(X^Y)>
\lambda_\phi(X_\infty)$ follows from lower semicontinuity and jumping.

$(II)$. Point $(I)$ implies $\lambda_{\phi|_C}(\Delta_C)\geq\lambda_\phi(\Delta)$ and $\lambda_{\phi|_C}(\Delta_C)\geq\lambda(X_\infty)$;
Proposition~\ref{oss} gives $\lambda_{\phi|_C}(\Delta_C)\leq\lambda_\phi(\Delta)$.

$(III)$. By Theorem~\ref{newjump}, if $X_\infty$ does not jump, then $\lambda_\phi(X_\infty)\geq
\lambda_\phi(\Delta)$, and point $(II)$ concludes.
\qed

\begin{thm}\label{fatto2}
  Let $\phi\in\Aut(\Gamma)$. Let $\Delta$ be a simplex of $\Og(\Gamma)$. 
  Then there is a min-point $X_{\min}$ in $\overline{\Delta}^\infty$ (\em{i.e.} a point so that
  $\lambda_\phi(X_{\min})=\lambda_\phi(\Delta)$; note that $X_{\min}$ does not jump in $\Delta$ by
  Theorem~\ref{newjump}). 

Moreover, suppose that $X_{\min}$ is {\em maximal} in the following sense: if $X' \in \overline{\Delta}^\infty$ such that $\lambda_\phi(X') =\lambda_\phi(X_{\min})=\lambda_\phi(\Delta)$, and $\Delta_{X_{\min}} \subseteq \overline{\Delta_{X'}}^\infty$, then $\Delta_{X_{\min}} = \Delta_{X'}$. ($X_{\min}$ is maximal with respect to the partial order induced by the faces of $\Delta$). Then: 
	\begin{itemize}
        \item $\lambda_{\phi}(X_{\min}) = \lambda_{\phi}(\Delta_{X_{\min}}) =
        \lambda_{\phi}(\Delta)$;
      \item any point $P$, such that $\Delta_{X_{\min }} \subseteq \overline{\Delta_P}^\infty \subseteq
        \overline{\Delta}^\infty$, satisfies $\lambda_{\phi}(P) \geq \lambda_\phi(\Delta)$
        (hence $P$ does not
        jump in $\Delta$ by Theorem~\ref{newjump});  
        \item for any $\epsilon >0$, there exist points $Z, W$ such that:
          \begin{itemize}
		\item $Z \in \Delta$,
		\item $\Delta_{X_{\min}} \subseteq \overline{\Delta_W}^\infty \subseteq
                  \overline{\Delta}^\infty$,  
		\item $\lambda_{\phi}(W), \lambda_{\phi}(Z) \leq \lambda_{\phi}(\Delta) + \epsilon$,
		\item $\lambda_{\phi}$ is continuous along the Euclidean segments, $ZW$ and
                  $WX_{\min}$, and any point $P$ along these segments satisfies the following:
                  $\lambda_{\phi}(\Delta) \leq \lambda_{\phi}(P)$.
          \end{itemize}
		
	\end{itemize}

(We allow degeneracies, meaning that $X_{\min}$ could equal $W$, or even $Z$).
 \end{thm}
\proof 
We start by proving the existence of an $X_{\min}$.

Supose first that at every point of $\overline{\Delta}^\infty$, there is a segment in $\Delta$ to that point such that $\lambda_{\phi}$ is continuous along the segment. Then the statement is clear, since we can choose a minimizing sequence in $\Delta$, whose displacements tend to $\lambda_{\phi}(\Delta)$. This sequence has a limit point, $X_{\min}$ whose displacement is bounded above by $\lambda_{\phi}(\Delta)$ by Theorem~\ref{fatto1}. But the continuity along the segment implies that $\lambda_{\phi}(X_{\min}) \geq \lambda_{\phi}(\Delta)$, and so $\lambda_{\phi}(X_{\min}) = \lambda_{\phi}(\Delta)$. 

Thus we may assume that there is some point
$X_\infty \in\overline{\Delta}^\infty$ whose displacement jumps along all segments in
$\Delta$.
\begin{lem}\label{Lemmareferee}
For any such point $X_\infty$ there exists a point in
$\Hor(X_\infty)\cap\overline{\Delta}^\infty$ that realises $\lambda_\phi(\Delta)$. (Possibly $X_\infty$ itself realises the displacement.)
\end{lem}
\proof
We argue by induction on the rank of $\Gamma$. In rank $1$ there is nothing to prove.
Since $\lambda_\phi$ jumps at $X_\infty$ along segments, then $\lambda_\phi(X_\infty)<\infty$.
Let $X_\infty$ be obtained by some $X$ by collapsing a $\phi$-invariant sub-graph $A$, and let
$C=\core(A)$. 

By induction there is $Y\in\overline{\Delta_C}^\infty$ such that
$\lambda_{\phi|_C}(Y)=\lambda_{\phi|_C}(\Delta_C)$. Let $A^Y$ be a metric graph so that
$\core(A^Y)=Y$ (we are replacing the core part $C=\core(A)$ with $Y$). Let $X_\infty^{tY}$ be the graph obtained by inserting a volume-$t$ copy
of $A^Y$ in  the collapsed part of $X_\infty$. Note that
$X_\infty^{tY}\in\Hor(X_\infty)\cap\overline{\Delta}^\infty$.
By Proposition~\ref{oss} 
\begin{equation}
\label{eq:ji}
\lambda_\phi(X_\infty^{tY})\geq\lambda_{\phi|_C}(Y)=\lambda_{\phi|_C}(\Delta_C) 
\end{equation}
and by Lemma~\ref{lemmanewlemma}, point $(II)$,
\begin{equation}
\label{eq:jii}
\lambda_{\phi|_C}(\Delta_C)=\lambda_\phi(\Delta)\geq\lambda_\phi(X_\infty).
\end{equation}

Thus, by Theorem~\ref{newjump}, for any $t$, $X_\infty^{tY}$ does not jump in $\Delta$. If
the displacement jumps at $X_\infty$ along the segment $X^{1Y}_{\infty}X_\infty$,
then\footnote{Even if $X_\infty$ is supposed to jump along all segments in $\Delta$, it may not
jump along some segment in the boundary, and $X^{tY}_\infty$ may belong to the boundary.} by Lemma~\ref{lconst} for small
enough $t$ we have
$\lambda_\phi(X^{tY}_\infty)=\lambda_{\phi|_C}(Y)=\lambda_\phi(\Delta)$. Otherwise, by $(1)$
and $(2)$ above we have
$\lambda_\phi(X_\infty)\leq\lambda_\phi(\Delta)\leq\lambda_\phi(X^{tY}_\infty)$, and this plus
non-jumping, force $\lambda_\phi(X_\infty)=\lambda_\phi(\Delta)$.
In any case, we found a point in $\Hor(X_\infty)$ that
realises  $\lambda_\phi(\Delta)$.\qed
\medskip

Lemma~\ref{Lemmareferee} proves the first claim of Theorem~\ref{fatto2}. 
Now choose $X_{\min}$ to be maximal, as in the statement of Theorem~\ref{fatto2} (always under the assumption that $X_{\min}$ is a minimizing point), and we shall verify the list of properties. 

If $X_{\min}$ does not minimise the displacement in its simplex, then there is a point
$X'\in\Delta_{X_{\min}}$ such that
$\lambda_\phi(X')<\lambda_\phi(X_{\min})=\lambda_\phi(\Delta)$. In particular $X'$ jumps in
$\Delta$. Lemma~\ref{Lemmareferee} tells us that there is a point $X''\in \Hor(X')\cap
\overline{\Delta}^\infty$ such that $\lambda_\phi(X'')=\lambda_\phi(\Delta)$ but this would
contradict our maximality assumption.  

In general, for any $P\in\overline{\Delta}^\infty$ with $\Delta_{X_{\min}}\subseteq
\overline{\Delta_P}^\infty\subseteq\overline{\Delta}^\infty$, if $P$ jumps in $\Delta$, 
then by  Lemma~\ref{Lemmareferee} 
there is $P'\in\Hor(P)\cap\overline{\Delta}^\infty$ such that
$\lambda_\phi(P')=\lambda_\phi(\Delta)$. By our maximality assumption, this cannot happen.
So $P$ does not jump in $\Delta$ and by Theorem~\ref{newjump} we have
$\lambda_\phi(P)\geq\lambda(\Delta)$.

We check now the last property. If there is a segment in $\Delta$ to $X_{\min}$ along which $\lambda_{\phi}$ is continuous, we are done by taking $W = X_{\min}$ and $Z$ sufficiently close to $X_{\min}$. Hence we may assume that $\lambda_{\phi}$ jumps along all segments to $X_{\min}$. 

In this case, we do the construction of Lemma~\ref{Lemmareferee} with $X_{\infty} = X_{\min}$. That is, we set $W_t$ to be the graph $X^{tY}_\infty$. The maximimality assumption ensures both that $\lambda_{\phi}$ is continuous along the segment $W_1 X_{\min}$ {\em and}, that $\lambda_{\phi}$ is continuous along any segment from a point in $\Delta$ to some $W_t$. By choosing points sufficiently close to each other, we may ensure that $\lambda_{\phi}(Z), \lambda_{\phi}(W_t) \leq \lambda(\Delta) + \epsilon$.

\medskip

\qed

\begin{ex}\label{exjumpseg}[A non-jumping point which jumps along segments]
Let $F_2=\langle a,b\rangle$ and $\phi\in\Aut(F_2)$ be any iwip (so $\lambda_\phi>1$). 
For $n\geq 2$, let $F_{2n+2}=\langle a_0,b_0,a_1,b_1,\dots a_n,b_n \rangle$.
For any $i$,  $\phi$ induces  $\phi_i\in\Aut(\langle a_i,b_i\rangle)$ by
identifying $\langle a_i,b_i\rangle$ with $\langle a,b\rangle$. For any $i>0$ choose a
non-trivial $w_i\in\langle a_{i-1},b_{i-1}\rangle$ and define $\psi\in\Aut(F_{2n+2})$ by
setting $\psi|_{\langle a_0,b_0\rangle}=\phi_0$, and for $i>0$, $\psi(a_i)=\phi_i(a_i)w_i$ and
$\psi(b_i)=\phi_i(b_i)w_i$.

In order to understand the displacement of the simplex, we can just give each $a_i, b_i$ the {\em projective} length coming from the mimimum displacement of $\phi_i$, but leave ourselves free to choose the volume of the pair $\{ a_i, b_i\}$, $i \leq n-1$. This freedom allows us to impose the extra condition that $w_i$ be as short as we like, thus showing that  $\lambda_\psi \leq \lambda_\phi + \epsilon$, for any $\epsilon > 0$. Hence, $\lambda_\phi=\lambda_\psi$.

Let $R$ be the rose whose petals are labelled $a_i,b_i$ and let $\Delta=\Delta_R$. 
For any $X\in\Delta$, the displacement of $\psi$ is strictly bigger than
$\lambda_\phi$ and the minimum is attained at the graph $X_\infty$ corresponding to the collapse of
$a_i,b_i,\ i=0,\dots,n-1$, with length of petals $a_n,b_n$ given by a train track for $\phi_n$. Nonetheless, $X_\infty$ jumps along all segments; in more detail, the stretching factor for the loop $a_1$ (for example) is strictly bigger than $\lambda_\phi$ for any $X\in\Delta$. Now if we consider the segment from $X$ to $X_\infty$, the stretching factor for $a_1$ is constant along this segment, except for the discontinuity at $X_\infty$, because on this segment
the thin part is shrunk uniformly to zero.

Hence the displacement of $\psi|_{\langle
  a_0,b_0,\dots,a_{n-1},b_{n-1}\rangle}$ equals $\lambda_\phi$ and it is attained at a boundary
point. Point $W$ of Theorem~\ref{fatto2} corresponds to the collapse of petals
$a_0,b_0,\dots, a_{n-2},b_{n-2}$ from a graph $Z_W$ of $\Delta$ whose petals $a_n,b_n$ are
stretched by $\psi$ more than any other.

\end{ex}

\section{Convexity properties of the displacement function}\label{section_conv}

We recall that we are using
the terminology ``simplex'' in a wide sense, as $\Delta_X$ is a standard simplex if we work in
$\mathbb P\O(\Gamma)$ and the cone over it if we work in $\O(\Gamma)$. (Remember we use
Definition~\ref{not:gamma} for $\Gamma$.)

The displacement function $\lambda$ is scale invariant on $\O(\Gamma)$ so it descends to a function on
$\mathbb P\O(\Gamma)$. In order to control the value of $\lambda$ on segments in terms of
its value on vertices, we would like to say that $\lambda$ is convex on segments.
A minor issue appears with projectivization; if $\Delta$ is a simplex of
$\O(\Gamma)$, then its euclidean segments are well defined, and their projections on $\mathbb
P(\O(\Gamma))$ are euclidean segments in the image of $\Delta$. However, the linear
parametrization is not a projective invariant (given $X,Y$, the points $(X+Y)/2$ and $(5X+Y)/2$
are in different projective classes).

It follows that convexity of a scale invariant function is not
well-defined. In fact if $\sigma$ is a segment in $\Delta$, $\pi:\Delta\to\mathbb
P\Delta$ is the projection, and $f$ is a convex function on
$\sigma$, then $f\circ\pi^{-1}$ may be not convex. It is convex only up to reparametrization of
the segment $\pi(\sigma)$.  Such functions are called quasi-convex, and this notion will be
enough for our purposes.

\begin{defn}
  A function $f:[A,B]\to \mathbb R$ is called {\em quasi-convex} if for all $[a,b]\subseteq[A,B]$
$$\forall t\in[a,b]\qquad f(t)\leq\max\{f(a),f(b)\}.$$
\end{defn}
Note that quasi-convexity is scale invariant.

\begin{lem}\label{lconvexity}
  For any $[\phi]\in\Out(\Gamma)$ and for any open simplex $\Delta$ in $\O(\Gamma)$ the function
  $\lambda=\lambda_\phi$ is quasi-convex on segments of $\Delta$. Moreover, if $\lambda(A)>\lambda(B)$ then
  $\lambda$ is strictly monotone near $A$.
\end{lem}
\proof
Let $X$ be a $\Gamma$-graph such that $\Delta=\Delta_X$.
We use the Euclidean coordinates of $\Delta$ labelled with edges of $X$, namely a point $P$ in
$\Delta$ is given by a vector whose $e^{th}$ entry is the length of edge $e$ in $P$. In the same way,
to any reduced loop $\eta$ in $X$ we associate its occurrence vector, whose $e^{th}$ entry is
the number of times that $\eta$ passes through the edge $e$. We will denote by $\eta$ both the
loop and its occurrence vector. With this notation, the length function is bilinear:
$$L_X(\gamma)=\langle X,\gamma\rangle$$
(where $\langle,\rangle$ denotes the standard scalar product on $\mathbb R^k$.)

Let $\sigma$ be a segment in $\Delta$ with endpoints $A,B$.
Let $\gamma$ be a candidate. We consider both $\gamma$ and $\phi\gamma$ as loops in $X$.
Up to switching $A$ and $B$, we may assume that
$$\frac{\langle A,\phi\gamma\rangle}{\langle A,\gamma\rangle}\geq
\frac{\langle B,\phi\gamma\rangle}{\langle B,\gamma\rangle}.$$

Such a condition is scale invariant, and since $\lambda$ is scale invariant, up to rescaling
$B$ we may assume that $\langle B,\gamma\rangle>\langle A,\gamma\rangle$. We now parametrize
$\sigma$ in $[0,1]$:
$$\sigma(t)=A_t=Bt+(1-t)A$$
We are interested in the function:

$$F_\gamma(t)=\frac{\langle A_t,\phi\gamma\rangle}{\langle A_t,\gamma\rangle}
=\frac{\langle Bt+(1-t)A,\phi\gamma\rangle}{\langle Bt+(1-t)A,\gamma\rangle}
=\frac{\langle A,\phi\gamma\rangle+t\langle
  B-A,\phi\gamma\rangle}{\langle A,\gamma\rangle+t\langle B-A ,\gamma\rangle} $$

A direct calculation shows that the second derivative of a function of the type
$f(t)=(a+tb)/(c+td)$ is given by $2(ad-bc)d/(c+td)^3$.

So the sign of $F_\gamma''(t)$ is given by $$\big(\langle A, \phi\gamma\rangle \langle
B,\gamma\rangle-\langle B,\phi\gamma\rangle\langle A,\gamma\rangle\big)
\big(\langle B-A,\gamma\rangle\big)$$
which is non-negative by our assumption on $A,B$. Hence $F_\gamma(t)$ is (weakly)-convex and
therefore quasi-convex:
$$F_\gamma(t)\leq\max\{F_\gamma(A),F_\gamma(B)\}.$$

Now, by the Sausage Lemma~\ref{sausagelemma} we have:

\begin{eqnarray*}
  \lambda_\phi(A_t)&=&\max_\gamma F_\gamma(t)\leq\max\{\max_\gamma F_\gamma(A),\max_\gamma F_\gamma(B)\}\\&=&\max\{\lambda_\phi(A),\lambda_\phi(B)\}.
\end{eqnarray*}

Finally, since there are finitely many {\em lengths} of candidates, there is a candidate $\gamma_o$ such that for $t$ sufficiently small
we have $\lambda_\phi(A_t)=F_{\gamma_o}(t)$. By convexity, if $F_{\gamma_o}$ is not strictly monotone near $A$, then it must be locally constant, and thus
$F_{\gamma_o}''(t)=0$. Hence

$$\lambda_\phi(A)=\lambda_\phi(A_0)=\frac{\langle A,\phi\gamma_o\rangle}{\langle A,\gamma_o\rangle}=
\frac{\langle B,\phi\gamma_o\rangle}{\langle B,\gamma_o\rangle}\leq\lambda_\phi(B).$$

\qed

\begin{lem}\label{lconv2}
  Let $[\phi]\in\Out(\Gamma)$, let $\lambda=\lambda_\phi$, and let $\Delta$ be a simplex in
  $\O(\Gamma)$. Let 
  $A,B\in \overline{\Delta}^\infty$ be two points that have not jumped in $\Delta$. Then for any
  $P\in\overline{AB}$ $$\lambda(P)\leq\max\{\lambda(A),\lambda(B)\}$$

  Moreover, if $\lambda(A)\geq \lambda(B)$, then $\lambda|_{\overline{AB}}$ is continuous at $A$.
\end{lem}
\proof Let $X$ be a graph of groups so that $\Delta_X\subseteq\overline{\Delta}^\infty$
contains the interior of the segment
$\overline{AB}$. By
Lemma~\ref{lconvexity}, the function $\lambda$ is quasi-convex on the interior of $\overline{AB}$ as
a segment in $\O(X)$. Let $\{A_i\}$ and $\{B_i\}$ sequences in $\Delta$ such that
 $A=\lim A_i$ and $B=\lim B_i$ with $\lim\lambda(A_i)=\lambda(A)$ and $\lim
 \lambda(B_i)=\lambda(B)$. Such sequences exists because of the non-jumping hypothesis.
For all points $P$ in the
segment $\overline{AB}$, there is a sequence of points $P_i$ in the segment $\overline{A_iB_i}$
such that $P_i\to P$. By Lemma~\ref{lconvexity} we know
$$\lambda(P_i)\leq\max\{\lambda(A_i),\lambda(B_i)\},$$
and by lower semicontinuity
(Theorem~\ref{fatto1}) of
$\lambda$ and the non-jumping assumption, such an inequality passes to the limit. In particular,
if $\lambda(A)\geq \lambda(B)$, then $\lambda(P)\leq\lambda(A)$ for any $P\in\overline{AB}$.

Now suppose that $P^j\to A$ is a sequence in the segment $\overline{AB}$. Then by lower semicontinuity Theorem~\ref{fatto1}
applied to the space  $\O(X)$ on the segment, $\overline{AB}$, we have
$$\lambda(A)\geq\lim_j\lambda(P_j)\geq\lambda(A).$$
\qed

We end this section with an estimate of the derivative of functions like the $F_\gamma(t)$ defined as
in Lemma~\ref{lconvexity}, which will be used in the sequel. As above, we use the formalism
$\langle X,\gamma\rangle=L_X(\gamma)$.

\begin{lem}\label{derivative} Let $X$ be a $\Gamma$-graph and let
  $\Delta=\Delta_X$ be its simplex in $\O(\Gamma)$. Let $A,B\in\overline{\Delta}^\infty$. Let
  $\gamma$ be a loop in $X$ which is not collapsed neither in $A$ nor in $B$ and set  $$C=\max\{\frac{L_A(\gamma)}{L_B(\gamma)},\frac{L_B(\gamma)}{L_A(\gamma)}\}$$
   Let $\phi$ be any  automorphism of $\Gamma$. Suppose that 
$\frac{\langle B,\phi\gamma\rangle}{\langle B,\gamma\rangle} \geq \frac{\langle
  A,\phi\gamma\rangle}{\langle A,\gamma\rangle}$. Let $A_t=tB+(1-t)A$ be the linear
parametrization of the segment $AB$ in $\Delta$ and define
$F_\gamma(t)=
\frac{\langle A_t,\phi\gamma\rangle}{\langle A_t,\gamma\rangle}$.
Then 
$$0\leq F'_\gamma(t)\leq C \frac{\langle B,\phi\gamma\rangle}{\langle B,\gamma\rangle}$$

In particular, for any point $P$ in the segment $AB$ we have
$$\lambda_\phi(P)\geq\frac{\langle P,\phi\gamma\rangle}{\langle P,\gamma\rangle}\geq \frac{\langle
  B,\phi\gamma\rangle}{\langle B,\gamma\rangle}-C\lambda_\phi(B)\frac{||P-B||}{||A-B||}$$
where $||X-Y||$ denotes the standard Euclidean metric on $\Delta$. 
\end{lem}
Before the proof, a brief comment on the statement is desirable. First, note that the constant
$C$ does not depend on $\phi$. Moreover, by taking the supremum where $\gamma$ runs over all candidates given
by the Sausage Lemma~\ref{sausagelemma}, then $C$ does not even depend on $\gamma$. 
Finally if $\gamma$ is a candidate that realizes $\lambda_\phi(B)$, then we
get a bound of the steepness of $F_\gamma$ which does not depend on $\phi$ nor on $\gamma$ but
just on $\lambda_\phi(B)$ and $||A-B||$.
\proof We have
$$F_\gamma(t)=
\frac{\langle A_t,\phi\gamma\rangle}{\langle A_t,\gamma\rangle}
=\frac{\langle Bt+(1-t)A,\phi\gamma\rangle}{\langle Bt+(1-t)A,\gamma\rangle}
=\frac{\langle A,\phi\gamma\rangle+t\langle
  B-A,\phi\gamma\rangle}{\langle A,\gamma\rangle+t\langle B-A ,\gamma\rangle}$$
and a direct calculation show that
\begin{equation}
  F'_\gamma(t)=\frac{\langle B,\gamma\rangle\langle A,\gamma\rangle}{(\langle
    A_t,\gamma\rangle)^2}\left(\frac{\langle B,\phi\gamma\rangle}{\langle B,\gamma\rangle}
-\frac{\langle A,\phi\gamma\rangle}{\langle A,\gamma\rangle}\right)
\end{equation}
The first consequence of this equation is that the sign of $F'_\gamma$ does not depend on
$t$, and since 
$\frac{\langle B,\phi\gamma\rangle}{\langle B,\gamma\rangle} \geq \frac{\langle
  A,\phi\gamma\rangle}{\langle A,\gamma\rangle}$, then $F'_\gamma\geq 0$. Moreover, 
since $\langle A_t,\gamma\rangle$ is linear on $t$, 
we have $\frac{\langle B,\gamma\rangle\langle A,\gamma\rangle}{(\langle
    A_t,\gamma\rangle)^2}\leq C$. Therefore we get 
$$F'_\gamma(t)\leq C
\frac{\langle B,\phi\gamma\rangle}{\langle B,\gamma\rangle}
$$
and the first claim is proved.
For the second claim, note that the parameter $t$ is nothing but $||A-A_t||/||A-B||$ and thus
$$F_\gamma(1)-F_\gamma(t)\leq 
(1-t)C\frac{\langle B,\phi\gamma\rangle}{\langle B,\gamma\rangle}=
\frac{||B-A_t||}{||B-A||}C\frac{\langle B,\phi\gamma\rangle}{\langle B,\gamma\rangle}.
$$
If $P=A_t$, we have $F_\gamma(1)=\frac{\langle B,\phi\gamma\rangle}{\langle B,\gamma\rangle}$
and $F_\gamma(t)=\frac{\langle P,\phi\gamma\rangle}{\langle P,\gamma\rangle}$.
 By taking in account $\lambda_\phi(B)\geq\frac{\langle B,\phi\gamma\rangle}{\langle
  B,\gamma\rangle}$ and $\lambda_\phi(P)\geq\frac{\langle P,\phi\gamma\rangle}{\langle
  P,\gamma\rangle}$ we get the result.\qed

\section{Existence of minimal displaced points and train tracks at the bordification}\label{SectionTTatinfinity}

The first question that naturally arises in the study of the displacement function of
automorphisms is about the existence of min-points.
The existence of points that minimize the displacement is proved in~\cite[Theorem 8.4]{FM13}
for irreducible automorphisms. The philosophy of the proof works in the general case, but we
are forced to pass to the boundary at infinity --- whence taking in account possible jumps. The
notion of ``train track at infinity'' will be introduced in order for deal with such situations. A
second issue that appears in general case is that,  since the bordification of $\O(\Gamma)$ is not
locally compact, one cannot use compactness for claiming that minimizing sequences have
accumulation points. We overcome that difficulties first by using a Sausage
Lemma  trick as in~\cite{FM13}, and then by proving  that the set of all possible
simplex-displacements form a well-ordered subset of $\R$.

As a product of this machinery we have also other interesting results, such as the fact that the collection of partial train tracks detect all invariant free factors.

We use the terminology of Definitions\ref{def20202.6} and~\ref{not:gamma}
for $G,\G$ and $\Gamma$.

\begin{lem}[Sausage lemma trick]\label{trick}
  For any $\Gamma$, for any $X\in\overline{\O(\Gamma)}^\infty$ the set
  $\{\lambda_\phi(X):[\phi]\in\Out(\Gamma)\}$ is discrete\footnote{We include the possibility
    that $\lambda_\phi(X)=\infty$; e.g. if $X$ has a collapsed part which is not $\phi$-invariant.}. In other words, given $X$, all possible
  displacements of $X$ with respect to all automorphisms (hence markings) run over a
  discrete set (plus possibly $\infty$). 
\end{lem}
\proof
This proof is similar to that of~\cite[Theorem 8.4]{FM13}, we include it for
completeness. By the Sausage Lemma~\ref{sausagelemma}, $\lambda_\phi(X)=\Lambda(X,\phi X)$ is
computed by the ratio of translation lengths of candidates
(we include the possibility that $\lambda_\phi(X)=\infty$, e.g. if $X$ has a
collapsed part which is not $\phi$-invariant). The possible values of
$L_X(\phi\gamma)$ (with $\ul{\gamma}$ any loop in $\ul X$) form a discrete set just because $X$ has
finitely many orbits of edges. 
 Candidates are in general
infinitely many in number, but there are only finitely many lengths arising from them. Thus the possible values of
$\Lambda(X,\phi X)$ runs over a discrete subset of $\mathbb R$.\qed

\begin{thm}\label{conj}
For any $\Gamma$ the global simplex-displacement spectrum
$$\operatorname{spec}(\Gamma)= \Big\{\lambda_\phi(\Delta): [\phi]\in\Out(\Gamma), \Delta\text{
  a simplex of } \overline{\O(\Gamma)}^\infty \text{such that } \lambda_\phi(\Delta)<+\infty\}$$ is well-ordered as a subset of $\mathbb R$.
In particular, for any $[\phi]\in\Out(\Gamma)$ the spectrum of possible minimal displacements $$\operatorname{spec}(\phi)= \Big\{\lambda_\phi(\Delta):\Delta\text{ a simplex of }
\overline{\O(\Gamma)}^\infty \text{such that }\lambda_\phi(\Delta)<+\infty\}$$ is well-ordered as a subset of $\mathbb R$.
\end{thm}
\proof Recall that
we defined $\lambda_\phi(\Delta)$ as $\inf_{X\in\Delta}\lambda_\phi(X)$. For this proof we work
with co-volume one graphs (so  we are in ${\Og}_1(\Gamma)$). In any simplex we use
the standard Euclidean norm, denoted by $||\cdot||$. 

We argue by induction on the rank of $\Gamma$ (See Definition~\ref{pr4_rank}). Clearly if the rank of $\Gamma$ is one there is nothing to
prove. We now assume the claim true for any $\Gamma'$ of rank smaller than $\Gamma$.

We will show that any
monotonically decreasing sequence in $\operatorname{spec}(\Gamma)$ has a (non trivial)
sub-sequence which is constant, whence the original sequence 
is eventually constant itself. This
implies  that $\operatorname{spec}(\Gamma)$ is well-ordered. For the second claim, since
$\operatorname{spec}(\phi)$ is a subset of a well-ordered set, it is  well-ordered.

We follow the line of reasoning of~\cite[Theorem 8.4]{FM13}.  Let $\lambda_i\in\operatorname{spec}(\Gamma)$
be a monotonically decreasing sequence. Note that displacements are non-negative so $\lambda_i$
converges to some number $L\in\R$.
For any $i$ we chose $\phi_i$ and a simplex $\Delta_i$ such that
$\lambda_i=\lambda_{\phi_i}(\Delta_i)$ (those exist by definition of
$\operatorname{spec}(\Gamma)$). By Theorem~\ref{fatto2} there exists
$X_i\in\overline{\Delta_i}^\infty$ such that
$\lambda_{\phi_i}(X_i)=\lambda_{\phi_i}(\Delta_i)$ and
$\lambda_{\phi_i}(X_i)=\lambda_{\phi_i}(\Delta_{X_i})$. In particular, up to replacing $\Delta_i$
with $\Delta_{X_i}$ we may assume $\Delta_i=\Delta_{X_i}$. Since the displacement is scale
invariant, point $X_i$ can be chosen in $\overline{{\Og}_1(\Gamma)}^\infty$.
Up to possibly passing to sub-sequences we may assume
that there is $[\psi_i]\in\Out(\Gamma)$ such that $\psi_iX_i$ belongs to a fixed simplex
$\Delta$. Therefore, by replacing $\phi_i$ with $\psi_i\phi_i\psi_i^{-1}$ we may assume that
the $X_i$ all belong to the same simplex $\Delta$.\footnote{Note that $X_i$ may be a boundary point of $\O(\Gamma)$ and that
  we have made no assumption about jumps, so $X_i$ may jump.}
So we have
$$\lambda_{\phi_i}(X_i)=\lambda_{\phi_i}(\Delta_{X_i})=\lambda_i\searrow L\qquad \qquad \qquad \Delta_{X_i}=\Delta_i=\Delta.$$
 
 Up to sub-sequences, $X_i$ converges to a point  $X_\infty$ in the simplicial
 closure of $\Delta$. We show now that $\lambda_{\phi_i}(X_\infty)<+\infty$ eventually on $i$.

 If $X_\infty\in\overline{\Delta}$ (the finitary closure), then it follows from the fact that
 $\lambda_{\phi_i}(X_i)<+\infty$. Otherwise $X_\infty\in\partial_\infty \overline{\Delta}$. In
 this case, let $M$ be the length of the shortest loop in $X_\infty$. For any $\varepsilon >0$,
 $X_i$ is eventually $(M,\varepsilon)$-collapsed (because $X_\infty$ is not in a finitary face of
 $\Delta$ and $X_i\to X_\infty$). For $\varepsilon$ small enough, and eventually on $i$, the
 $\varepsilon$-thin part $(X_i)_\varepsilon$ is the core-subgraph of $X_i$ which is collapsed
 in order to reach the deformation space $\O(X_\infty)$ where $X_\infty$ lives.

 Since $\lambda_i\to L$, in particular $\lambda_i<L+1$ eventually on $i$. Thus, if
 $\varepsilon$ is small enough to satisfy the hypothesis of Proposition~\ref{propfinvariance}
 (with $C=L+1$) then $(X_i)_\varepsilon$ is $\phi_i$-invariant. This implies that
 $\lambda_{\phi_i}$ is not infinite on $\O(X_\infty)$, and in particular
 $\lambda_{\phi_i}(X_\infty)<+\infty$ as claimed.

 \medskip
 
Since $X_i$ is a  min-point for the function $\lambda_{\phi_i}$ on $\Delta_{X_i}=\Delta$, by Lemma~\ref{lconvexity} the function
 $\lambda_{\phi_i}$ either is constant on the segment $X_iX_\infty$ or it is not locally
 constant near $X_\infty$. By Lemma~\ref{lconst} in the latter case $X_\infty$ has not jumped
 w.r.t. $\lambda_{\phi_i}$ along the segment $X_iX_\infty$. 

Therefore we have the following three cases, and 
 up to subsequences we may assume that we are in the same case for any $i$:
 \begin{enumerate}
 \item $\lambda_{\phi_i}$ is constant and continuous on $X_iX_\infty$;
 \item $\lambda_{\phi_i}$ is constant on the interior of $X_iX_\infty$ and there is a jump at
   $X_\infty$, hence  $\lambda_{\phi_i}(X_\infty)<\lambda_{\phi_i}(X_i)$ by lower
   semicontinuity Theorem~\ref{fatto1};
 \item $\lambda_{\phi_i}$ is monotone increasing on the segment $X_iX_\infty$, and continuous at $X_\infty$.
 \end{enumerate}

 In the first case $\lambda_i=\lambda_{\phi_i}(X_\infty)$ is a (bounded) converging sequence of
 displacements of the single point $X_\infty$. By  Lemma~\ref{trick}, it must be eventually
 constant.
 
 In the second case we use the inductive hypothesis. Since the displacement is continuous on
 $\O(X)$, and since we have a jump at $X_\infty$, in this case $X_\infty$ is a point at
 infinite of $\O(X)$.  Let $C_i=(X_i)_\varepsilon$ the
 $\varepsilon$-thin part of $X_i$. We choose $\varepsilon$ as above so that $C_i$ is the
 $\phi_i$-invariant core sub-graph of $X_i$ which is collapsed to reach the deformation space
 $\O(X_\infty)$. Up to sub-sequences we may assume that $C_i$ is topologically the same graph for any $i$. Since $\lambda_{\phi_i}$ jumps at $X_\infty$ along the segment $X_iX_\infty$, by   
Lemma~\ref{lconst} we have $$\lambda_i=\lambda_{\phi_i}(Y)=\lambda_{\phi_i|_{C_i}}(C_i)$$
for any $Y$ in the interior of the segment $X_iX_\infty$.

If $C_i$ would not locally minimise the function $\lambda_{\phi_i|_{C_i}}$ on its simplex
$\Delta_{C_i}$, then we could perturb a little $C_i$ and strictly decrease the displacement
$\lambda_{\phi_i}(X_i)$ contradicting the fact that
$\lambda_{\phi_i}(X_i)=\lambda_{\phi_i}(\Delta_{X_i})$. So $C_i$ locally minimises the
displacement on its simplex. By quasi-convexity Lemma~\ref{lconvexity}, in any simplex local
minima are minima. Therefore $$\lambda_i=\lambda_{\phi_i|_{C_i}}(C_i)=\lambda_{\phi_i|_{C_i}}(\Delta_{C_i})\in\operatorname{spec}(C_i).$$ 
By induction $\operatorname{spec}(C_i)$ is well-ordered, hence the monotonically decreasing
sequence $\lambda_i$ must be eventually constant.

\vskip2cm

All that remains is case $(3)$. In this case $$\lambda_{\phi_i}(X_i)<\lambda_{\phi_i}(X_\infty).$$
 Let $R>0$ be such that for any face $\Delta'$ of $\Delta$ such that $X_\infty\notin
 \overline{\Delta'}^\infty$, the ball $B(X_\infty,2R)$ is disjoint from $\Delta'$. In other
 words, if $P\in B(X_\infty,2R)$ and it is obtained form $X$ by collapsing a sub-graph $P_0$, then
 $P_0$ is collapsed also in $X_\infty$. Eventually on $i$, $X_i\in B(X_\infty,R)$. Let $Y_i$ be
 the point on the Euclidean half-line from $X_\infty$ toward $X_i$, at distance exactly $R$ from $X_\infty$. 

 The stretching factor $\Lambda(X_\infty,\phi_i X_\infty)$ is realised by some
 candidates (Theorem~\ref{sausagelemma}), and  since $\lambda_{\phi_i}$ is strictly decreasing from
 $X_\infty$ to $X_i$, among such candidates there is at least one, say $\gamma_i$, whose
 stretching factor $\frac{L_X(\phi_i(\gamma_i))}{L_X(\gamma_i)}$ locally decreases near
 $X_\infty$ as a function of $X$ in the segment from $X_\infty$ to $X_i$. Moreover, the function
 $\frac{L_X(\phi_i(\gamma_i))}{L_X(\gamma_i)}$ is the ratio of two functions that are linear on
 $X$ (it is the same function
 $F_\gamma(t)$ of the proof of Lemma~\ref{lconvexity}); then its derivative on any Euclidean
 segment has constant sign, and in particular it is strictly decreasing from $X_\infty$ to $Y_i$.

By Lemma~\ref{derivative} applied with $A=Y_i$, $B=X_\infty$ (and $P=X_i$) we have 
 $$\lambda_{\phi_i}(X_i)\geq \lambda_{\phi_i}(X_\infty)\left(1-C\frac{||X_i-X_\infty||}{R}\right).$$
where
$C=\max\{\frac{L_{Y_i}(\gamma_i)}{L_{X_\infty}(\gamma_i)},\frac{L_{X_\infty}(\gamma_i)}{L_{Y_i}(\gamma_i)}\}$. Let $\varepsilon_i=||X_i-X_\infty||$.
Since there are finitely many lengths of candidates and by our choice of $R$, the
constant $C$ is uniformly bounded independently on $i$. Since $X_i\to X_\infty$ we have
$\varepsilon_i=||X_i-X_\infty||\to 0$ and thus
\begin{equation}\label{eqNref}
  \lambda_{\phi_i}(X_\infty)(1-C\varepsilon_i)\leq \lambda_{\phi_i}(X_i)\leq \lambda_{\phi_i}(X_\infty).
\end{equation}

Since $\lambda_i\searrow L$, left-hand side inequality tells us that
$\lambda_{\phi_i}(X_\infty)$ is uniformly bounded. Therefore, by Theorem~\ref{trick}, up to
sub-sequences we may assume that 
$\lambda_{\phi_i}(X_\infty)$ is a constant not depending on $i$. Moreover, $(\ref{eqNref})$ implies that
$|\lambda_{\phi_i}(X_\infty)-\lambda_i|\to 0$, hence $\lambda_{\phi_i}(X_\infty)=L$.

Now  $(\ref{eqNref})$ implies $\lambda_i\leq L$, and since $\lambda_i\searrow L$ and  it is
bounded above by $L$, it must be constant equal to $L$.\qed

\medskip
We suspect that $\operatorname{spec}(\phi)$ and $\operatorname{spec}(\Gamma)$ are not only well-ordered but in fact
discrete. However, Theorem~\ref{conj} will be enough for our purposes.

\begin{thm}[Existence of minpoints]\label{thmminptE} Let $\Gamma$ be as in Notation~\ref{not:gamma}.
  Let $[\phi]$ be any element in $\Out(\Gamma)$. Then there exists $X\in\overline{\O(\Gamma)}^\infty$
  that has not jumped and such that $$\lambda_\phi(X)=\lambda(\phi).$$
\end{thm}
\proof Let $X_i\in\O(\Gamma)$ be a minimizing sequence for $\lambda_\phi$. Without loss of
generality we may assume that the sequence $\lambda_\phi(\Delta_{X_i})$ is monotone
decreasing and Theorem~\ref{conj} implies that it is eventually constant. Therefore $X_i$ can be
chosen in a fixed simplex $\Delta$. Theorem~\ref{fatto2} concludes.\qed

\medskip

An interesting corollary of Theorem~\ref{thmminptE} is that we can characterize (global) jumps
extending equivalence ``$(3)\Leftrightarrow(1)$'' of Theorem~\ref{newjump} from a local to a global statement.

\begin{thm}\label{thmjump}
  Let $[\phi]\in\Out(\Gamma)$. Let $X\in\Og(\Gamma)$ and let
  $X_\infty\in\partial_\infty\Delta_X$ be obtained from $X$ by 
  collapsing a $\phi$-invariant core graph $A$. Then $X_\infty$ has not jumped if and only
  if $$\lambda(\phi|_A)\leq \lambda_\phi(X_\infty).$$
\end{thm}
\proof
Suppose that $X_\infty$ has not jumped (Definition~\ref{defnojumpr}). Then there
is a simplex $\Delta$ where $X_\infty$ has not jumped, and the claim follows from
Theorem~\ref{newjump} because $\lambda_{\phi|_A}\leq\lambda_{\phi|_A}(\Delta_A)$.

On the other hand, suppose $\lambda(\phi|_A)\leq\lambda_\phi(X_\infty)$. By
Theorem~\ref{thmminptE} there is a simplex in $\O(A)$ containing a minimizing sequence for $\phi|_A$.
Let $A_\e$ be an element in that simplex so that $\lambda_{\phi|_A}(A_\e)<\lambda(\phi|_A)+\e$, and let
$f_A:A_\e\to A_\e$ be an optimal map representing $\phi|_A$. Note that $A_\e$ and $A$ may be not
homeomorphic. Let $\widehat X$ be a $\Gamma$-graph obtained by inserting a copy of $A_\e$ in
$X_\infty$. (We notice that since $A_\e$ may be not homeomorphic to $A$, we can have $\Delta_{\widehat
  X}\neq\Delta_X$. We also notice that such $\Delta_{\widehat X}$ is not unique as we have plenty
of freedom of attaching the edges of $X_\infty$ to $A_\e$.)
By Lemma~\ref{lemma9}, for any $\e>0$ there is an element $X_\e\in\Delta_{\widehat X}$ and an
optimal map $f_\e:X_\e\to X_\e$ representing $\phi$ so that $X_\e\to X_\infty$ and
 $\Lip(f_\e)\leq \lambda_\phi(X_\infty)+\e$, hence $\lambda_\phi(X_\e)\leq
\lambda_\phi(X_\infty)+\e$. Thus $X_\infty$ has no jump in $\Delta_{\widehat X}$, and therefore has not jumped.\qed

\medskip

{\bf Warning: Differences between min-points at infinity and partial train tracks.} By Theorem~\ref{Theoremtt} we know that minimally displaced points and partial
train tracks coincide. But some care is needed here, as that theorem is stated for points of $\O(\Gamma)$, and not for points at infinity. In fact, given $\phi\in\Aut(\Gamma)$,
$X\in\O(\Gamma)$, and $A\subset X$ a $\phi$-invariant sub-graph, a priori
it may happen that $\lambda(\phi)$ is different if we consider $\phi$ as an element of
$\Aut(X)$ or of $\Aut(X/A)$. That is to say, we may have $X_\infty=X/A$ such that
$\lambda_\phi(X_\infty)=\lambda(\phi)$ but $X_\infty$ is not a train track point in $\O(X/A)$.

For instance, consider the case where $X=A\cup B$, with both $A$ and $B$ invariant. Suppose
that $\lambda(\phi)=\lambda(\phi|_A)>\lambda(\phi|_B)$. Now suppose that
$\lambda(\phi)=\lambda_{\phi}(X)=\lambda_{\phi|_A}(A)=\lambda_{\phi|_B}(B)$. Collapse $A$. Then the resulting
point $X_\infty$ is a min point for $\phi$ in $\overline{\O(\Gamma)}^\infty$ which has not
jumped, but since $\lambda(\phi|_B)<\lambda(\phi)$, it is not a min point for $\phi$ on
$\O(X/A)$.

We want to avoid such a pathology. Here we need to make a distinction between $\lambda(\phi)$
computed in different spaces, so we will specify the space over which we take the infimum.

\begin{lem}[Min-points vs relative min-points]\label{lemmattE}
  Let $[\phi]\in\Out(\Gamma)$. Let $X_\infty\in\overline{\O(\Gamma)}^\infty$ be such that:
  \begin{itemize}
  \item There is $X\in\O(\Gamma)$ such that $\ul{X_\infty}$ is obtained from $\ul{X}$ by collapsing a
    (possibly empty) core sub-graph $\ul{A}$ in $\ul{X}$,  and such that $X_\infty$ has not jumped in $\Delta_X$;
  \item $\lambda_\phi(X_\infty)=\inf_{Y\in\O(\Gamma)}\lambda_\phi(Y)$.
  \end{itemize}
Suppose moreover that $X_\infty$  maximizes the dimension of $\Delta_{X_\infty}$ among the set
of elements in $\overline{\O(\Gamma)}^\infty$ satisfying such conditions (such a set is not empty
by Theorem~\ref{thmminptE}). Then
  $$\lambda_\phi(X_\infty)=\inf_{Y\in\O(X/A)}\lambda_\phi(Y).$$
(Hence it is in $\TT(\phi)\subset\O(X/A)=\O(X_\infty)$.)
\end{lem}
\proof If $A$ is empty this is an instance of Theorem~\ref{Theoremtt}. In general, again by
Theorem~\ref{Theoremtt}, it suffices to show that $X_\infty$ is a partial train track point for $\phi$
in $\O(X_\infty)$.  
Suppose the contrary. In this case, near $X_\infty$
there is a point $X'_\infty\in\O(X_\infty)$ such that $\lambda_\phi(X_\infty')<\lambda_\phi(X_\infty)$.
Indeed, if $X_\infty$ is not a local
min point in $\Delta_{X_\infty}\subset\O(X_\infty)$, then we can find $X'_\infty$ just near
$X_\infty$ in $\Delta_{X_\infty}$. Otherwise, by Lemma~\ref{LemmaX}
there is a point $X'_\infty$ obtained form $X_\infty$ by folds directed by optimal maps (and such
that $\dim(\Delta_{X'_\infty})>\dim(\Delta_{X_\infty})$)  such that
$\lambda_\phi(X_\infty')<\lambda_\phi(X_\infty)$. 

Let $\e=(\lambda_\phi(X_\infty)-\lambda_\phi(X_\infty'))/2$.

Since $X_\infty$ has not jumped, by Theorem~\ref{thmjump} we have
$\lambda(\phi|_A)\leq \lambda_\phi(X_\infty)$. If $\lambda(\phi|_A)< \lambda_\phi(X_\infty)$,
let $A'\in\O(A)$ be a point such that $\lambda_{\phi|_A}(A')<\lambda_\phi(X_\infty)$. Now
Lemma~\ref{lemma9} provides an element of $\O(\Gamma)$ which is displaced less or equal than
$\max\{\lambda_{\phi|_A}(A'),\lambda_\phi(X_\infty')+\e\}$, contradicting the fact that
$X_\infty$ is a minpoint for $\lambda_\phi$ in $\O(\Gamma)$. Therefore $\lambda(\phi|_A)=\lambda(X_\infty)$.

By Theorem~\ref{thmminptE} there is
$A_\infty\in\overline{\O(A)}^\infty$ such that
$\lambda_{\phi|_A}(A_\infty)=\lambda(\phi|_A)$ and which has not jumped in $\O(A)$.
Thus $A_\infty$ is obtained, without jumps, from a point
$A'_\infty\in\O(A)$ by collapsing a (possibly empty) invariant core sub-graph $B$. So $A_\infty\in\O(A'_\infty/B)$.

Let $Y$ be a $\Gamma$-graph obtained by inserting a copy of $A_\infty'$ in $X'_\infty$. Let
$Y'$ be the graph obtained collapsing $B$. $Y'$ belongs to the simplicial boundary of
$\Delta_Y$ and, since $A_\infty$ has no jump, then so does $Y'$. Now, observe that
$Y'\in\O(Y/B)$ and $A_\infty$ is a $\phi$-invariant subgraph of $Y'$ so that $Y'/A_\infty=X_\infty'$.
Lemma~\ref{lemma9}
provides an element in $Y'_\infty\in\O(Y/B)$, in the same simplex of $Y'$, which is displaced no more than
$\lambda_{\phi|_A}(A_\infty)$ (because
$\lambda_\phi(X_\infty')<\lambda_\phi(X_\infty)=\lambda_{\phi|_A}(A_\infty)$).
Now, $Y'_\infty$ is a new minpoint for $\lambda_\phi$ with
$\dim(\Delta_{Y'_\infty})>\dim(\Delta_{X_\infty})$ contradicting the maximality hypothesis on
$X_\infty$. It follows that $X_\infty$ is a train track point in $\O(X_\infty)$ as desired.\qed

\medskip

We have just seen that, even if non-jumping min-points are not necessarily partial train
tracks, some of them are. Conversely, we see now non-jumping partial train tracks at the bordification are always min-points for $\lambda_\phi$.

\begin{lem}\label{triangin}
Let $\phi\in\Aut(\Gamma)$ and let $X\in\O(\Gamma)$. If there is $k$ so that there is a constant
$A>0$ such that for any $n>>1$
$$Ak^n\leq \Lambda(X,\phi^n X)$$ then $k\leq \lambda(\phi)$.
\end{lem}
\proof This follows from the multiplicative triangular inequality.
For any $Y\in\O(\Gamma)$ we have $\Lambda(Y,\phi^nY)\leq \Lambda(Y,\phi Y)^n$. Define a constant
$C=\Lambda(X,Y)\Lambda(Y,X)$ and notice that we also have $C=\Lambda(X,Y)\Lambda(\phi
Y,\phi X)$. Then,
$$Ak^n\leq\Lambda(X,\phi^n X)\leq\Lambda(X,Y)\Lambda(Y,\phi^n Y)
\Lambda(\phi^nY,\phi^n X)\leq C\Lambda(Y,\phi Y)^n$$
whence, for any $n$ $$\left(\frac{k}{\Lambda(Y,\phi Y)}\right)^n\leq \frac{C}{A}.$$
This implies $k\leq \Lambda(Y,\phi(Y))$. By choosing a minimizing sequence of points
$Y_i$ we get $k\leq \lambda(\phi)$.\qed

\begin{lem}\label{dom5}
  Let $\phi\in\Aut(\Gamma)$. Let $X_\infty\in\overline{\Og(\Gamma)}$ which has not jumped.
  Suppose that there is a loop $\gamma \in X_\infty$ and $k>0$ such that
  $L_{X_\infty}(\phi^n)(\gamma)\geq k^n L_{X_\infty}(\gamma)$. Then
  $$k\leq\lambda(\phi).$$

  In particular, if $X_\infty$ is a partial train track for $\phi$ as an element of $\Aut(X_\infty)$,
  then it is a min-point  for $\phi$ as an element of $\Aut(\Gamma)$.
\end{lem}
\proof Let $X\in\Og(\Gamma)$ so that $X_\infty$ is obtained from $X$ by collapsing a core sub-graph
$A\subset X$. Let $X_\e$ be a point of $X$ where $\vol(A)<\e$.
Let $\gamma$ be as in the hypothesis.
For $\e$ small enough we have
$L_{X_\e}(\gamma)\leq 10 L_{X_\infty}(\gamma)$, and therefore
\begin{eqnarray*}
\Lambda(X_\e,\phi^nX_\e) \geq  \frac{L_{X_\e}(\phi^n\gamma)}{L_{X_\e}(\gamma)}\geq
\frac{L_{X_\infty}(\phi^n\gamma)}{10L_{X_\infty}(\gamma)}\geq
\frac{k^nL_{X_\infty}(\gamma)}{10L_{X_\infty}(\gamma)}=
\frac{k^n}{10}.
\end{eqnarray*}
By Lemma~\ref{triangin} we have $\lambda(\phi)\geq k$.

For the second claim it suffice to choose let $\gamma$ a legal candidate that
realizes $\Lambda(X_\infty,\phi X_\infty)$. So
$L_{X_\infty}(\phi^n(\gamma))=\lambda_\phi(X_\infty)^nL_{X_\infty}(\gamma)$.

Hence $\lambda(\phi)\geq \lambda_\phi(X_\infty)$ and since $X_\infty$
has not jumped $\lambda(\phi)\leq\lambda_\phi(X_\infty)$.\qed

We are now in position to complete extension of Theorem~\ref{newjump} from a local to a
global statement (see~Theorem~\ref{thmjump}).

\begin{thm}\label{corlalx}
  Let $\phi\in\Aut(\Gamma)$. Let $X\in\O(\Gamma)$ and $X_\infty$ be such that $\ul{X_\infty}$ is obtained from $\ul{X}$ by
  collapsing a $\phi$-invariant core sub-graph $\ul{A}$. Then $$\lambda(\phi|_A)\leq \lambda(\phi).$$

Moreover, if $\lambda(\phi|_A)=\lambda_\phi(X_\infty)$, then $$\lambda(\phi)=\lambda(\phi|_A).$$

In particular $X_\infty$ has not jumped
if and only if $$\lambda(\phi)\leq \lambda(X_\infty).$$
\end{thm}
\proof Let $\lambda=\lambda(\phi|_A)$.
By  Lemma~\ref{lemmattE} and Theorem~\ref{Theoremtt}, there is $\widehat A\in\overline{\O(A)}^\infty$ which is a min-point for
$\phi|_A$, which has not  jumped in $\O(A)$, and which is a partial train track for $\phi|_A$ as an
element of $\Aut(\widehat A)$. Let $f_A$ be a partial
train track map $f_A:\widehat A\to \widehat A$ representing $\phi|_A$.
Therefore, there is a periodic line $\gamma$ in $\widehat A_{\max}$ with legal images in $\widehat A_{\max}$ and
stretched exactly by $\lambda$. Let now $\widehat X\in\O(\Gamma)$ be obtained by
inserting a copy of $\widehat A$ in $X_\infty$. Since $\widehat A$ has not jumped in $\O(A)$, then
$\widehat X$ has not jumped in $\O(\Gamma)$.

Let $f:\widehat X\to\widehat X$ be any straight map
representing $\phi$ so that $f|_A=f_A$. Therefore $f_A^n(\gamma)$ is immersed for any $n$ and the length of $f_A^n(\gamma)$ is $\lambda^n$
times the length of $\gamma$. It follows that $L_{\widehat X}((\phi^n)\gamma)=\lambda^n(L_{\widehat
  X}(\gamma))$.

By Lemma~\ref{dom5}
$\lambda(\phi|_A)=\lambda\leq \lambda(\phi)$, and the first claim is proved. Moreover, if
$\lambda(\phi|_A)=\lambda_\phi(X_\infty)$, then
$$\lambda(\phi)\leq\lambda(X_\infty)=\lambda(\phi|_A)=\lambda\leq\lambda(\phi)$$ and therefore
all inequalities are equalities.
Finally, if $X$ has not jumped then $\lambda(X)\geq \lambda(\phi)$ just because this inequality is true by definition for points in
$\O(\Gamma)$ and clearly passes to limits of non-jumping sequences. The converse inequality
follows from the second claim and Theorem~\ref{thmjump}.
\qed

\medskip
Note that Theorem~\ref{corlalx} implies that {\em a posteriori} we can remove the non-jumping
requirement from Theorem~\ref{thmminptE} and Lemma~\ref{lemmattE}.

\begin{cor}[Min-points don't jump]\label{corminpoindonjump}
  Let $\phi$ be any element in $\Aut(\Gamma)$. If $X\in\overline{\O(\Gamma)}^\infty$
  is such that $\lambda_\phi(X)=\lambda(\phi)$, then it has not jumped.
\end{cor}
\proof This is a direct consequence of Theorem~\ref{corlalx}.\qed

\medskip

We introduce the notion of partial train track at infinity.

\begin{defn}[Partial Train track at infinity]\label{dttinfty}
  Let $\phi\in\Aut(\Gamma)$. The set $\TT^\infty(\phi)$ is defined as the set of points
  $X\in\overline{\O(\Gamma)}^\infty$ such that  $X$ has not jumped, and $X$ is a partial train track
  point for $\phi$ in $\O(X)$. (Hence $\lambda_\phi(X)=\lambda(\phi)$ by Lemma~\ref{dom5}.)
\end{defn}

Note that $\TT(\phi)\subset\TT^{\infty}(\phi)$. The main differences are that $\TT(\phi)$ may
be empty, while any $\phi$ has a partial train track in $\TT^\infty(\phi)$. On the
other hand, $\TT(\phi)$ coincides with the set of minimally displaced points, while
$\TT^\infty(\phi)$ may be strictly contained in the set of minimally displaced points.

With this definition  we can collect some of the above results
in the following simple statement, which is a straightforward
consequence of Theorems~\ref{Theoremtt},~\ref{thmminptE} and Lemmas~\ref{lemmattE},~\ref{dom5}.
\begin{thm}[Existence of partial train tracks at infinity]\label{corttE}
  For any $[\phi] \in\Out(\Gamma)$, $TT^\infty(\phi)\neq\emptyset$. For any
  $X\in\TT^{\infty}(\phi)$, $\lambda_\phi(X)=\lambda(\phi)$.
\end{thm}

We end this section by discussing some interesting consequences of the theory developed so far. 
In particular we show that if $\phi$ is reducible then there is a train track showing reducibility.

\begin{thm}[Detecting reducibility]\label{corred}
  Let $\phi\in\Aut(\Gamma)$ be reducible. Then there is $T\in\TT^\infty(\phi)$ such that either
  $T\in\partial_\infty\O(\Gamma)$ or there is an optimal
  map $f_T:T\to T$ representing $\phi$ such that there is a proper sub-graph of $T$ which is $f_T$-invariant.
\end{thm}
\proof Since $\phi$ is reducible there is  $X\in\O(\Gamma)$, a straight map $f:X\to X$
representing $\phi$ and a proper non trivial sub-graph $A\subset \ul{X}$ such that $f(A)=A$. We can therefore
collapse $A$ and $\lambda$ won't explode. By Theorem~\ref{corttE} there is a partial train
track point $Z$ for $\phi$ in $\overline{\O(X/A)}^\infty$ and a partial train track point $Y$
for $\phi|_A$ in $\overline{\O(A)}^\infty$.  If
$\lambda_{\phi|_A}(Y)\leq\lambda_\phi(Z)$, then
$Z\in\TT^\infty(\phi)\cap\partial_\infty\O(\Gamma)$ and  we are done.
Otherwise, since $Z$ has not jumped (as a point of $\partial_\infty\O(X/A)$),
we can regenerate it to a point $Z'\in\O(X/A)$ with
$\lambda_\phi(Z')<\lambda_{\phi|_A}(Y)$. We now apply regeneration Lemma~\ref{lemma9} to $Y$ and
$Z'$. If $Y\in\partial_\infty\O(A)$, then we get a partial train track for $\phi$ in
$\partial_\infty\O(\Gamma)$. If $Y\in\O(A)$ we get a partial train track for $\phi$ in $\O(\Gamma)$
admitting $Y$ as an invariant sub-graph.\qed

\medskip

In fact, the proof of Theorem~\ref{corred} proves more: that the set of partial train tracks detect any
(maximal) invariant free factor system. Roughly, if $A$ is an invariant free factor syetem, then there is
$B\supseteq A$ (possibly $B=A$) and a partial train track point which shows $B$ as the
fundamental group of an invariant sub-graph. The precise statement is the following.

\begin{thm}[Strong reformulation of Corollary~\ref{corred}]\label{strongcorred}
  Let $\phi\in\Aut(\Gamma)$. Let $X$ be a $\Gamma$-graph having a $\phi$-invariant core
  sub-graph $A$. Then there is $Z\in\overline{\O(X/A)}^\infty$ and $W\in\Hor_{\O(\Gamma)}(Z)$
  such that the simplex $\Delta_W$ 
  contains a minimizing sequence for $\lambda$. Moreover if $Y$ is the graph used to
  regenerate $W$ from $Z$, then the minimizing sequence can be chosen with straight maps $f_i$
  such that $f_i(Y)=Y$ and $\Lip(f_i)\to\lambda(\phi)$.
\end{thm}
\proof Follows from the proof of Theorem~\ref{corred} (and Lemma~\ref{lemma9}).\qed

\medskip

Finally, as in the case of irreducible automorphisms, the existence of partial train tracks
gives the following fact.
\begin{cor}\label{corollaryphin}
  For any $\phi\in\Aut(\Gamma)$ we have $\lambda(\phi^n)=\lambda(\phi)^n.$
\end{cor}
\proof This follows from Theorem~\ref{corttE} and Lemma~\ref{Lemmatt2}.\qed

\providecommand{\bysame}{\leavevmode\hbox to3em{\hrulefill}\thinspace}
\providecommand{\MR}{\relax\ifhmode\unskip\space\fi MR }
\providecommand{\MRhref}[2]{%
  \href{http://www.ams.org/mathscinet-getitem?mr=#1}{#2}
}
\providecommand{\href}[2]{#2}

\end{document}